\documentclass[journal,11pt,draftclsnofoot,onecolumn,twoside]{IEEEtran}
\usepackage{cite}
\ifCLASSINFOpdf
   \usepackage[pdftex]{graphicx}
\else
   \usepackage[dvips]{graphicx}
\fi
\usepackage[cmex10]{amsmath}
\usepackage{amsfonts}
\usepackage{array}

\usepackage[caption=false]{caption}
\usepackage[font=footnotesize]{subfig}

\newcommand{\define}{\stackrel{\text{\tiny def}}{=}}

\newcommand{\finesys}{A}
\newcommand{\finegrid}{\Omega}

\newcommand{\rcoarsesys}{\bar{A}}
\newcommand{\rcoarsegrid}{\bar{\Omega}}
\newcommand{\rcgc}{\bar{K}}
\newcommand{\rdown}{\bar{D}}
\newcommand{\rup}{\bar{U}}
\newcommand{\rhap}{\bar{N}}
\newcommand{\rinterp}{\bar{I}_I}
\newcommand{\rrestri}{\bar{I}_R}
\newcommand{\rifilter}{\bar{F}_I}
\newcommand{\rrfilter}{\bar{F}_R}
\newcommand{\rifilterev}{\bar{\Pi}_I}
\newcommand{\rifilterevL}{\bar{\Pi}_{I,L}}
\newcommand{\rifilterevH}{\bar{\Pi}_{I,H}}
\newcommand{\rrfilterev}{\bar{\Pi}_R}
\newcommand{\rrfilterevL}{\bar{\Pi}_{R,L}}
\newcommand{\rrfilterevH}{\bar{\Pi}_{R,H}}
\newcommand{\rdelta}{\bar{\Delta}}
\newcommand{\rmodal}{\bar{\Gamma}}

\newcommand{\bcoarsesys}{\widetilde{A}}
\newcommand{\bcoarsegrid}{\widetilde{\Omega}}
\newcommand{\bcgc}{\widetilde{K}}
\newcommand{\bdown}{\widetilde{D}}
\newcommand{\bup}{\widetilde{U}}
\newcommand{\bhap}{\widetilde{N}}
\newcommand{\binterp}{\widetilde{I}_I}
\newcommand{\brestri}{\widetilde{I}_R}
\newcommand{\bifilter}{\widetilde{F}_I}
\newcommand{\brfilter}{\widetilde{F}_R}
\newcommand{\bifilterev}{\widetilde{\Pi}_I}
\newcommand{\bifilterevL}{\widetilde{\Pi}_{I,L}}
\newcommand{\bifilterevH}{\widetilde{\Pi}_{I,H}}
\newcommand{\brfilterev}{\widetilde{\Pi}_R}
\newcommand{\brfilterevL}{\widetilde{\Pi}_{R,L}}
\newcommand{\brfilterevH}{\widetilde{\Pi}_{R,H}}
\newcommand{\bdelta}{\widetilde{\Delta}}
\newcommand{\bmodal}{\widetilde{\Gamma}}

\newtheorem{thm}{Theorem}
\newtheorem{pro}{Proposition}
\newtheorem{lem}{Lemma}
\newtheorem{cor}{Corollary}
\newtheorem{defn}{Definition}
\newtheorem{rmk}{Remark}

\begin{document}
%
% paper title
% can use linebreaks \\ within to get better formatting as desired
\title{Direct Multi--grid Methods for Linear Systems\\ with Harmonic Aliasing Patterns}
%
%
% author names and IEEE memberships
% note positions of commas and nonbreaking spaces ( ~ ) LaTeX will not break
% a structure at a ~ so this keeps an author's name from being broken across
% two lines.
% use \thanks{} to gain access to the first footnote area
% a separate \thanks must be used for each paragraph as LaTeX2e's \thanks
% was not built to handle multiple paragraphs
%

\author{Pablo~Navarrete~Michelini% <-this % stops a space
\thanks{P. Navarrete Michelini is with the Department
of Electrical Engineering, Universidad de Chile, Santiago,
RM, 8370451 Chile.}% <-this % stops a space
%\thanks{Manuscript received April 19, 2005; revised January 11, 2007.}}
\thanks{Manuscript received September X, 2009.}}

% note the % following the last \IEEEmembership and also \thanks - 
% these prevent an unwanted space from occurring between the last author name
% and the end of the author line. i.e., if you had this:
% 
% \author{....lastname \thanks{...} \thanks{...} }
%                     ^------------^------------^----Do not want these spaces!
%
% a space would be appended to the last name and could cause every name on that
% line to be shifted left slightly. This is one of those "LaTeX things". For
% instance, "\textbf{A} \textbf{B}" will typeset as "A B" not "AB". To get
% "AB" then you have to do: "\textbf{A}\textbf{B}"
% \thanks is no different in this regard, so shield the last } of each \thanks
% that ends a line with a % and do not let a space in before the next \thanks.
% Spaces after \IEEEmembership other than the last one are OK (and needed) as
% you are supposed to have spaces between the names. For what it is worth,
% this is a minor point as most people would not even notice if the said evil
% space somehow managed to creep in.

% The paper headers
%\markboth{{IEEE} Trans. Signal Process.,~Vol.~V, No.~N, September~2009}%
\markboth{{IEEE} Transactions On Signal Processing,~preprint, September~2009}%
{Navarrete Michelini: Direct Multi--grid Methods}
% The only time the second header will appear is for the odd numbered pages
% after the title page when using the twoside option.
% 
% *** Note that you probably will NOT want to include the author's ***
% *** name in the headers of peer review papers.                   ***
% You can use \ifCLASSOPTIONpeerreview for conditional compilation here if
% you desire.

% If you want to put a publisher's ID mark on the page you can do it like
% this:
%\IEEEpubid{0000--0000/00\$00.00~\copyright~2007 IEEE}
% Remember, if you use this you must call \IEEEpubidadjcol in the second
% column for its text to clear the IEEEpubid mark.

% use for special paper notices
%\IEEEspecialpapernotice{(Invited Paper)}

% make the title area
\maketitle

\begin{abstract}
% 200 words maximum
Multi--level numerical methods that obtain the exact solution of a linear system are presented. The methods are devised by combining ideas from the full multi--grid algorithm and perfect reconstruction filters. The problem is stated as whether a direct solver is possible in a full multi--grid scheme by avoiding smoothing iterations and using different coarse grids at each step. The coarse grids must form a partition of the fine grid and thus establishes a strong connection with domain decomposition methods. An important analogy is established between the conditions for direct solution in multi--grid solvers and perfect reconstruction in filter banks. Furthermore, simple solutions of these conditions for direct multi--grid solvers are found by using mirror filters. As a result, different configurations of direct multi--grid solvers are obtained and studied.
\end{abstract}
% IEEEtran.cls defaults to using nonbold math in the Abstract.
% This preserves the distinction between vectors and scalars. However,
% if the journal you are submitting to favors bold math in the abstract,
% then you can use LaTeX's standard command \boldmath at the very start
% of the abstract to achieve this. Many IEEE journals frown on math
% in the abstract anyway.

% Note that keywords are not normally used for peerreview papers.

% For peer review papers, you can put extra information on the cover
% page as needed:
% \ifCLASSOPTIONpeerreview
% \begin{center} \bfseries EDICS Category: 3-BBND \end{center}
% \fi
%
% For peerreview papers, this IEEEtran command inserts a page break and
% creates the second title. It will be ignored for other modes.
\IEEEpeerreviewmaketitle

%%%%%%%%%%%%%%%%%%%%%%%%%%%%%%%%%%%%%%%%%%%%%%%%%%%%
\section{Introduction}
\label{sec:introduction}
% The very first letter is a 2 line initial drop letter followed
% by the rest of the first word in caps.
% 
% form to use if the first word consists of a single letter:
% \IEEEPARstart{A}{demo} file is ....
% 
% form to use if you need the single drop letter followed by
% normal text (unknown if ever used by IEEE):
% \IEEEPARstart{A}{}demo file is ....
% 
% Some journals put the first two words in caps:
% \IEEEPARstart{T}{his demo} file is ....
% 
% Here we have the typical use of a "T" for an initial drop letter
% and "HIS" in caps to complete the first word.
\IEEEPARstart{T}{his} study focuses on the problem of solving the linear system of equations
% You must have at least 2 lines in the paragraph with the drop letter
% (should never be an issue)
\begin{equation} \label{eq:fine_system}
\finesys u=f
\end{equation}
over the field of complex numbers. The problem is restricted to the case when the number of equations, $n$, is the same as the number of unknowns, and the number $n$ is even (in some cases a power of $2$). The \emph{system matrix} $A\in\mathbb{C}^{n\times n}$ is sparse and it will be assumed to be invertible, with special attention to ill--conditioned cases. The problem becomes challenging when $n$ scales to large numbers (e.g. thousands of unknowns). This situation arises frequently in scientific and engineering computations, most notably in the solution of PDEs \cite{YSaad_1996a,WHackbusch_1993a} and other areas like simulation of stochastic models \cite{JRNorris_1998a} and the solution of optimization problems \cite{JEDennis_RBSchnabel_1983a}.

A vast amount of numerical methods exist to solve this problem efficiently. They vary from \emph{direct} (or \emph{exact}) solvers that compute the exact solution, $u$, to \emph{iterative} solvers that compute a sequence of approximations that converges to the exact solution, $v_k \rightarrow u$. Here, convergence must be defined to minimize some norm of the approximation error, $e_k=u-v_k$. Direct solvers are often based on some type of matrix factorization, being the most popular those based on $LU$ decomposition \cite{ISDuff_AMErisman_JKReid_1986a,JJDongarra_ISDuff_DCSorensen_HAvan_der_Vorst_1998a}. Among the iterative solvers, the most common are: stationary iterative methods (e.g. Gauss--Seidel, Jacobi and Richardson iterations), Krylov subspace methods (e.g. conjugate gradients, GMRES and BiCG) and multi--level methods (e.g. multi--grid and domain decomposition) \cite{YSaad_1996a,WHackbusch_1993a}.

In this study, multi--level numerical methods working as direct solvers are obtained. In the same category there are other direct multi--level solvers like: total reduction methods \cite{JSchroder_UTrottenberg_1973a,JSchroder_UTrottenberg_HReutersberg_1976a}, partial (cyclic) reduction methods \cite{BLBuzbee_GHGolub_CWNielson_1970a} and $LU$ factorization of non--standard forms \cite{DGinesa_GBeylkin_JDunnc_1998a}. Besides their structural differences, each method works under certain limitations. Total reduction and partial (cyclic) reduction methods are specifically designed for Poisson's equation, and $LU$ factorization of non--standard forms works for elliptic problems. In this study, the limitations are not described in terms of categories of PDEs but in terms of two additional properties on the system. These are: $\finesys$ has to be diagonalizable, and ignoring some of the unknowns should produce a specific aliasing pattern.

First, by assuming that the system matrix is diagonalizable, we have the eigendecomposition
\begin{equation}
\finesys = W \Lambda V^H \;,
\end{equation}
where the columns of $W$ form the set of right eigenvectors, $\Lambda$ is a diagonal matrix with the eigenvalues of $A$, and the columns of $V$ form the set of left eigenvectors. Here, the right and left eigenvectors form a \emph{biorthogonal basis} so that $V^H W=I$ and they do not need to be equal, which allows the case of a non--symmetric system matrix. This restriction limits the applications to non--defective problems.

% An example of a double column floating figure using two subfigures.
% (The subfig.sty package must be loaded for this to work.)
% The subfigure \label commands are set within each subfloat command, the
% \label for the overall figure must come after \caption.
% \hfil must be used as a separator to get equal spacing.
% The subfigure.sty package works much the same way, except \subfigure is
% used instead of \subfloat.
%
%\begin{figure*}[!t]
%\centerline{\subfloat[Case I]\includegraphics[width=2.5in]{subfigcase1}%
%\label{fig_first_case}}
%\hfil
%\subfloat[Case II]{\includegraphics[width=2.5in]{subfigcase2}%
%\label{fig_second_case}}}
%\caption{Simulation results}
%\label{fig_sim}
%\end{figure*}
%
% Note that often IEEE papers with subfigures do not employ subfigure
% captions (using the optional argument to \subfloat), but instead will
% reference/describe all of them (a), (b), etc., within the main caption.
\begin{figure*}[!t]
\centerline{
\subfloat[$w_1$ and down--sampled $w_1$]{\includegraphics[width=0.122\textwidth]{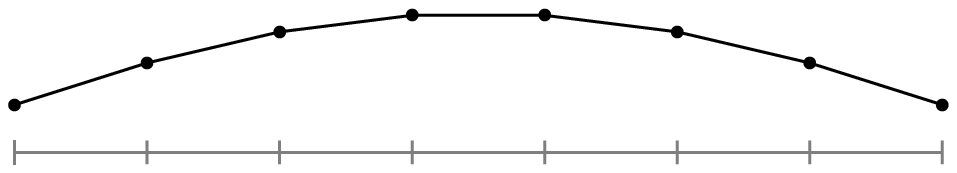} \includegraphics[width=0.122\textwidth]{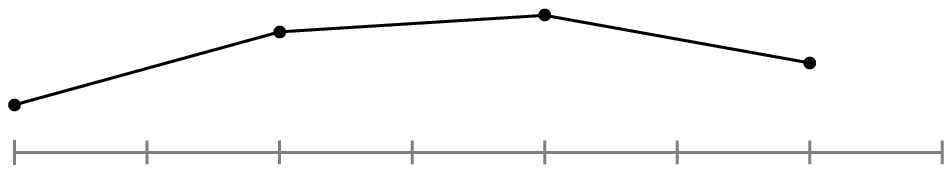}} \hfil
\subfloat[$w_2$ and down--sampled $w_2$]{\includegraphics[width=0.122\textwidth]{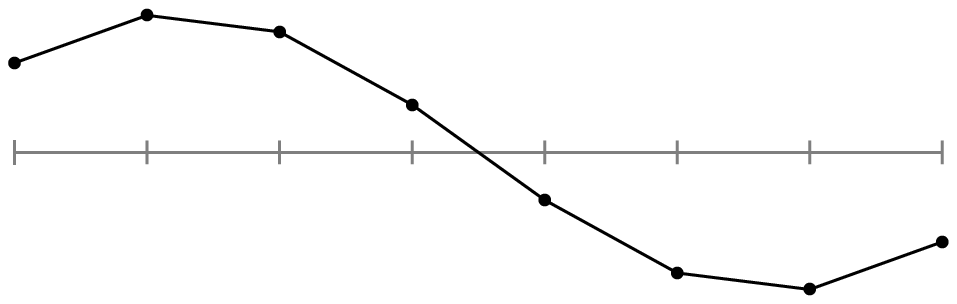} \includegraphics[width=0.122\textwidth]{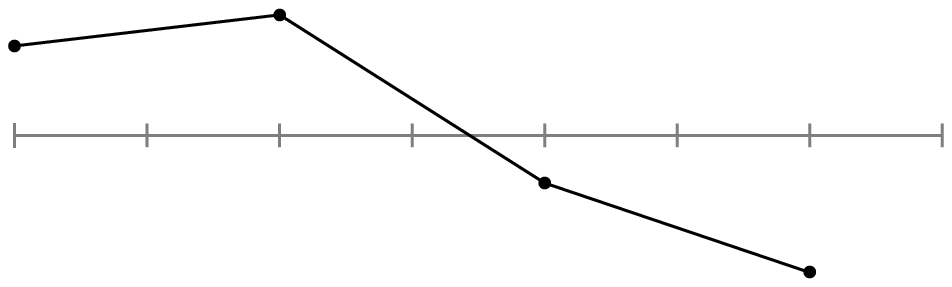}} \hfil
\subfloat[$w_3$ and down--sampled $w_3$]{\includegraphics[width=0.122\textwidth]{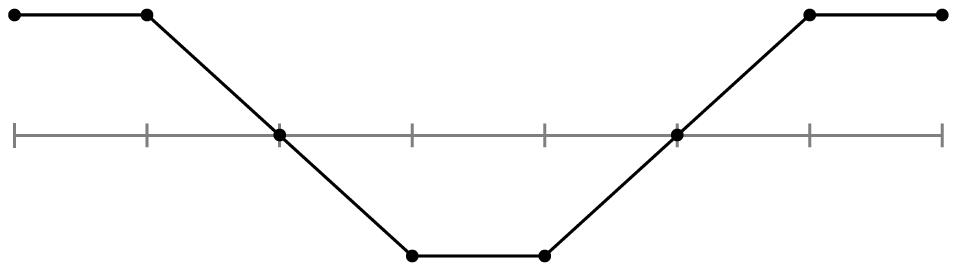} \includegraphics[width=0.122\textwidth]{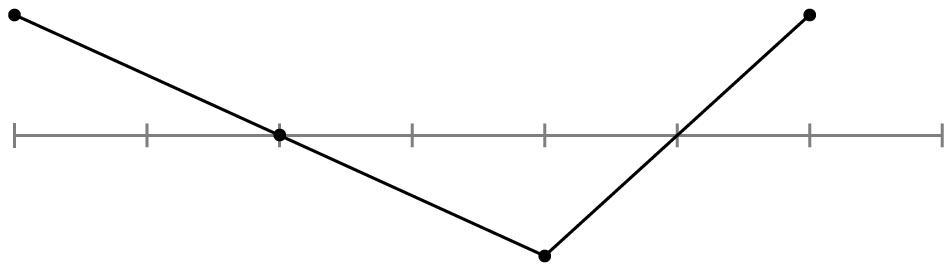}} \hfil
\subfloat[$w_4$ and down--sampled $w_4$]{\includegraphics[width=0.122\textwidth]{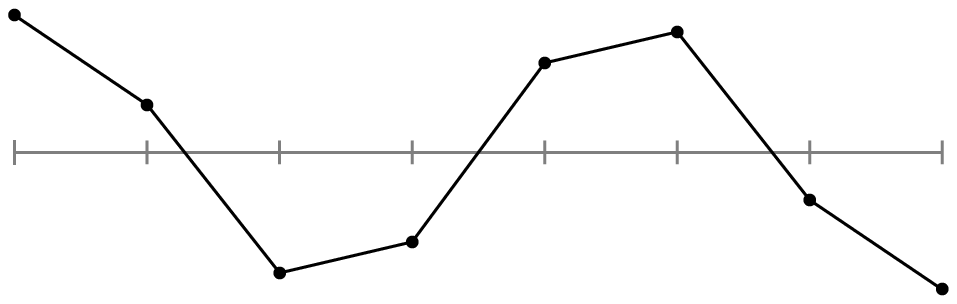} \includegraphics[width=0.122\textwidth]{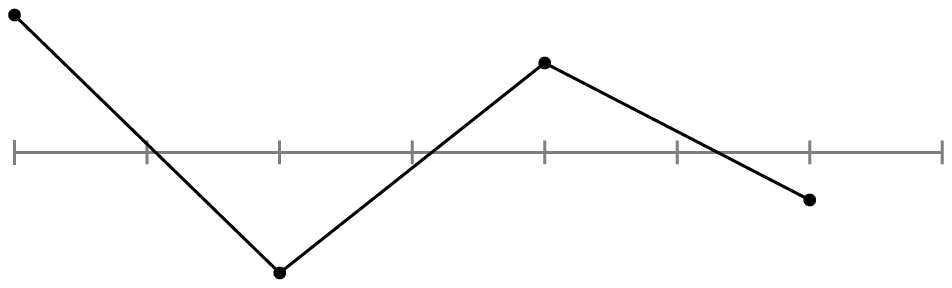}}}
\centerline{
\subfloat[$w_8$ and down--sampled $w_8$]{\includegraphics[width=0.122\textwidth]{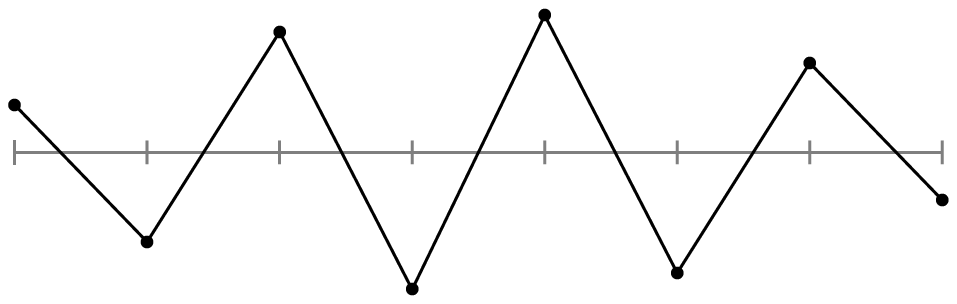} \includegraphics[width=0.122\textwidth]{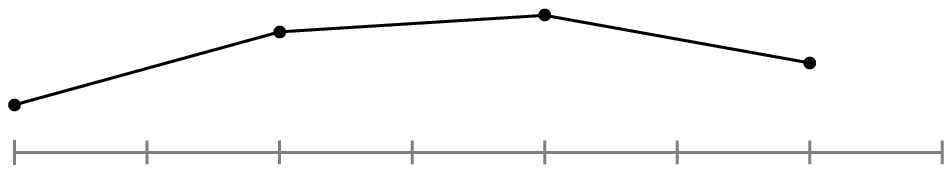}} \hfil
\subfloat[$w_7$ and down--sampled $w_7$]{\includegraphics[width=0.122\textwidth]{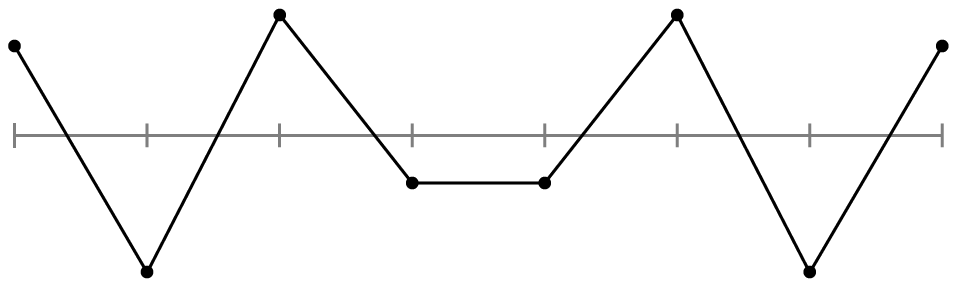} \includegraphics[width=0.122\textwidth]{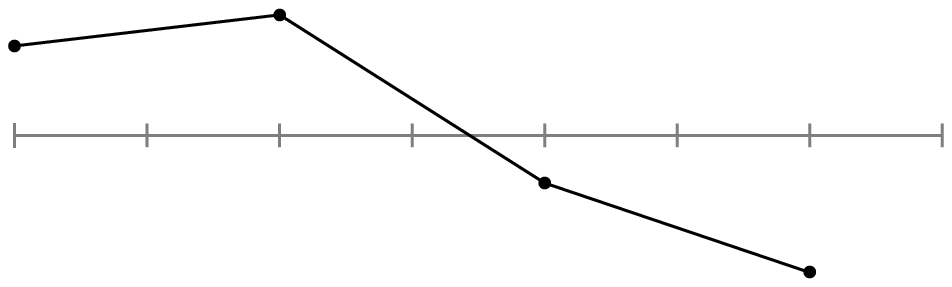}} \hfil
\subfloat[$w_6$ and down--sampled $w_6$]{\includegraphics[width=0.122\textwidth]{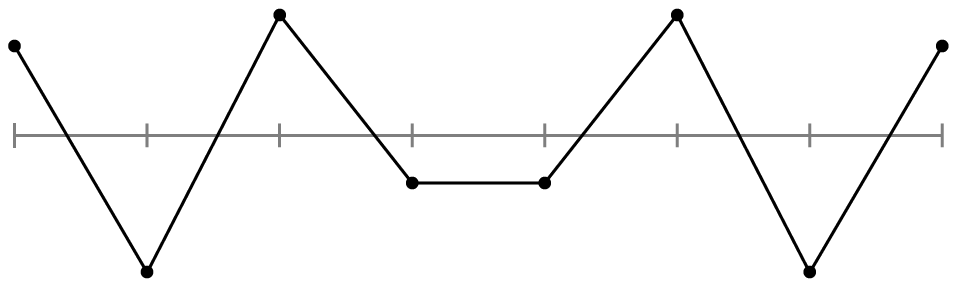} \includegraphics[width=0.122\textwidth]{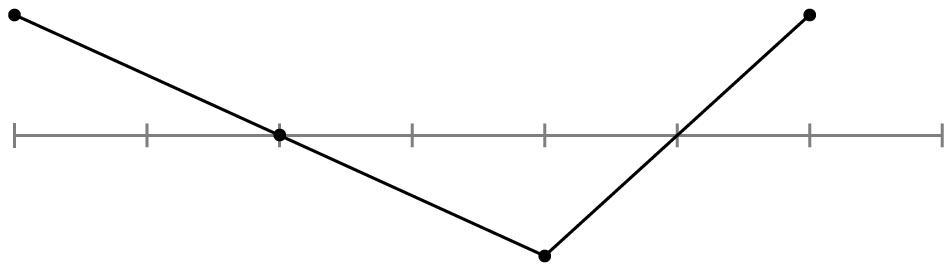}} \hfil
\subfloat[$w_5$ and down--sampled $w_5$]{\includegraphics[width=0.122\textwidth]{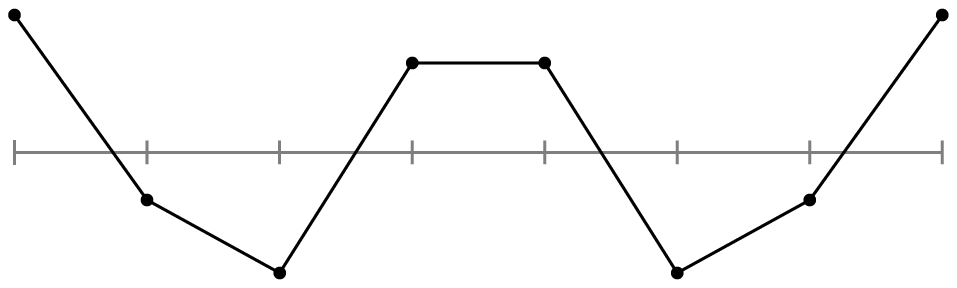} \includegraphics[width=0.122\textwidth]{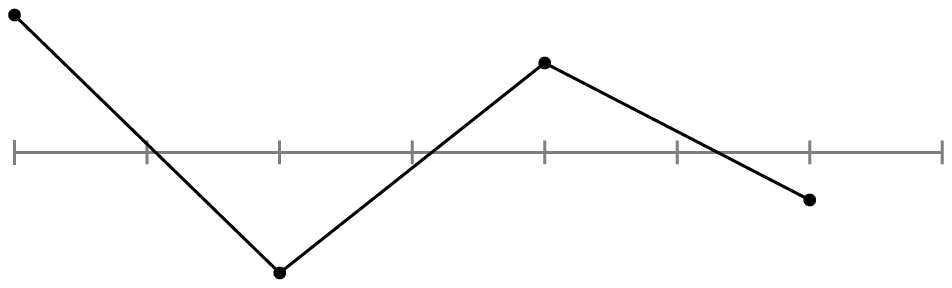}}}
\caption{Example of harmonic aliasing pattern in a harmonic sine basis. The set of vectors $w_j\in\mathbb{R}^8$, $j=1,\ldots,9$, with components $(w_j)_i=\tfrac{2}{9}\sin(\tfrac{ij\pi}{9})$, $i=1,\ldots,9$ forms an orthonormal basis of $\mathbb{R}^8$. A down--sampling operation drops the components with even $i$. The resulting down--sampled eigenvectors are not linearly independent. The down--sampling of $w_k$ is equal to the down--sampling of $w_{9-k}$ with $k=1,\ldots,4$.} \label{fig:red_aliasing_example}
\end{figure*}
The second restriction is on the effects of \emph{down--sampling} the eigenvectors of the system. This operation drops a number of components when applied to vectors and keeps the remaining components untouched. An example is shown in Fig. \ref{fig:red_aliasing_example} where one every two samples of harmonic functions are dropped. By down--sampling the eigenvectors of the system their linear independence is lost, because the dimension of the down--sampled space is less than the number of eigenvectors. This is a general description of the phenomenon of \emph{aliasing}. Here, it will be assumed that there is a subset of $n/2$ components defining a down--sampling operation that makes each down--sampled eigenvector equal (up to sign) to only one of the other down--sampled eigenvectors (see Fig. \ref{fig:red_aliasing_example}). This specific pattern is called \emph{harmonic aliasing} pattern \cite{PNavarrete_2008a} and is found in harmonic functions, which are eigenvectors of \emph{linear space invariant} (LSI) systems\footnote{In numerical analysis a different terminology is used. The \emph{stencil} of an unknown is defined as a geometric arrangements of the non--zero coefficients in the correspondent row in $A$, centered at the diagonal element. This is equivalent to the concept of \emph{impulse response} in signal processing and an LSI system is equivalent to a system with constant stencil coefficients.}, and other basis including at least: Hadamard matrices and eigenvectors of coupled systems of equations \cite{PNavarrete_2008a,PNavarrete_2008c}.

The generality of a system with harmonic aliasing patterns is not known at its largest extent. They were introduced in \cite{PNavarrete_2008a} in order to extend the strong convergence analysis of multi--grid algorithms based on \emph{local Fourier analysis} (LFA). This analysis was introduced by Achi Brandt in the late 70's and remains as the main rigorous tool for the design of multi--grid methods \cite{ABrandt_1994, ABrandt_1977b}. LFA is based on Fourier analysis and is thus restricted to LSI systems, which makes it difficult to use it in many applications. The extended convergence analysis in \cite{PNavarrete_2008a} has not overcome this problem drastically. Nevertheless, in this study its algebraic framework will allow not only the analysis of traditional multi--grid methods but also the introduction of new direct solvers with connections to other important numerical methods.

The contributions of this paper are thus both practical and conceptual. From the practical side, multi--grid algorithms have been successfully applied in many practical problems but the theory behind does not reach the same level. The most common implementation is \emph{algebraic multi--grid} (AMG) which obtains a multi--grid configuration based on heuristics \cite{UTrottenberg_CWOosterlee_ASchueller_2000a}. The numerical methods obtained in this paper represent a step forward on what the theory can achieve. This is, under the assumptions just mentioned, a completely algebraic configuration is found that solves the problem exactly with no use of heuristics. The algorithms are computationally efficient and adaptable to the computational resources (single--core or multi--core).

The conceptual contribution of this paper is the introduction of numerical methods in a setting that establishes a clear connection between multi--rate systems and different multi--level numerical methods. The analysis exploits the structural similarities of the \emph{full multi--grid algorithm} \cite{WLBriggs_VEHenson_SFMcCormick_2000a} and problems of signal reconstruction in \emph{multi--rate systems} \cite{PPVaidyanathan_1993a}. On one hand, both classical and extended convergence analysis for multi--grid solvers have shown the importance of aliasing phenomena to explain how the algorithm converges \cite{PNavarrete_2008a,UTrottenberg_CWOosterlee_ASchueller_2000a}. Here, aliasing appears as an additional source of error that has to be controlled by the algorithm. On the other hand, the problem of signal reconstruction in multi--rate systems follows a different approach. A signal is decomposed in different coarse levels and the question is posed as whether the original signal can be reconstructed from these different pieces of information. Here, each coarse level component has aliasing effects and perfect reconstruction is possible because these effects cancel each other. Given the structural similarities between multi--grid methods and multi--rate system, the central question is whether aliasing cancellation can be used in multi--grid to obtain direct solvers, where the analogue of ``perfect reconstruction'' is ``exact solution.''

The main goal, under the assumptions just mentioned, is to modify a full multi--grid algorithm and obtain a direct multi--grid solver in direct analogy with the problem of perfect reconstruction filters. In order to establish this analogy the problem of perfect reconstruction needs to be generalized for systems with harmonic aliasing patterns, which allows to configure perfect reconstruction filters that are not necessarily LSI systems. On the other hand, the full multi--grid algorithm also requires modifications. Multiple coarse grids are needed in order to keep all the information from the original problem in coarse levels. A partition of the complete set of unknowns defines several coarse levels and represents a particular type of \emph{domain decomposition} \cite{AToselli_OWidlund_2004a}.

Similar approaches can be found in the literature. The use of multiple coarse grids in multi--grid has been introduced by Frederickson and McBryan in \cite{PFrederickson_OMcBryan_1988a} and Hackbusch in \cite{WHackbusch_1988a,WHackbusch_1989a}. It is known that aliasing cancellation helps these methods to converge fast \cite{TFChan_RSTuminaro_1989a}. Nevertheless, none of them work as direct solvers and they still use smoothing iterations. Total reduction methods are the closest in structure to the direct solvers obtained and can be seen as a more restrictive version of one of the algorithms presented.

In section \ref{sec:preliminaries}, full two--grid algorithms and multi--rate systems are reviewed. In section \ref{sec:aliasing}, harmonic aliasing patterns are introduced. In section \ref{sec:qm_filters}, perfect reconstruction filters are studied in the context of harmonic aliasing patterns. In section \ref{sec:convergence}, the convergence analysis of two--grid methods based on grid partitions is studied. In section \ref{sec:d2g}, the problem of finding direct two--grid solvers is stated and solved. In section \ref{sec:dmg}, the multi--grid case is considered. Finally, in section \ref{sec:example} some examples are presented.

%%%%%%%%%%%%%%%%%%%%%%%%%%%%%%%%%%%%%%%%%%%%%%%%%%%%
\section{Preliminaries}
\label{sec:preliminaries}
\begin{figure*}[!t]
\centerline{
\subfloat[Full two--grid algorithm]{{\label{fig:f2g}\includegraphics[width=4.3in]{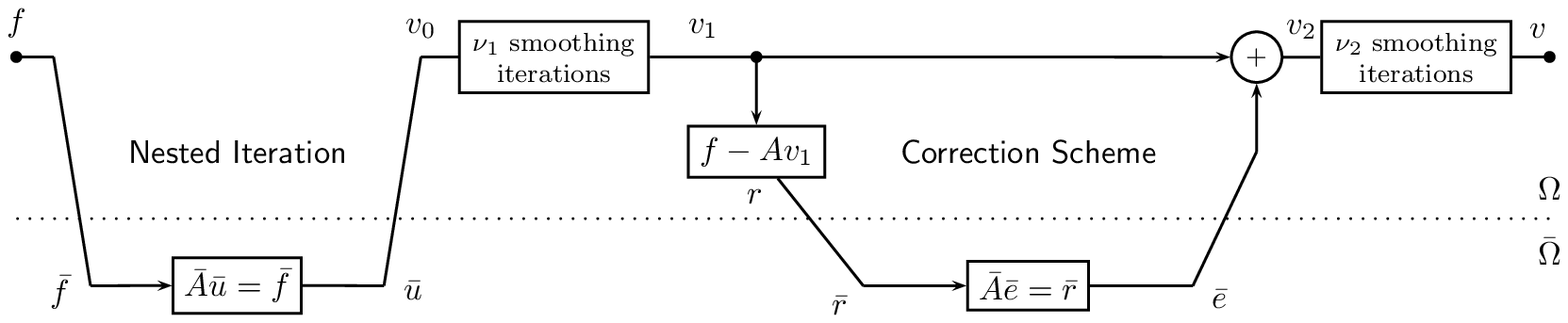}}} \hfil
\subfloat[Two--channel multi--rate system]{{\label{fig:multirate_2ch_infinite}\includegraphics[width=2.6in]{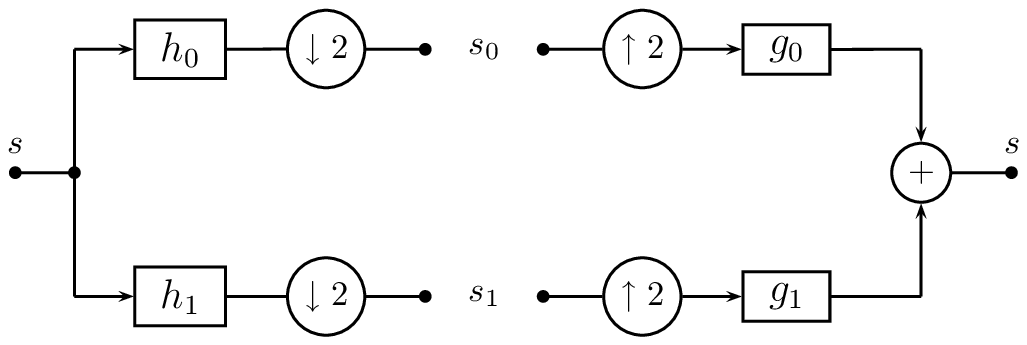}}}}
\caption{The numerical methods introduced in this paper are based on two systems. First, the full two--grid algorithm in Fig. \ref{fig:f2g}. This is an iterative algorithm to obtain an approximate solution of $\finesys u=f$. The dotted line separates vectors from the fine and coarse grid domains. The interpolation (restriction) operation is applied to vectors crossing the dotted line from below (above). Second, the two--channel multi--rate system in Fig. \ref{fig:multirate_2ch_infinite}. Here, each box represents an LSI system and inside the stationary impulse response is shown. The circles with ``$\downarrow 2$'' represent down--sampling operators that drop one every two samples. Similarly, circles with ``$\uparrow 2$'' represent up--sampling operators that insert one zero every two samples.}
\end{figure*}

%---------------------------------------------------
\subsection{Full two--grid algorithm}
The full two--grid algorithm is an iterative solver used to obtain approximate solutions of (\ref{eq:fine_system}). In order to simplify the problem, the algorithm uses the concept of \emph{grids} and \emph{coarse grids}. Whatever the nature of the problem is, the system (\ref{eq:fine_system}) can always be associated with a graph in which the unknowns of the system are the nodes of the graph. The nodes of the graph are associated with a set of labels $\finegrid$ which is called the \emph{fine grid}. A \emph{coarse grid}, $\rcoarsegrid$, is a proper subset of the fine grid; i.e., $\rcoarsegrid\subset\finegrid$.

The so--called \emph{inter--grid operators} are defined as any linear transformation between scalar fields on $\finegrid$ and $\rcoarsegrid$. That is $\rinterp\in\mathbb{C}^{n\times|\rcoarsegrid|}$ and $\rrestri\in\mathbb{C}^{|\rcoarsegrid|\times n}$, where $\rinterp$ is the \emph{interpolation operator} and $\rrestri$ is the \emph{restriction operator}. In addition to these operations, and following the standard of most multi--grid applications, a coarse system matrix is defined following the \emph{Galerkin condition} \cite{WLBriggs_VEHenson_SFMcCormick_2000a}
\begin{equation}
\rcoarsesys \define \rrestri\; \finesys\; \rinterp\;. \label{eq:galerkin}
\end{equation}

A \emph{full two--grid algorithm} solves (\ref{eq:fine_system}) by using the system shown in Fig. \ref{fig:f2g}. Here, there are three steps involved. First, the two boxes perform fixed numbers of \emph{smoothing} or \emph{stationary iterative} iterations. These iterations --typically Gauss--Seidel, Jacobi, Richardson, etc.-- are known to obtain good local approximations of the solution
\cite{WLBriggs_VEHenson_SFMcCormick_2000a}. This means that high--frequency components of the \emph{approximation error}, $e_k=u-v_k$, are efficiently reduced. In each smoothing iteration the approximation error evolves as $e_{k+1} = S e_k$. The matrix $S$ is thus called the \emph{smoothing operator}. For better understanding of this step is convenient to assume that  $S$ is a \emph{filter}, or \emph{Fourier multiplier}. If this is the case then $S$ has same eigenvectors as $\finesys$ and it can be decomposed as $S=W\Sigma V^H$. The matrix of eigenvalues, $\Sigma$, contains the \emph{frequency response}, or \emph{symbols}, of the smoothing operator. The smoothing effect on the approximation error is seen in the frequency domain as a damping effect concentrated on high--frequency eigenvectors. In other words, the eigen--decomposition of $S$ can be written in block--form as
\begin{equation}
S = W \left[\begin{array}{cc}
\Sigma_L & \\
 & \Sigma_H
\end{array}\right] V^H \;,
\end{equation}
where $\Sigma_L$ and $\Sigma_H$ correspond to low-- and high--frequency eigenvalues close to $1$ and $0$, respectively.

The remaining steps in Fig. \ref{fig:f2g} make use of the coarse grid. First, the so--called \emph{nested iteration} step takes $f$ and computes an initial approximation $v_0$. This step solves a coarse grid equation using a restricted version of the source vector, $\bar{f}$. And second, the so--called \emph{correction scheme} improves the approximation $v_1$ to obtain $v_2$. This step computes an approximation of $e_k$ from the error equation $\finesys e_k=r_k$, where $r_k=f-\finesys v_k$ is the \emph{residual vector} taking the role of the source vector. A coarse grid equation is solved using a restricted version of the source vector, $\bar{r}_k$. Once the approximation is obtained, it is added to correct the current approximation.

It is observed that the approximation of $e_k$ obtained in the coarse grid equations effectively represents its low--frequency components and fails to represents its high--frequency components. This is a consequence of the computations in coarse grids where the big picture (low--frequencies) of the solution is clear but details (high--frequencies) are lost.

The approximation error in the coarse grid steps evolves as
\begin{equation}
e_0 = \rcgc u \quad\;\mbox{and}\quad e_{k+1} = \rcgc e_k \;, \label{eq:two_grid_error}
\end{equation}
for nested iterations and the correction scheme, respectively. Here, the matrix
\begin{equation} \label{eq:coarse_grid_correction_matrix}
\rcgc = I-\rinterp\rcoarsesys^{-1}\rrestri \finesys
\end{equation}
determines the evolution of the approximation error and is called the \emph{coarse grid correction} matrix \cite{PWesseling_1992a}.

\begin{figure*}[!t]
\centerline{\subfloat[Two--grid filtering effect]{\label{fig:2grid_filter}\includegraphics[width=0.48\textwidth]{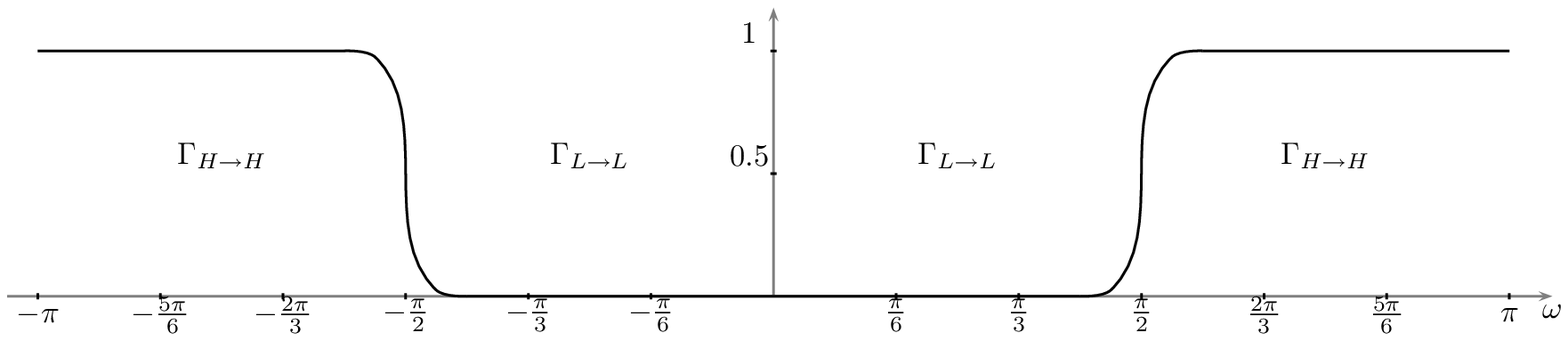}} \hfil
\subfloat[Two--grid aliasing effect]{\label{fig:2grid_alias}\includegraphics[width=0.48\textwidth]{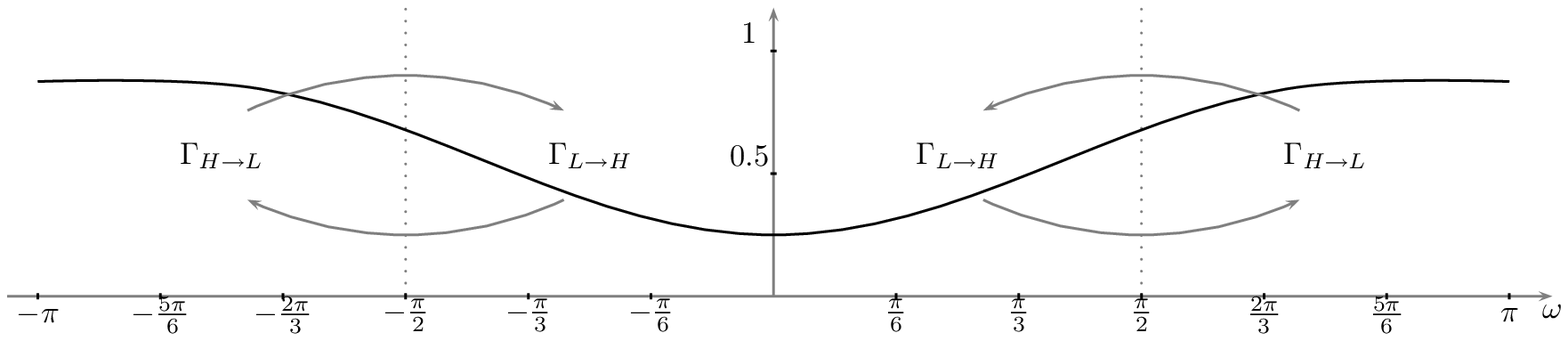}}}
\caption{Representation of the frequency effects of a coarse grid correction operator in a two--grid algorithm. The operator is not a filter and thus is not completely represented by filtering effects. The operator is completely represented by filtering and aliasing effects. $\rmodal_{L\rightarrow L}$ indicates the region of low--frequencies where the symbols are close to $0$. $\rmodal_{H\rightarrow H}$ indicates the region of high--frequencies where the symbols are close to $1$. $\rmodal_{L\rightarrow H}$ and $\rmodal_{H\rightarrow L}$ multiply low-- and high--frequencies and take them into high-- and low--frequencies, respectively.} \label{fig:cgc_frequency}
\end{figure*}

The reduction of low--frequency components of the error in two--grid steps suggests that $\rcgc$ is a \emph{high--pass filter}. If this would be the case then there would be a decomposition $\rcgc=W\rmodal V^H$ where $\rmodal$ is a diagonal matrix. One would expect a frequency response of the filter with the shape shown in Fig. \ref{fig:2grid_filter}. Unfortunately, this is not the case. $\rcgc$ is not a filter because of aliasing effects.

Although $\rcgc$ is not technically a \emph{filter} (or \emph{Fourier multiplier}), under the assumption of harmonic aliasing patterns a decomposition $\rcgc=W\rmodal V^H$ exists where $\rmodal$ is sparse (but not diagonal). In the forthcoming sections it will be shown that
\begin{equation} \label{eq:cgc_decomposition}
\rcgc = W \left[\begin{array}{cc}
\rmodal_{L\rightarrow L} & \rmodal_{H\rightarrow L} \\
\rmodal_{L\rightarrow H} & \rmodal_{H\rightarrow H}
\end{array} \right] V^H \;,
\end{equation}
where $\rmodal_{L\rightarrow L}$, $\rmodal_{H\rightarrow L}$, $\rmodal_{L\rightarrow H}$ and $\rmodal_{H\rightarrow H}$ are all diagonal matrices if one assumes harmonic aliasing patterns. The filtering effect is contained in $\rmodal_{L\rightarrow L}$ and $\rmodal_{H\rightarrow H}$, and are expected to be as shown in Fig. \ref{fig:2grid_filter}. As opposed to a proper filter, this figure does not tell everything about the coarse grid correction matrix. A second graphic, shown in Fig. \ref{fig:2grid_alias}, must show the aliasing effect from $\rmodal_{H\rightarrow L}$ and $\rmodal_{L\rightarrow H}$.

%---------------------------------------------------
\subsection{Two--channel multi--rate systems}
In multi--rate systems one is interested to decompose a discrete signal in different components at lower sampling rates \cite{PPVaidyanathan_1993a,MVetterli_JKovacevic_1995a}. A two--channel multi--rate system is shown in Fig. \ref{fig:multirate_2ch_infinite}. Two restriction operations are performed by filtering and down--sampling. These operations split the original signal, $s$, into two signals, $s_0$ and $s_1$, each with half of the original samples. The original signal is recovered by summing two interpolation operations performed by inserting zeros (up--sampling) and then filtering.

Here, the boxes represent LSI systems which have stationary impulse responses and are filters with respect to a harmonic basis. Their frequency responses are given by the Fourier transforms of their impulse responses: $H_0(\omega)$, $H_1(\omega)$, $G_0(\omega)$ and $G_1(\omega)$. In practical applications the problem is whether the support of the frequency responses can overlap and still be able to recover the original signal. This is possible when the filters fulfill the following conditions by Vetterli \cite{MVetterli_1986a}
\begin{align}
G_0(\omega) H_0(\omega) + G_1(\omega) H_1(\omega) & = 2 \;, \label{eq:multirate_2ch_infinite_perfect_reconstruction} \\
G_0(\omega) H_0(\omega-\pi) + G_1(\omega) H_1(\omega-\pi) & = 0 \;. \label{eq:multirate_2ch_infinite_aliasing_free}
\end{align}
Here, condition (\ref{eq:multirate_2ch_infinite_aliasing_free}) causes the aliasing effects from different channels to cancel each other and (\ref{eq:multirate_2ch_infinite_perfect_reconstruction}) causes the final sum to be equal to the original signal. The more general result with an arbitrary number of decompositions is due to Vaidyanathan \cite{PPVaidyanathan_1987a}.

%%%%%%%%%%%%%%%%%%%%%%%%%%%%%%%%%%%%%%%%%%%%%%%%%%%%
\section{Red--Black Harmonic Aliasing}
\label{sec:aliasing}

\begin{figure*}[!t]
\centerline{
\subfloat[Fine grid with nodes colored red and black.]{\label{fig:fine_grid}\raisebox{0.15in}{\includegraphics[width=0.5\textwidth]{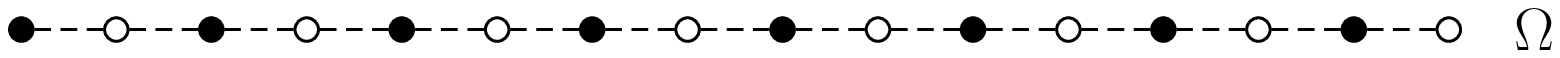}}} \hfil
\subfloat[Red and black coarse grids]{\label{fig:partition}\includegraphics[width=0.5\textwidth]{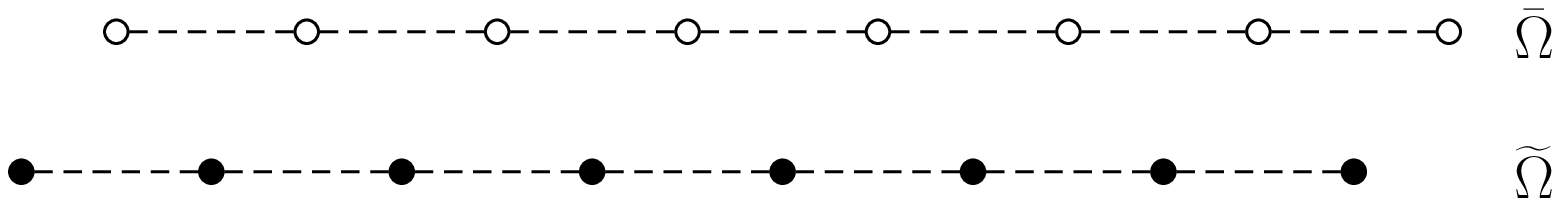}}}
\caption{Red--black partition of a grid in 1D. The fine grid $\finegrid$ is partitioned into a red coarse grid $\rcoarsegrid$ and a black coarse grid $\bcoarsegrid$.} \label{fig:red_black_partition}
\end{figure*}

A \emph{red} coarse grid, $\rcoarsegrid$, and a \emph{black} coarse grid, $\bcoarsegrid$, are defined by a partition of the fine grid. This is, $\finegrid=\rcoarsegrid+\bcoarsegrid$, where the sum denotes the disjoint union of two sets. An example is shown in Fig. \ref{fig:red_black_partition}. The motivation of the red--black partition is to keep track of all the fine grid nodes in coarser grids. In this way the partition represents a particular type of \emph{domain decomposition} \cite{AToselli_OWidlund_2004a}.

The selection of nodes to the red and black partition will be represented by \emph{down--sampling} operators according to the following definition.
\begin{defn}[Down/Up--sampling matrices]
The red and black \emph{down--sampling} matrices are defined as $\rdown \in \{0,1\}^{|\rcoarsegrid| \times n}$ and $\bdown \in \{0,1\}^{|\bcoarsegrid| \times n}$ such that
\begin{align}
(\rdown)_{i,j} & \define \left\{\begin{array}{ll}
1 & \mbox{if node $j\in\finegrid$ is the $i^{th}$ red node in $\rcoarsegrid$} \\
0 & \mbox{otherwise}
\end{array}\right.
\end{align}
and
\begin{align}
(\bdown)_{i,j} & \define \left\{\begin{array}{ll}
1 & \mbox{if node $j\in\finegrid$ is the $i^{th}$ black node in $\bcoarsegrid$} \\
0 & \mbox{otherwise}
\end{array}\right.,
\end{align}
respectively. Similarly, the red and black \emph{up--sampling} matrices are defined as $\rup \in \{0,1\}^{n \times |\rcoarsegrid|}$ and $\bup \in \{0,1\}^{n \times |\bcoarsegrid|}$ such that $\rup=\rdown^T$ and $\bup=\bdown^T$, respectively.
\end{defn}

The down--sampling matrix of a certain color represents a linear transformation that takes a fine grid vector and drops all the values that correspond to a node of different color. The up--sampling matrix takes a coarse grid vector and inserts zeros at the new nodes (of different color) in the fine grid. From this interpretation, a set of basic properties follows
\begin{align}
& \rdown \rup=I \quad\mbox{and}\quad \bdown\bup=I \;, \label{eq:DU_identity} \\
& \bdown \rup=0 \quad\mbox{and}\quad \rdown\bup=0 \;,\quad\mbox{and} \label{eq:DU_zero} \\
& \rup\rdown+\bup\bdown=I \;. \label{eq:UD_completion} 
\end{align}

In \cite{PNavarrete_2008a} a so--called \emph{harmonic aliasing} pattern was defined only for the red grid. The following definition considers both red and black grids. This addition has important consequences in the results to come.

\begin{defn}[Red and Black Harmonic Aliasing Patterns] \label{def:rbHAP}
A matrix $M\in\mathbb{C}^{n\times n}$ is said to have \emph{red} and \emph{black harmonic aliasing patterns} if it is diagonalizable, with biorthogonal eigenvectors $W$ and $V$, and there exists a red--black partition which divides the domain into two halves, with down--sampling matrices $\rdown  \in \{0,1\}^{\tfrac{n}{2} \times n}$ and $\bdown  \in \{0,1\}^{\tfrac{n}{2} \times n}$, such that
\begin{equation} \label{eq:rbHAP}
V^H \rup \rdown W = \rhap \quad\mbox{and}\quad V^H \bup \bdown W = \bhap \;,
\end{equation}
respectively. Here, $\rhap$ and $\bhap$ are the \emph{red} and \emph{black harmonic aliasing patterns} defined, respectively, as
\begin{equation} \label{eq:redblack_aliasing_patterns}
\rhap \define \frac{1}{2}
\left[\begin{array}{cc}
I & I \\
I & I
\end{array}\right] \quad\mbox{and}\quad
\bhap \define \frac{1}{2}
\left[\begin{array}{cc}
I & -I \\
-I & I
\end{array}\right] \;.
\end{equation}
\end{defn}

In this definition the red and black harmonic aliasing patterns appear as independent properties. The following statement shows how these definitions are equivalent.

\begin{pro}
A matrix with red harmonic aliasing pattern has black harmonic aliasing pattern, and vice versa. Therefore, a matrix with these properties is said to have a red--black harmonic aliasing pattern.
\end{pro}
\begin{proof}
Taking (\ref{eq:UD_completion}), pre--multiplied by $V^H$ and post--multiplied by $W$ gives
\begin{equation}
V^H \rup \rdown W + V^H \bup \bdown W = I \;.
\end{equation}
Thus, if $V^H \rup \rdown W=\rhap$ then $V^H \bup \bdown W=\bhap$, and vice versa.
\end{proof}

The definition of red--black harmonic aliasing pattern is convenient for algebraic manipulation but the connection with the common concept of aliasing is not clear yet. In the following theorem an alternative and equivalent definition is given which makes this connection explicit.

\begin{figure*}[!t]
\centerline{
\subfloat[Low--frequency harmonic function $w_L$.]{\includegraphics[width=0.5\textwidth]{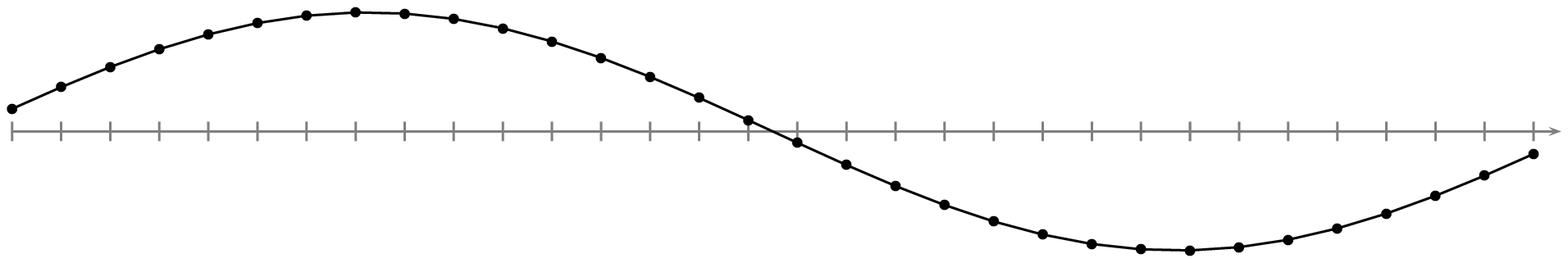}} \hfil
\subfloat[High--frequency harmonic function $w_H$.]{\includegraphics[width=0.5\textwidth]{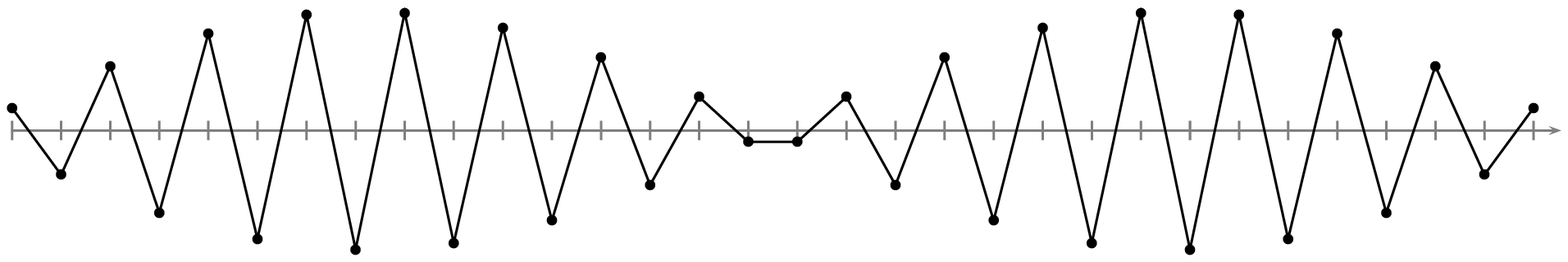}}}
\centerline{
\subfloat[$\rdown w_L$]{\includegraphics[width=0.25\textwidth]{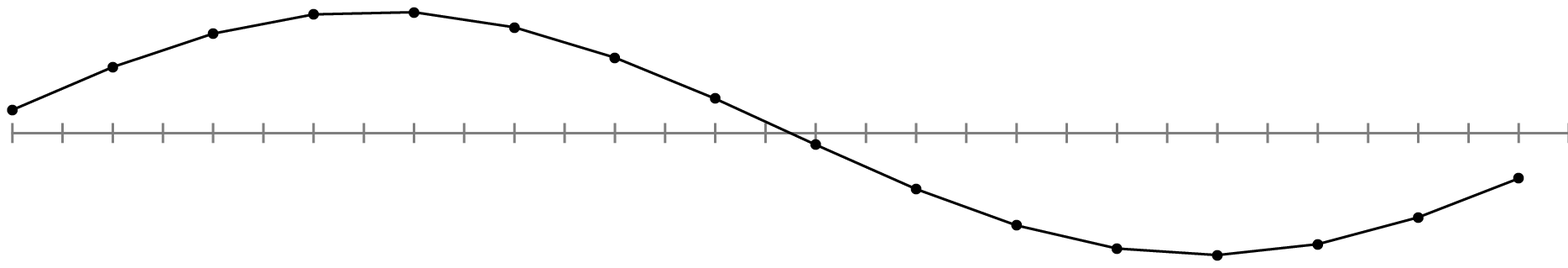}} \hfil
\subfloat[$\bdown w_L$]{\includegraphics[width=0.25\textwidth]{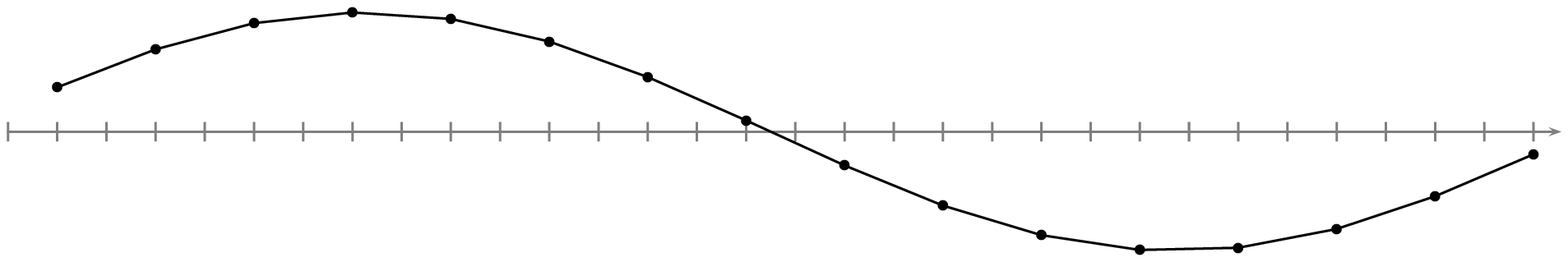}} \hfil
\subfloat[$\rdown w_H$]{\includegraphics[width=0.25\textwidth]{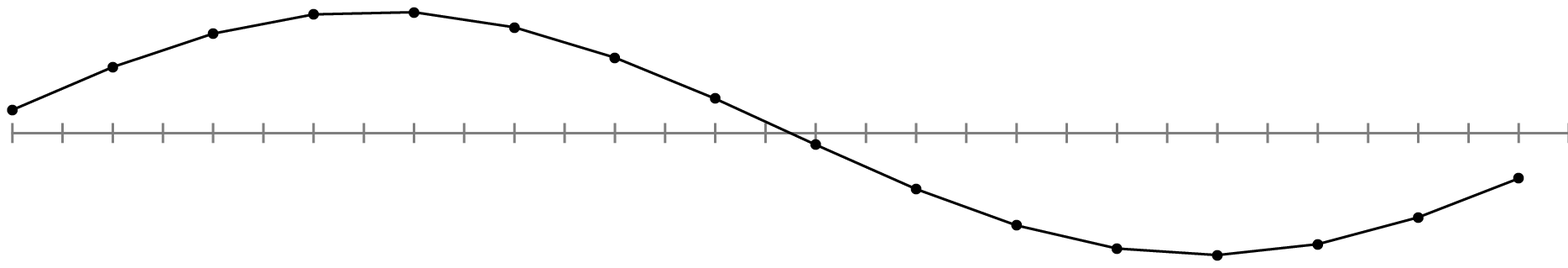}} \hfil
\subfloat[$\bdown w_H$]{\includegraphics[width=0.25\textwidth]{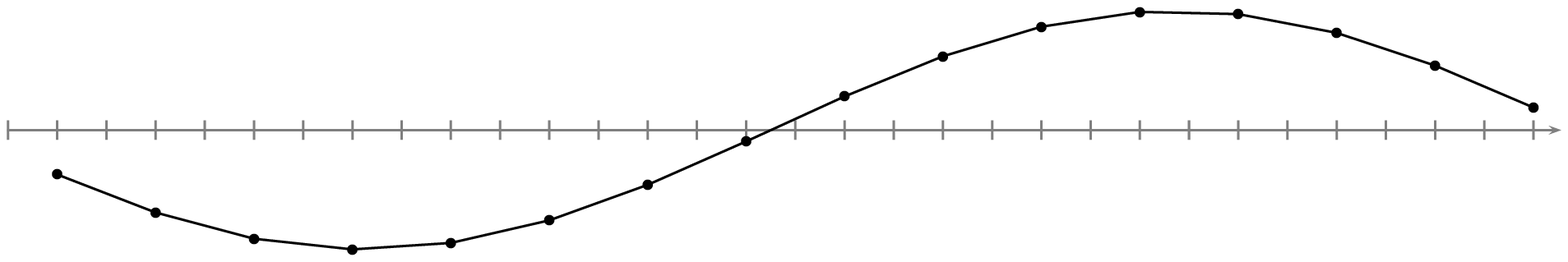}}}
\caption{Difference between red and black harmonic aliasing patterns. The red down--sampling of $w_L$ and $w_H$ are equal, $\rdown w_L = \rdown w_H$. Whereas the black down--sampling of $w_L$ is the negative of the black down--sampling of $w_H$, $\bdown w_L = -\bdown w_H$.} \label{fig:red_black_aliasing_example}
\end{figure*}
\begin{thm}[Navarrete and Coyle] \label{thm:red_black_surjective}
A matrix $M\in\mathbb{C}^{n\times n}$ with red--black harmonic aliasing pattern is equivalent to have a diagonalizable matrix, with biorthogonal eigenvectors $W$ and $V$, for which there exists a red--black partition dividing the domain into two halves, with down--sampling matrices $\rdown  \in \{0,1\}^{\tfrac{n}{2} \times n}$ and $\bdown  \in \{0,1\}^{\tfrac{n}{2} \times n}$, and such that there is an ordering of the eigenvectors for which the partitions $W=\left[W_L W_H\right]$ and $V=\left[V_L V_H\right]$ fulfill the conditions
\begin{align}
\rdown W_L & = \rdown W_H \;, & \rdown V_L & = \rdown V_H \;, \label{eq:red_surjective} \\
\bdown W_L & = - \bdown W_H \quad\mbox{and} & \bdown V_L & = - \bdown V_H. \label{eq:black_surjective}
\end{align}
\end{thm}

The proof of this theorem is partially contained in \cite{PNavarrete_2008a} were only the red coarse grid was considered. The result for the black coarse grid is shown in Appendix \ref{app:proof_red_black_surjective}.

The sign difference between the red and black harmonic aliasing patterns is explained in Fig. \ref{fig:red_black_aliasing_example} and it represents the fact that harmonic basis vectors are composed of two envelopes which, intermixed with the same sign form a low frequency and, intermixed with opposed signs form a high--frequency.

The definition of red--black harmonic aliasing patterns for a given matrix does not involve its eigenvalues. But, in the algebra derived from the partition of eigenvectors in Theorem \ref{thm:red_black_surjective} it will be necessary to specify which eigenvalues are associated with each partition. The following remark introduces the notation to make this distinction clear.

\begin{rmk} \label{rem:partition_eigenvalues}
For a matrix $M\in\mathbb{C}^{n\times n}$ with red--black harmonic aliasing patterns and eigen--decomposition $M=W E V^H$, the partition of eigenvectors, $W=\left[W_L W_H\right]$ and $V=\left[V_L V_H\right]$, induces a partition of eigenvalues such that
$E = \left[\begin{smallmatrix} E_L & \\ & E_H \end{smallmatrix}\right]$,
where $E_L$ and $E_H$ are the diagonal matrices of eigenvalues associated with $L$ and $H$ eigenvectors, respectively.
\end{rmk}

%%%%%%%%%%%%%%%%%%%%%%%%%%%%%%%%%%%%%%%%%%%%%%%%%%%%
\section{Perfect Reconstruction Filters for Systems with Harmonic Aliasing Patterns}
\label{sec:qm_filters}

In the context of finite discrete systems we want to extend the idea of perfect reconstruction filters for systems with harmonic aliasing patterns, which are more general than LSI systems. First, we define an object that will extend the operation of quadrature and conjugate mirror filters \cite{ACroisier_DEsteban_CGaland_1976a,FMintzer_1985a,MJSmith_TPBarnwellIII_1984a}.

\begin{defn}[Mirror Matrix] \label{def:mirror_operator}
The \emph{mirror} of a matrix $M\in\mathbb{C}^{n\times n}$ with respect to a red--black partition represented by down/up--sampling matrices $\rdown$, $\rup=\rdown^T$, and $\bdown$, $\bup=\bdown^T$, is defined as
\begin{equation}
M^\star \define (\rup\rdown-\bup\bdown) M (\rup\rdown-\bup\bdown) \;.
\end{equation}
\end{defn}

Matrices $\rup\rdown\in\{0,1\}^{n\times n}$ and $\bup\bdown\in\{0,1\}^{n\times n}$ are diagonal with one's whenever $i=j$ is a red or black node, respectively, and zero otherwise. Therefore, $\rup\rdown-\bup\bdown$ is a diagonal matrix that takes the value $1$ when $i=j$ is a red node, and takes the value $-1$ when $i=j$ is a black node.

For an LSI system with stationary impulse response $h[n]$, the mirror operator with respect to a uniform down--sampling by factor of $2$ (where red nodes are even nodes) produces an impulse response $h^\star[n]=(-1)^n h[n]$. This filter has symbols $H^\star(\omega)=H(\omega-\pi)$ and thus swaps low and high frequencies. The following proposition generalizes this result to systems with red--black harmonic aliasing patterns.

\begin{pro} \label{thm:frequency_inversion}
The mirror of a matrix $M\in\mathbb{C}^{n\times n}$ with red--black harmonic aliasing patterns and eigen--decomposition $M=W\left[\begin{smallmatrix}E_L & \\ & E_H \end{smallmatrix}\right]V^H$, is a filter with eigen--decomposition
\begin{equation}
M^\star = W \left[\begin{array}{cc} E_H & \\ & E_L \end{array}\right] V^H \;.
\end{equation}
\end{pro}
\begin{proof}
Using the eigen--decomposition of $M$ and the definition of red and black harmonic aliasing patterns, the result is obtained as follows
\begin{equation}
\begin{aligned}
M^\star & = W V^H(\rup\rdown-\bup\bdown) W\; E\; V^H (\rup\rdown-\bup\bdown) W V^H \\
 & = W \left[\begin{smallmatrix} 0 & I \\ I & 0 \end{smallmatrix}\right] \; \left[\begin{smallmatrix} E_L & 0 \\ 0 & E_H \end{smallmatrix}\right] \; \left[\begin{smallmatrix} 0 & I \\ I & 0 \end{smallmatrix}\right] V^H \\
 & = W \left[\begin{smallmatrix} E_H & \\ & E_L \end{smallmatrix}\right] V^H \;.
\end{aligned}
\end{equation}
\end{proof}

\begin{figure}[!t]
\centering
\includegraphics[width=3in]{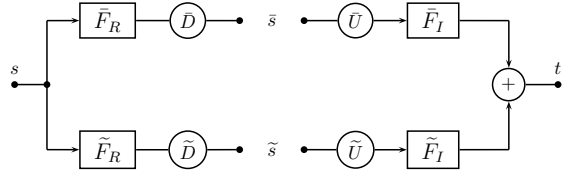}
\caption{Two--channel finite multi--rate system. The system takes a signal $s\in\mathbb{C}^n$, decomposes it into $\bar{s}\in\mathbb{C}^{n/2}$ and $\widetilde{s}\in\mathbb{C}^{n/2}$, and recovers $t\in\mathbb{C}^n$. A perfect reconstruction system is such that $t=s$. Here, $\rrfilter$ and $\brfilter$ are restriction filters; $\rifilter$ and $\bifilter$ are interpolation filters; $\rdown$ and $\bdown$ represent down--sampling operations; and, $\rup$ and $\bup$ represent up--sampling operations.} \label{fig:multirate_2ch_finite}
\end{figure}

These results are all we need to proceed into the results of the next sections. Now, in order to show the analogy between perfect reconstruction filters and direct multi--grid solvers we have to look back to the problem of perfect reconstruction and see how the main results look for systems with harmonic aliasing patterns. First, we need to restrict the problem to discrete signals in $\mathbb{C}^n$. The two--channel multi--rate system in Fig. \ref{fig:multirate_2ch_finite} is then the finite and discrete version of the system in Fig. \ref{fig:multirate_2ch_infinite}. The problem is to find interpolation and restriction matrices $\rifilter$, $\rrfilter$, $\bifilter$ and $\brfilter$ in $\mathbb{C}^{n\times n}$ that are filters with respect to a biorthogonal basis $W$ and $V$, and such that we have perfect reconstruction; i.e., $t=s$.

The following theorem restates Vetterli's conditions (\ref{eq:multirate_2ch_infinite_perfect_reconstruction}) and (\ref{eq:multirate_2ch_infinite_aliasing_free}) for this new context.

\begin{thm}[Vetterli]
Let $\rifilter$, $\rrfilter$, $\bifilter$ and $\brfilter$ be filters with respect to a biorthogonal basis $W$ and $V$ in $\mathbb{C}^{n\times n}$. Let their matrices of eigenvalues be $\rifilterev$, $\rrfilterev$, $\bifilterev$ and $\brfilterev$, respectively, and the eigenvectors have red--black harmonic aliasing patterns with respect to the down--sampling matrices $\rdown$ and $\bdown$. Then, the multi--rate system in Fig. \ref{fig:multirate_2ch_finite} has the perfect reconstruction property, $t=s$, if and only if
\begin{align}
\rifilterev \rrfilterev + \bifilterev \brfilterev & = I \;, \label{eq:multirate_2ch_finite_perfect_reconstruction}\\
\rifilterevL \rrfilterevH - \bifilterevL \brfilterevH & = 0 \;,\;\mbox{and} \label{eq:multirate_2ch_finite_aliasing_freeLH}\\
\rifilterevH \rrfilterevL - \bifilterevH \brfilterevL & = 0 \;, \label{eq:multirate_2ch_finite_aliasing_freeHL}
\end{align}
where the $L$ and $H$ subindexes follow the notation introduced in Remark \ref{rem:partition_eigenvalues}.
\end{thm}
\begin{proof}
From Fig. \ref{fig:multirate_2ch_finite}, perfect reconstruction is obtained if and only if
\begin{equation}
\rifilter \rup \rdown \rrfilter + \bifilter \bup \bdown \brfilter = I \;.
\end{equation}
Pre--multiplying by $V^H$, post--multiplying by $W$ and using the definitions in (\ref{eq:rbHAP}) gives
\begin{equation}
\left[\begin{smallmatrix}
\rifilterevL \rrfilterevL + \bifilterevL \brfilterevL & \rifilterevL \rrfilterevH - \bifilterevL \brfilterevH \\
\rifilterevH \rrfilterevL - \bifilterevH \brfilterevL & \rifilterevH \rrfilterevH + \bifilterevH \brfilterevH
\end{smallmatrix}\right] =
\left[\begin{smallmatrix}
I & 0 \\
0 & I
\end{smallmatrix}\right] \;.
\end{equation}
The blocks in the diagonal are equivalent to (\ref{eq:multirate_2ch_finite_perfect_reconstruction}) and the off--diagonals give (\ref{eq:multirate_2ch_finite_aliasing_freeLH}) and (\ref{eq:multirate_2ch_finite_aliasing_freeHL}).
\end{proof}

Next, we are interested in a particular solution of these conditions using mirror filters. The following corollary generalizes the solutions using \emph{quadrature} and \emph{conjugate mirror} filters (QMF and CMF) for perfect reconstruction in a two--channel multi--rate system with respect to a biorthogonal basis with red--black harmonic aliasing patterns.

\begin{cor}[Croisier et al. / Smith and Barnwell III / Mintzer]
Perfect reconstruction conditions (\ref{eq:multirate_2ch_finite_perfect_reconstruction}), (\ref{eq:multirate_2ch_finite_aliasing_freeLH}) and (\ref{eq:multirate_2ch_finite_aliasing_freeHL}) are fulfilled if $\rifilter$ is such that
\begin{equation} \label{eq:mirror_property}
\rifilterevL^2+\rifilterevH^2=I \;,\quad\mbox{and}
\end{equation}
\begin{equation} \label{eq:qmf_perfect_reconstruction}
\rrfilter = \rifilter \;, \quad \bifilter = \rifilter^\star \;,\; \mbox{and}\quad \brfilter = \rifilter^\star \;.
\end{equation}
\end{cor}
\begin{proof}
Using Theorem \ref{thm:frequency_inversion}, if (\ref{eq:qmf_perfect_reconstruction}) is assumed then the conditions (\ref{eq:multirate_2ch_finite_aliasing_freeLH}) and (\ref{eq:multirate_2ch_finite_aliasing_freeHL}) are fulfilled. The mirror condition (\ref{eq:mirror_property}) is necessary to fulfill (\ref{eq:multirate_2ch_finite_perfect_reconstruction}).
\end{proof}

Quadrature mirror filters were introduced as a solution of the perfect reconstruction problem in \cite{ACroisier_DEsteban_CGaland_1976a}. In its original formulation, only a Haar filter gives a sparse quadrature mirror filter \cite{SMallat_1998a}. Later, this problem was solved by the introduction of conjugate mirror filters \cite{MJSmith_TPBarnwellIII_1984a,FMintzer_1985a}, and later generalized for biorthogonal and multi--channel paraunitary systems \cite{MVetterli_1986a,PPVaidyanathan_1987a}. The problem to obtain perfect reconstruction FIR filters has to do with the delays between the filters. From a signal processing perspective the red--black partition corresponds to a \emph{polyphase decomposition} \cite{PPVaidyanathan_1993a,MVetterli_JKovacevic_1995a,GStrang_TNguyen_1996a} that introduces the delays in the down/up--sampling operations and allows to condense the different solutions for perfect reconstruction.

%%%%%%%%%%%%%%%%%%%%%%%%%%%%%%%%%%%%%%%%%%%%%%%%%%%%
\section{Two--grid Convergence}
\label{sec:convergence}

Since the configuration of a two--grid algorithm involves many parameters and assumptions, it is convenient to introduce a single terminology that refers to this configuration.
\begin{defn}[Red--Black Harmonic Two--grid Configuration]
A \emph{red--black harmonic two--grid configuration} is a set of matrices depending on a system matrix $A\in\mathbb{C}^{n\times n}$ with red--black harmonic aliasing patterns and biorthogonal eigenvectors $W$ and $V$. The red--black partition of nodes is represented by down--sampling operators $\rdown$ and $\bdown$. The configuration is completed with the inter--grid filters $\rrfilter$, $\rifilter$, $\brfilter$ and $\bifilter$, all in $\mathbb{C}^{n\times n}$. Then, a series of matrices associated with the configuration are defined. First, the interpolation and restriction operators, and the two coarse system matrices given by the Galerkin condition, are defined as
\begin{align}
\rrestri    & = \rdown \rrfilter \;,     & \rinterp        & = \rifilter \rup \;, & \brestri    & = \bdown \brfilter \;,         & \binterp & = \bifilter \bup \;, \label{eq:igrid_decomposition} \\
\rcoarsesys & = \rrestri\finesys\rinterp & \quad\mbox{and} &                      & \bcoarsesys & = \brestri\finesys\binterp \;. &          &
\end{align}
The system matrix and the inter--grid filters have eigen--decompositions
\begin{equation*}
\begin{aligned}
\finesys & = W\; \Lambda\; V^H \;,\quad & \rrfilter & = W\; \rrfilterev\; V^H \;,                & \rifilter & = W\; \rifilterev\; V^H \;, \\
         &                              & \brfilter & = W\; \brfilterev\; V^H \;\;\;\;\mbox{and} & \bifilter & = W\; \bifilterev\; V^H \;,
\end{aligned}
\end{equation*}
and the partition of eigenvectors, $W=\left[W_L W_H\right]$ and $V=\left[V_L V_H\right]$, leads to the partitions of eigenvalues
\begin{equation*}
\begin{aligned}
\Lambda     & = \left[\begin{smallmatrix} \Lambda_L &    \\ & \Lambda_H    \end{smallmatrix}\right] , &
\rrfilterev & = \left[\begin{smallmatrix} \rrfilterevL & \\ & \rrfilterevH \end{smallmatrix}\right] , &
\rifilterev & = \left[\begin{smallmatrix} \rifilterevL & \\ & \rifilterevH \end{smallmatrix}\right] , \\
 & &
\brfilterev & = \left[\begin{smallmatrix} \brfilterevL & \\ & \brfilterevH \end{smallmatrix}\right] \;\;\mbox{and} &
\bifilterev & = \left[\begin{smallmatrix} \bifilterevL & \\ & \bifilterevH \end{smallmatrix}\right] .
\end{aligned}
\end{equation*}
Finally, the two coarse grid correction matrices are defined as
\begin{align}
\rcgc & = I-\rinterp\rcoarsesys^{-1}\rrestri \finesys \quad\mbox{and} \label{eq:rcgc_matrix} \\
\bcgc & = I-\binterp\bcoarsesys^{-1}\brestri \finesys \;. \label{eq:bcgc_matrix}
\end{align}
\end{defn}

Now, based on these definitions, the goal is to obtain an eigen--decomposition of the coarse grid correction matrices. The following lemma takes the first step by giving useful expressions for the inverse of the coarse system matrices in (\ref{eq:rcgc_matrix}) and (\ref{eq:bcgc_matrix}).

\begin{lem}[Navarrete and Coyle] \label{lem:galerkin_inverses}
In a red--black harmonic two--grid configuration the inverses of the red and black coarse grid matrices are given by
\begin{align}
\rcoarsesys^{-1} & = 4\; (\rdown W_L) \rdelta^{-1} (\rdown V_L)^H \quad\mbox{and} \label{eq:red_galerkin_inverse} \\
\bcoarsesys^{-1} & = 4\; (\bdown W_L) \bdelta^{-1} (\bdown V_L)^H \;, \label{eq:black_galerkin_inverse}
\end{align}
respectively, with
\begin{align}
\rdelta & = \rrfilterevL \Lambda_L \rifilterevL + \rrfilterevH \Lambda_H \rifilterevH \quad\mbox{and} \label{eq:red_delta} \\
\bdelta & = \brfilterevL \Lambda_L \bifilterevL + \brfilterevH \Lambda_H \bifilterevH \;. \label{eq:black_delta}
\end{align}
\end{lem}

The proof of this lemma is partially contained in \cite{PNavarrete_2008a} were only the red coarse grid was considered. The proof for the black grid is shown in Appendix \ref{app:proof_galerkin_inverses}.

The result in Lemma \ref{lem:galerkin_inverses} reflects the structure of a system with harmonic aliasing patterns. The eigenvectors of coarse system matrices are given by a linear--independent subset of the down--sampling eigenvectors of the system matrix. The coarse eigenvalues are expressed as sums of low-- and high--frequency eigenvalues as a result of aliasing effects.

Finally, the following theorem gives the eigen--decompositions of coarse grid correction matrices.

\begin{thm}[Navarrete and Coyle] \label{thm:cgc_decomposition}
In a red--black harmonic two--grid configuration the red and black coarse grid correction matrices (\ref{eq:rcgc_matrix}) and (\ref{eq:bcgc_matrix}) can be decomposed as
\begin{align}
\rcgc & = W \left[\begin{array}{cc}
\rmodal_{L\rightarrow L} & \rmodal_{H\rightarrow L} \\
\rmodal_{L\rightarrow H} & \rmodal_{H\rightarrow H}
\end{array} \right] V^H \quad\mbox{and} \label{eq:rcgc_decomposition} \\
\bcgc & = W \left[\begin{array}{cc}
\bmodal_{L\rightarrow L} & \bmodal_{H\rightarrow L} \\
\bmodal_{L\rightarrow H} & \bmodal_{H\rightarrow H}
\end{array} \right] V^H \;, \label{eq:bcgc_decomposition}
\end{align}
with
\begin{equation*}
\begin{aligned}
\rmodal_{L\rightarrow L} & = I - \rifilterevL \rdelta^{-1} \rrfilterevL \Lambda_L \;, & \rmodal_{H\rightarrow L} & = -\rifilterevL \rdelta^{-1} \rrfilterevH \Lambda_H \;, \\
\rmodal_{L\rightarrow H} & = -\rifilterevH \rdelta^{-1} \rrfilterevL \Lambda_L \;, & \rmodal_{H\rightarrow H} & = I - \rifilterevH \rdelta^{-1} \rrfilterevH \Lambda_H
\end{aligned}
\end{equation*}
and
\begin{equation*}
\begin{aligned}
\bmodal_{L\rightarrow L} & = I - \bifilterevL \bdelta^{-1} \brfilterevL \Lambda_L \;, & \bmodal_{H\rightarrow L} & = \bifilterevL \bdelta^{-1} \brfilterevH \Lambda_H \;, \\
\bmodal_{L\rightarrow H} & = \bifilterevH \bdelta^{-1} \brfilterevL \Lambda_L \;, & \bmodal_{H\rightarrow H} & = I - \bifilterevH \bdelta^{-1} \brfilterevH \Lambda_H \;.
\end{aligned}
\end{equation*}
\end{thm}

The proof of this theorem is partially contained in \cite{PNavarrete_2008a} were only the red coarse grid was considered. The result for the black coarse grid is shown in Appendix \ref{app:proof_cgc_decomposition}.

The results for the red and black coarse grids carry the difference in sign from the definition of harmonic aliasing patterns, which can be seen in the cross--modal symbols ($H\rightarrow L$ and $L\rightarrow H$).

%%%%%%%%%%%%%%%%%%%%%%%%%%%%%%%%%%%%%%%%%%%%%%%%%%%%
\section{Direct Two--grid Methods}
\label{sec:d2g}

In this section the full two--grid algorithm shown in Fig. \ref{fig:f2g} will be modified to obtain direct two--grid solvers for systems with harmonic aliasing patterns. The motivation is to eliminate smoothing iterations in the full two--grid scheme and base the algorithm purely on nested iterations and/or correction schemes.

A mere elimination of smoothing iterations in Fig. \ref{fig:f2g} would make it impossible for the algorithm to converge since partial information in a single coarse grid is not enough to get all the information from the fine grid. The algebra is very clear on this point because, based on the Galerkin condition, a coarse grid correction matrix is idempotent (or projection matrix). This is,
\begin{equation}
\rcgc^2=\rcgc \,\quad\mbox{and}\quad \bcgc^2=\bcgc \;,
\end{equation}
which means that several iterations of two--grid steps with a single coarse grid do nothing more than a single iteration. On the other hand, the red--black partition keeps all the information from the fine grid. Therefore, combining red and black coarse grids at different steps of the algorithm has a chance to converge depending on the configuration of the algorithm.

%---------------------------------------------------
\subsection{Multiplicative Approach}
\label{ssec:d2g_multiplicative}

Two modifications of the full two--grid algorithm from Fig. \ref{fig:f2g} are considered in Fig. \ref{fig:d2g_multiplicative}. First, smoothing iterations are removed. And second, red and black coarse grids are considered in the nested iteration and correction scheme steps, respectively.

If this algorithm works as a direct solver then it means that the problem is being factorized into coarse grid sub--problems. The following proposition shows this in terms of a factorization of $\finesys^{-1}$.

\begin{pro} \label{thm:multiplicative_decomposition}
The system shown in Fig. \ref{fig:d2g_multiplicative} works as a direct solver; i.e., $v=u$, if and only if the following decomposition applies on the inverse of the system matrix:
\begin{equation}
\finesys^{-1} = \rinterp \rcoarsesys^{-1} \rrestri + \binterp \bcoarsesys^{-1} \brestri - \rinterp \rcoarsesys^{-1} (\rrestri \finesys \binterp) \bcoarsesys^{-1} \brestri \;. \label{eq:multiplicative_decomposition}
\end{equation}
\end{pro}
\begin{proof}
Using (\ref{eq:two_grid_error}) in Fig. \ref{fig:d2g_multiplicative} gives
\begin{equation} \label{eq:error_multiplicative}
e = \bcgc \rcgc u \;.
\end{equation}
Then, an exact solution is obtained if and only if $\bcgc\rcgc=0$. Using (\ref{eq:rcgc_matrix}) and (\ref{eq:bcgc_matrix}) gives (\ref{eq:multiplicative_decomposition}).
\end{proof}

In this matrix factorization, the first two terms indicate the components of the inverse that come from the red and coarse grids independently. The third term indicates the dependence between the two solutions. In fact, the correction scheme in Fig. \ref{fig:d2g_multiplicative} works on top of the solution given by nested iteration and thus mixes the solutions from both coarse grids.

This factorization is known in domain decompositions as the \emph{multiplicative Schwartz procedure} \cite{AToselli_OWidlund_2004a}. In this procedure the evolution of the approximation error in $m$ iteartions is represented by the \emph{multiplicative operator}, $P_{mu}=I-E_{mu}$. Here, $E_{mu}=K_m\cdots K_1 K_0$ and $K_i$ are coarse grid correction matrices. The classical domain decomposition approach does not work as a direct solver and therefore $P_{mu}\neq I$. Thus, the two--grid configuration in Fig. \ref{fig:d2g_multiplicative} corresponds to a direct multiplicative Schwartz configuration.

The following theorem establishes conditions on the inter--grid filters to obtain a direct solver.

\begin{thm} \label{thm:multiplicative_qmf_solution}
A red--black harmonic two--grid configuration arranged as shown in Fig. \ref{fig:d2g_multiplicative}, with non--singular $\finesys$, $\rcoarsesys$, $\rrfilter$, $\bcoarsesys$ and $\bifilter$, works as a direct solver; i.e., $v=u$, if and only if
\begin{equation} \label{eq:multiplicative_condition}
\rrfilterevL \Lambda_L \bifilterevL = \rrfilterevH \Lambda_H \bifilterevH \;.
\end{equation}
\end{thm}
\begin{proof}
From (\ref{eq:error_multiplicative}), the direct solution is obtained if and only if
\begin{equation}
\bcgc \rcgc = \left[\begin{smallmatrix} 0 & 0 \\ 0 & 0 \end{smallmatrix}\right] \;.
\end{equation}
Using Theorem \ref{thm:cgc_decomposition}, the conditions in the diagonal blocks, which guarantee a perfect recovery of the solution, give
\begin{align}
(I & - \bifilterevL \bdelta^{-1} \brfilterevL \Lambda_L) (I - \rifilterevL \rdelta^{-1} \rrfilterevL \Lambda_L) & \nonumber \\
 & = \brfilterevH \brfilterevL (\rdelta\bdelta)^{-1} \bifilterevL \rifilterevH \Lambda_L \Lambda_H  & \mbox{and} \\
(I & - \bifilterevH \bdelta^{-1} \brfilterevH \Lambda_H) (I - \rifilterevH \rdelta^{-1} \rrfilterevH \Lambda_H) & \nonumber \\
 & = \brfilterevL \rrfilterevH (\rdelta\bdelta)^{-1} \bifilterevH \rifilterevL \Lambda_L \Lambda_H \;. &
\end{align}
The conditions in the off--diagonal blocks, which take care of aliasing cancellation, give
\begin{align}
(I & - \bifilterevL \bdelta^{-1} \brfilterevL \Lambda_L) \rifilterevL \rdelta^{-1} \rrfilterevH \Lambda_H & \nonumber \\
 & = (I - \rifilterevH \rdelta^{-1} \rrfilterevH \Lambda_H) \bifilterevL \bdelta^{-1} \brfilterevH  \Lambda_H & \mbox{and}\\
(I & - \bifilterevH \bdelta^{-1} \brfilterevH \Lambda_H) \rifilterevH \rdelta^{-1} \rrfilterevL \Lambda_L & \nonumber \\
 & = (I - \rifilterevL \rdelta^{-1} \rrfilterevL \Lambda_L) \bifilterevH \bdelta^{-1} \brfilterevL \Lambda_L \;. &
\end{align}
Here, all matrices are diagonal and therefore their products commute. None of the eigenvalues of inter--grid filters are zero because they are non--singular. Then, it is safe to multiply by any inverse of these matrices if necessary for simplifications. Since the coarse grid matrices are non--singular, each of the conditions above can be multiplied by $\rdelta\bdelta$. After this multiplication, using (\ref{eq:red_delta}) and (\ref{eq:black_delta}), the simplifications for the four equations independently give the same condition shown in (\ref{eq:multiplicative_condition}).
\end{proof}

It is clear from (\ref{eq:multiplicative_condition}) that inter--grid filters satisfying this condition should swap the low-- and high--frequency eigenvalues of the system matrix. The mirror of the system matrix fits perfectly for this purpose. The following corollary gives a particular solution which tries to keep the algorithm as simple as possible.

\begin{cor} \label{thm:multiplicative_qmf_configuration}
A red--black harmonic two--grid configuration arranged as shown in Fig. \ref{fig:d2g_multiplicative}, with non--singular $A$ and such that $\det(\Lambda_L+\Lambda_H)\neq 0$, works as a direct solver with the following configuration:
\begin{equation}
\begin{aligned}
\rifilter & = I \;,          & \quad           & \quad & \rrfilter & = I \;, \\
\bifilter & = \finesys^\star & \mbox{and}\quad & \quad & \brfilter & = I \;.
\end{aligned}
\end{equation}
The coarse grid correction matrices for this configuration are
\begin{align}
\rcgc & = (\Lambda_L+\Lambda_H)^{-1}\; W \left[\begin{array}{cc}
\Lambda_H & -\Lambda_H \\
-\Lambda_L & \Lambda_L
\end{array} \right] V^T  \quad\mbox{and} \label{eq:black_multiplicative_cgc} \\
\bcgc & = \frac{1}{2}\; W \left[\begin{array}{cc}
I & I \\
I & I
\end{array} \right] V^T \;. \label{eq:red_multiplicative_cgc}
\end{align}
\end{cor}
\begin{proof}
Using the eigen--decomposition of filters and Theorem \ref{thm:frequency_inversion} gives
\begin{equation} \label{eq:multiplicative_qmf_eigenvalues}
\begin{aligned}
\rifilterev & = I \;,                                     &            &       & \rrfilterev & = I \;, \\
\bifilterevL & = \Lambda_H \;,\; \bifilterevH = \Lambda_L & \mbox{and} & \quad & \brfilterev & = I \;.
\end{aligned}
\end{equation}
Using these eigenvalues in (\ref{eq:red_delta}) and (\ref{eq:black_delta}) gives $\rdelta=\Lambda_L+\Lambda_H$ and $\bdelta=2\Lambda_L\Lambda_H$. Therefore, by Lemma \ref{lem:galerkin_inverses} the coarse grid system matrices are invertible. Then, Theorem \ref{thm:multiplicative_qmf_solution} can be applied and the eigenvalues of inter--grid filters above fulfill the condition (\ref{eq:multiplicative_condition}). Finally, using (\ref{eq:multiplicative_qmf_eigenvalues}) in Theorem \ref{thm:cgc_decomposition} gives (\ref{eq:black_multiplicative_cgc}) and (\ref{eq:red_multiplicative_cgc}).
\end{proof}

From (\ref{eq:red_multiplicative_cgc}) we see that the nested iteration step is not filtering nor reducing any frequency component of the error. It just equals the gain of all the effects. On the other hand, in (\ref{eq:black_multiplicative_cgc}) we see that the correction scheme acts as a mirror filter by swapping the low-- and high--frequency eigenvalues of the system matrix in the diagonal blocks. The aliasing effect in the off--diagonals is adjusted to cancel the symbols at the same row in $\bcgc$, so that $\bcgc\rcgc=0$.

This solution is particularly simple in the nested iteration step, since only down/up--sampling operations are used. The coarse grid matrix $\rcoarsesys=\rdown\finesys\rup$ has better sparseness than the system matrix $\finesys$. The sparseness of $\bcoarsesys=\bdown\finesys^\star\finesys\bup$ depends on the structure of down--sampling and non--zeros in $\finesys$.

When solving Poisson's equation, this solution is equivalent to the total reduction method by Schr{\"o}der and Trottenberg \cite{JSchroder_UTrottenberg_1973a,JSchroder_UTrottenberg_HReutersberg_1976a}. The derivation of total reduction methods was based on the structure of a stationary impulse response for LSI systems and thus imposes stronger assumptions on the system.

%---------------------------------------------------
\subsection{Additive Approach}
\label{ssec:d2g_additive}

\begin{figure*}[!t]
\centerline{
\subfloat[Multiplicative scheme]{{\label{fig:d2g_multiplicative}\includegraphics[width=4.2in]{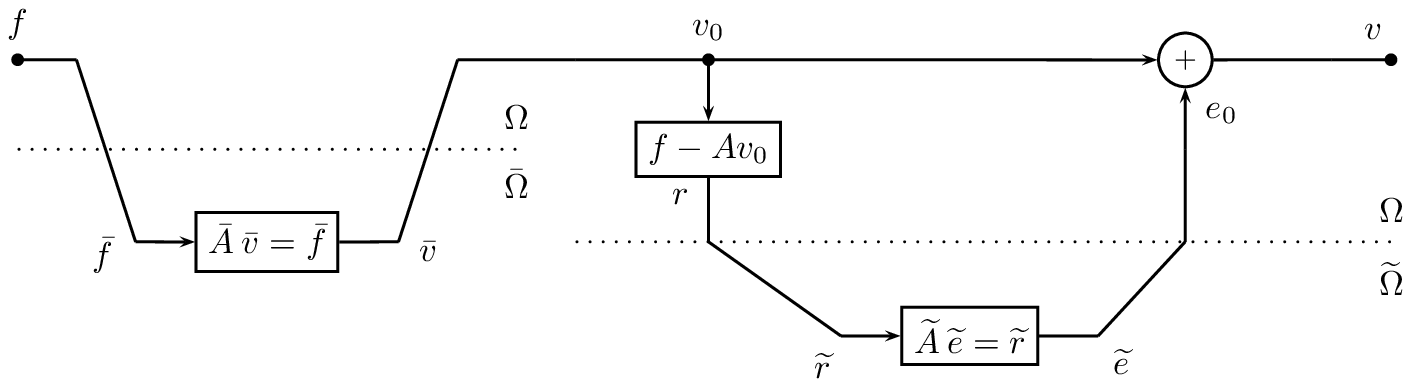}}} \hfil
\subfloat[Additive scheme]{{\label{fig:d2g_additive}\includegraphics[width=2.2in]{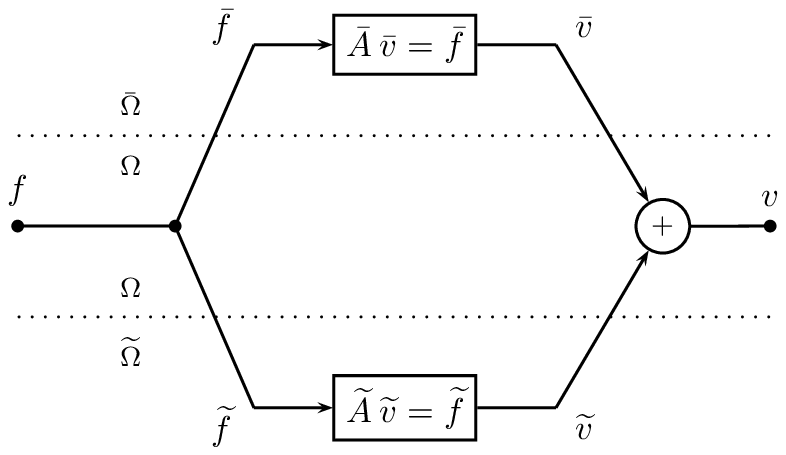}}}}
\caption{Two--grid schemes approaches to solve $\finesys u=f$ exactly. The multiplicative scheme is different than the full two--grid scheme shown in Fig. \ref{fig:f2g} because no smoothing iterations are used, and different coarse grids are used at each step. The additive scheme is the two--grid algorithm's version of the two--channel multi--rate system in Fig. \ref{fig:multirate_2ch_finite}. In both cases the two coarse grids form a partition of the fine grid and thus capture all the information from the fine grid. This fact allows these schemes to work as direct solvers; i.e., allows $v=u$.} \label{fig:d2g}
\end{figure*}

The second approach is the two--grid algorithm shown in Fig. \ref{fig:d2g_additive} which uses the structure of the two--channel multi--rate system from Fig. \ref{fig:multirate_2ch_finite}. This is a two--grid scheme with two nested iterations working in parallel at red and black coarse grids. The red and black coarse grids work as the two channels of a multi--rate system, but now, linear systems of equations are solved before the interpolation and addition of the two approximations. In a multi--rate system the idea is to reconstruct the input signal (in this case $f$) but now the idea is to transform this signal into the solution of the linear system.

If this algorithm works as a direct solver then it means that the problem is being factorized into coarse grid sub--problems. The following proposition shows this in terms of a factorization of $\finesys^{-1}$.

\begin{pro} \label{thm:additive_decomposition}
The system shown in Fig. \ref{fig:d2g_additive} works as a direct solver; i.e., $v=u$, if and only if the following decomposition applies on the inverse of the system matrix:
\begin{equation}
\finesys^{-1} = \rinterp \rcoarsesys^{-1} \rrestri + \binterp \bcoarsesys^{-1} \brestri \;. \label{eq:additive_decomposition}
\end{equation}
\end{pro}
\begin{proof}
The approximation error in Fig. \ref{fig:d2g_additive} is given by $e=u-(\bar{v}+\widetilde{v})$. Using (\ref{eq:two_grid_error}) gives
\begin{equation} \label{eq:error_additive}
e = \left(\rcgc + \bcgc - I\right) u
\end{equation}
Then, an exact solution is obtained if and only if $\rcgc + \bcgc = I$. Using (\ref{eq:rcgc_matrix}) and (\ref{eq:bcgc_matrix}) gives (\ref{eq:additive_decomposition}).
\end{proof}

As opposed to (\ref{eq:multiplicative_decomposition}), here no cross--terms appear in the decomposition. This reflects the fact that nested iterations run independent of each other. This makes this scheme better suited for parallelization.

This factorization is known in domain decompositions as the \emph{additive Schwartz procedure} \cite{AToselli_OWidlund_2004a}. In this procedure the evolution of the approximation error in $m$ iteartions is represented by the \emph{additive operator}, $P_{ad}=K_m + \cdots + K_1 + K_0$, where $K_i$ are coarse grid correction matrices. The classical domain decomposition approach does not work as a direct solver and therefore $P_{ad}\neq I$. Thus, the two--grid configuration in Fig. \ref{fig:d2g_additive} corresponds to a direct additive Schwartz configuration.

The following theorem establishes conditions on the inter--grid filters to obtain a direct solver.

\begin{thm} \label{thm:additive_qmf_solution}
A red--black harmonic two--grid configuration arranged as shown in Fig. \ref{fig:d2g_additive}, with non--singular $\finesys$, $\rcoarsesys$, $\rrfilter$, $\bcoarsesys$ and $\bifilter$, works as a direct solver; i.e., $v=u$, if and only if
\begin{align}
\rrfilterevL & \brfilterevL \Lambda_L^2 \rifilterevL \bifilterevL = \rrfilterevH \brfilterevH \Lambda_H^2 \rifilterevH \bifilterevH \;, \label{eq:additive_condition_1} \\
\rrfilterevL & \brfilterevH \Lambda_L^2 \rifilterevL \bifilterevL + \rrfilterevH \brfilterevH \Lambda_H^2 \rifilterevH \bifilterevL = \nonumber \\
 & \rrfilterevH \brfilterevL \Lambda_L^2 \rifilterevL \bifilterevL + \rrfilterevH \brfilterevH \Lambda_H^2 \rifilterevL \bifilterevH \quad\mbox{and} \label{eq:additive_condition_2} \\
\rrfilterevL & \brfilterevL \Lambda_L^2 \rifilterevL \bifilterevH + \rrfilterevH \brfilterevL \Lambda_H^2 \rifilterevH \bifilterevH = \nonumber \\
 & \rrfilterevL \brfilterevL \Lambda_L^2 \rifilterevH \bifilterevL + \rrfilterevL \brfilterevH \Lambda_H^2 \rifilterevH \bifilterevH \;. \label{eq:additive_condition_3}
\end{align}
\end{thm}
\begin{proof}
From (\ref{eq:error_additive}) the direct solution is obtained if and only if
\begin{equation}
\rcgc + \bcgc = \left[\begin{smallmatrix} I & 0 \\ 0 & I \end{smallmatrix}\right] \;.
\end{equation}
Using Theorem \ref{thm:cgc_decomposition}, the conditions in the diagonal blocks, which guarantee a perfect recovery of the solution, give
\begin{align}
\rifilterevL \rdelta^{-1} \rrfilterevL \Lambda_L + \bifilterevL \bdelta^{-1} \brfilterevL \Lambda_L & = I \quad\mbox{and} \label{eq:additive_condition_proof_1}\\
\rifilterevH \rdelta^{-1} \rrfilterevH \Lambda_H + \bifilterevH \bdelta^{-1} \brfilterevH \Lambda_H & = I \;. \label{eq:additive_condition_proof_2}
\end{align}
The conditions in the off--diagonal blocks, which take care of aliasing cancellation, give
\begin{align}
\rifilterevL \rdelta^{-1} \rrfilterevH - \bifilterevL \bdelta^{-1} \brfilterevH & = 0 \quad\mbox{and} \label{eq:additive_condition_proof_3} \\
\rifilterevH \rdelta^{-1} \rrfilterevL - \bifilterevH \bdelta^{-1} \brfilterevL & = 0 \;. \label{eq:additive_condition_proof_4}
\end{align}
Here, all matrices are diagonal and therefore their products commute. None of the eigenvalues of inter--grid filters are zero because they are non--singular. Then, it is safe to multiply by any inverse of these matrices if necessary for simplifications. Since the coarse grid matrices are non--singular, each of the conditions above can be multiplied by $\rdelta\bdelta$. After this multiplication, using (\ref{eq:red_delta}) and (\ref{eq:black_delta}), the algebra on the conditions (\ref{eq:additive_condition_proof_1}) and (\ref{eq:additive_condition_proof_2}) simplifies to the same condition in (\ref{eq:additive_condition_1}), and the algebra on (\ref{eq:additive_condition_proof_3}) and (\ref{eq:additive_condition_proof_4}) simplifies to (\ref{eq:additive_condition_2}) and (\ref{eq:additive_condition_3}), respectively.
\end{proof}

These conditions are analogous to Vetterli's conditions (\ref{eq:multirate_2ch_finite_perfect_reconstruction}--\ref{eq:multirate_2ch_finite_aliasing_freeHL}). The analogy is more clear from (\ref{eq:additive_condition_proof_1}--\ref{eq:additive_condition_proof_4}) where the only difference with (\ref{eq:multirate_2ch_finite_perfect_reconstruction}--\ref{eq:multirate_2ch_finite_aliasing_freeHL}) is the existence of the matrices $\rdelta$, $\bdelta$, $\Lambda_L$ and $\Lambda_H$. This indicates the fact that a linear system of equations is being solved.

Again, the mirror of the system matrix can be used to find a solution of these conditions. The following corollary gives a particular solution which tries to keep the algorithm as simple as possible.

\begin{cor} \label{thm:additive_qmf_configuration}
A red--black harmonic two--grid configuration arranged as shown in Fig. \ref{fig:d2g_additive}, with non--singular $A$, works as a direct solver with the following configuration:
\begin{equation}
\begin{aligned}
\rifilter & = \finesys^\star \;, & \quad           & \quad & \rrfilter = I \;, \\
\bifilter & = \finesys^\star     & \mbox{and}\quad & \quad & \brfilter = I \;.
\end{aligned}
\end{equation}
The coarse grid correction matrices for this configuration are
\begin{align}
\rcgc & = \frac{1}{2} \: W
\left[\begin{array}{cc}
I & -I \\
-I & I
\end{array} \right] V^H \quad\mbox{and} \label{eq:red_additive_cgc} \\
\bcgc & = \frac{1}{2} \: W
\left[\begin{array}{cc}
I & I \\
I & I
\end{array} \right] V^H \;. \label{eq:black_additive_cgc}
\end{align}
\end{cor}
\begin{proof}
Using the eigen--decomposition of filters and Theorem \ref{thm:frequency_inversion} gives
\begin{equation} \label{eq:additive_qmf_eigenvalues}
\begin{aligned}
\rifilterevL & = \Lambda_H \;, & \rifilterevH & = \Lambda_L \;, &            & \quad & \rrfilterev & = I \;, \\
\bifilterevL & = \Lambda_H \;, & \bifilterevH & = \Lambda_L     & \mbox{and} & \quad & \brfilterev & = I \;.
\end{aligned}
\end{equation}
Using these eigenvalues in (\ref{eq:red_delta}) and (\ref{eq:black_delta}) gives $\rdelta=\bdelta=2\Lambda_L\Lambda_H$. Therefore, by Lemma \ref{lem:galerkin_inverses}, the coarse grid system matrices are invertible. Then, Theorem \ref{thm:additive_qmf_solution} can be applied and the eigenvalues of inter--grid filters fulfill conditions (\ref{eq:additive_condition_1}--\ref{eq:additive_condition_3}). Using the eigenvalues from (\ref{eq:additive_qmf_eigenvalues}) in Theorem \ref{thm:cgc_decomposition} gives (\ref{eq:red_additive_cgc}) and (\ref{eq:black_additive_cgc}).
\end{proof}

The coarse grid correction matrices in (\ref{eq:red_additive_cgc}) and (\ref{eq:black_additive_cgc}) equal the gain of filtering and aliasing effects. Interestingly, the symbols of $\rcgc$ and $\bcgc$ correspond to the black and red harmonic aliasing pattern in (\ref{eq:redblack_aliasing_patterns}), $\bhap$ and $\rhap$, respectively. The aliasing effect in the off--diagonals have opposed signs between red and black coarse grids, so that $\bcgc+\rcgc=I$.

Compared with the solution for the multiplicative approach, the additive approach involves more computations. This is because both red and black coarse grid matrices use a mirror filter. On the other hand, this approach is better suited for parallelization. Therefore, both the multiplicative an additive approaches become useful depending on the computational resources available. The multiplicative approach is more convenient in a single processor and the additive approach is more convenient in a multi--core architecture.

%%%%%%%%%%%%%%%%%%%%%%%%%%%%%%%%%%%%%%%%%%%%%%%%%%%%
\section{Direct Multi--grid Methods}
\label{sec:dmg}

The convergence analysis of previous sections will be valid at each coarse level if harmonic aliasing patterns exist for each coarse system matrix. The following definition introduces the multi--scale property needed on the biorthogonal eigenvectors of the system.

\begin{defn}[Multi--grid Harmonic Basis]
A \emph{multi--grid harmonic basis} of level $0$ is any biorthogonal basis of $\mathbb{C}^n$. A \emph{multi--grid harmonic basis} of level $l>0$ is a biorthogonal basis of $\mathbb{C}^n$, $n=2^l n_0$ and $n_0\in\mathbb{N}^+$, with red--black harmonic aliasing patterns and such that the down--sampled low--frequency eigenvectors form a \emph{multi--grid harmonic basis} of level $l-1$.
\end{defn}

The two--grid methods in section \ref{sec:d2g} can be extended to multiple grids if the eigenvectors of a system matrix $\finesys\in\mathbb{C}^{n\times n}$, with $n=2^l n_0$ and $n_0=\mathcal{O}(1)$, form a multi--grid harmonic basis of level $l$.

%---------------------------------------------------
\subsection{Multiplicative direct multi--grid algorithm}
\label{ssec:dmg_multiplicative}

\begin{table*}
\renewcommand{\arraystretch}{1.3}
\caption{Pseudocode for the direct multi--grid algorithm following multiplicative and additive approaches to solve $\finesys u=f$. The matrix $\finesys^\star$ represents the mirror of $\finesys$ according to definition \ref{def:mirror_operator}, which only involves changing the sign of some of the entries in $\finesys$. It follows from both cases that the number of multiplications is $m(n)=\mathcal{O}(n \log n)$.} \label{tab:dmg}
\centering
\begin{tabular}{lc}
Task & Multiplications\\ \hline
$u = \mathtt{DMG\_multiplicative}(\finesys,f)$ & $m(n)$ \\
\hspace{.2in} $\mathtt{if}$ $n>n_0\;,\; n_0=\mathcal{O}(1)$ & \\
\hspace{.4in} $\bar{f} = \rdown f$ & $0$ \\
\hspace{.4in} $\rcoarsesys = \rdown \finesys\rup$ & $0$ \\
\hspace{.4in} $\bar{v} = \mathtt{DMG\_multiplicative}(\rcoarsesys,\bar{f})$ & $m(\tfrac{n}{2})$ \\
\hspace{.4in} $v_0 = \rup \bar{v}$ & $0$ \\
\hspace{.4in} $r = f-\finesys v$ & $\mathcal{O}(n)$ \\
\hspace{.4in} $\widetilde{r} = \bdown r$ & $0$ \\
\hspace{.4in} $\binterp = \finesys^\star \bup$ & $0$ \\
\hspace{.4in} $\bcoarsesys = \bdown \finesys \binterp$ & $\mathcal{O}(n)$ \\
\hspace{.4in} $\widetilde{v} = \mathtt{DMG\_multiplicative}(\bcoarsesys,\widetilde{r})$ & $m(\tfrac{n}{2})$ \\
\hspace{.4in} $e_0 = \binterp \widetilde{v}$ & $\mathcal{O}(n)$ \\
\hspace{.4in} $v = v_0+e_0$ & $0$ \\
\hspace{.4in} $\mathtt{return}(v)$ & \\
\hspace{.2in} $\mathtt{else}$ & \\
\hspace{.4in} $\mathtt{return}(\finesys^{-1} f)$ & $\mathcal{O}(1)$
\end{tabular} \hfil
\begin{tabular}[hbt]{lc}
Task & Multiplications\\ \hline
$u = \mathtt{DMG\_additive}(\finesys,f)$ & $m(n)$ \\
\hspace{.2in} $\mathtt{if}$ $n>n_0\;,\; n_0=\mathcal{O}(1)$ & \\
\hspace{.4in} $\bar{f} = \rdown f$ & $0$ \\
\hspace{.4in} $\rinterp = \finesys^\star \rup$ & $0$ \\
\hspace{.4in} $\rcoarsesys = \rdown \finesys \rinterp$ & $\mathcal{O}(n)$ \\
\hspace{.4in} $\bar{v} = \mathtt{DMG\_additive}(\rcoarsesys,\bar{f})$ & $m(\tfrac{n}{2})$ \\
\hspace{.4in} $\binterp = \finesys^\star \bup$ & $0$ \\
\hspace{.4in} $\bcoarsesys = \bdown \finesys \binterp$ & $\mathcal{O}(n)$ \\
\hspace{.4in} $\widetilde{v} = \mathtt{DMG\_additive}(\bcoarsesys,\widetilde{f})$ & $m(\tfrac{n}{2})$ \\
\hspace{.4in} $u = \rinterp \bar{v} + \binterp \widetilde{v}$ & $\mathcal{O}(n)$ \\
\hspace{.4in} $\mathtt{return}(u)$ & \\
\hspace{.2in} $\mathtt{else}$ & \\
\hspace{.4in} $\mathtt{return}(\finesys^{-1} f)$ & $\mathcal{O}(1)$
\end{tabular}
\end{table*}

The recursive implementation of the multiplicative approach shown in Fig. \ref{fig:d2g_multiplicative} is explained in Table \ref{tab:dmg}. Here, the difference equation for the number of multiplications performed gives $m(n)=\mathcal{O}(n \log n)$.

Starting from two coarse grids in the first level, the algorithm creates coarser grids which form a partition of the fine grid into more subsets of nodes. In Fig. \ref{fig:dmg_multiplicative} the sequence of coarse grids visited by the algorithm is shown for $l=3$. The structure is a W--cycle, well known in multi--grid methods \cite{WLBriggs_VEHenson_SFMcCormick_2000a}. At each level the partition of nodes is duplicated. For instance, in the coarsest level the grid partition gives
\begin{equation} \label{eq:dmg_partition_level3}
\finegrid=\finegrid_{rrr}+\finegrid_{rrb}+\finegrid_{rbr}+\finegrid_{rbb}+\finegrid_{brr}+\finegrid_{brb}+\finegrid_{bbr}+\finegrid_{bbb} \;,
\end{equation}
where the subindexes from left (fine) to right (coarse) indicate if red or black grids were chosen.

This algorithm has the nice property that the red coarse system matrix $\rcoarsesys=\rdown\finesys\rup$ reduces the sparseness of the system matrix $\finesys$. In coarser levels the system matrix might soon become diagonal and the system is solved in linear time. Thus, there are good chances that the structure in Fig. \ref{fig:dmg_multiplicative} changes from a W--cycle into a V--cycle \cite{WLBriggs_VEHenson_SFMcCormick_2000a}, reducing the computational complexity from $\mathcal{O}(n \log n)$ to $\mathcal{O}(n)$. This is actually what happens when solving Poisson's equation, where this method is equivalent to total reduction methods \cite{JSchroder_UTrottenberg_1973a,JSchroder_UTrottenberg_HReutersberg_1976a}. In general, this depends on the structure of the system and it will not be study here in depth. The example in section \ref{sec:example} will show a case where the computational complexity is effectively reduced.

%---------------------------------------------------
\subsection{Additive direct multi--grid algorithm}
\label{ssec:dmg_additive}

The recursive implementation of the additive approach shown in Fig. \ref{fig:d2g_additive} is explained in Table \ref{tab:dmg}. Here, the difference equation for the number of multiplications performed gives $m(n)=\mathcal{O}(n \log n)$.

\begin{figure*}[!t]
\centerline{
\subfloat[Multiplicative algorithm]{\raisebox{0.2in}{\label{fig:dmg_multiplicative}\includegraphics[width=4.2in]{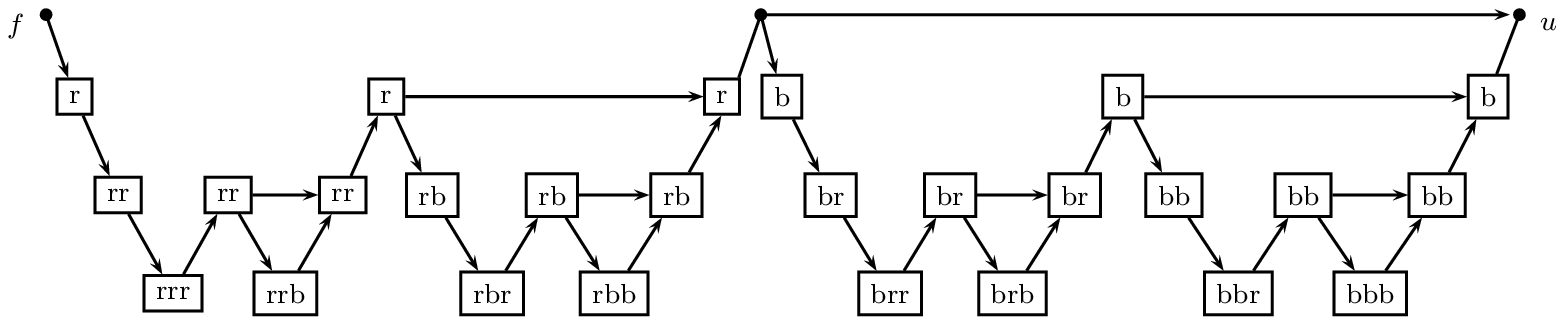}}} \hfil
\subfloat[Additive algorithm]{{\label{fig:dmg_additive}\includegraphics[width=3in]{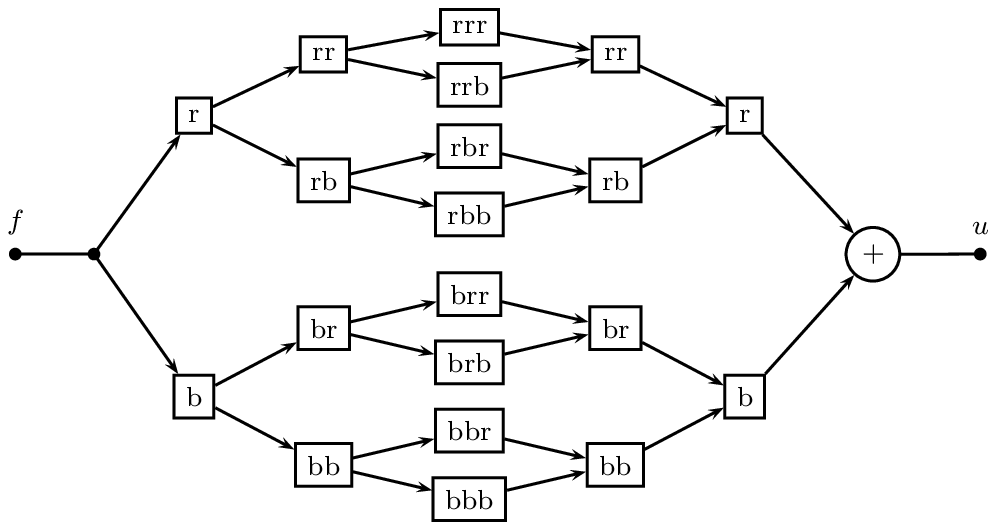}}}}
\caption{Diagram of the sequence of computations in a direct multi--grid scheme following multiplicative and additive approaches. Moving from fine to coarse levels, each box adds an `$r$' or `$b$' at the right to indicate whether the red or black coarse grid is used, respectively. In the multiplicative approach, all the possible coarse grids are visited sequentially through a W--cycle in order to obtain the exact solution. In the additive approach, all the possible coarse grids are visited in parallel through a binary tree in order to obtain the exact solution.} \label{fig:dmg}
\end{figure*}
Same as for the multiplicative approach, the algorithm creates coarser grids which form a partition of the fine grid into more subsets of nodes. In Fig. \ref{fig:dmg_additive} the sequence of coarse grids visited by the algorithm is shown for $l=3$. The structure is a binary tree that duplicates the number of partitions at each level. The partition of the coarsest grid is the same as in (\ref{eq:dmg_partition_level3}). Here, each of the coarsest grids is visited only once, which is possible because the solutions from different partitions do not depend on each other.

\begin{figure}[!t]
\centering
\includegraphics[width=3.5in]{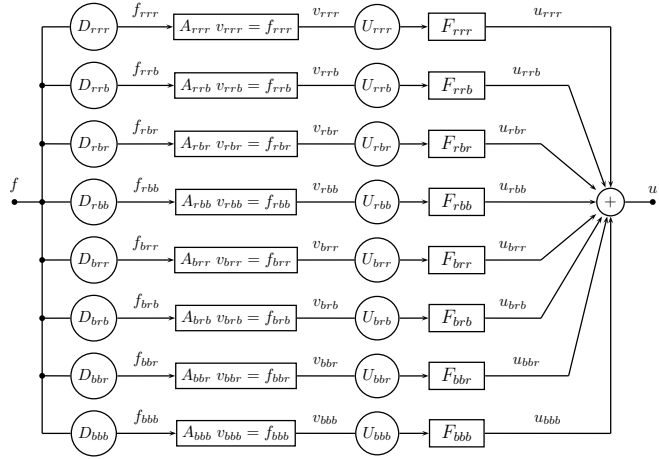}
\caption{Multi--grid scheme following a multi--channel additive approach to solve $\finesys u=f$. Compared with the binary tree path shown in Fig. \ref{fig:dmg_additive}, here the source vector $f$ is directly taken to the coarsest levels. The coarsest grids form a partition of the fine grid with the eight possible combinations of red ($r$) and black ($b$) coarsening.} \label{fig:dmg_multichannel}
\end{figure}
This algorithm is more convenient for parallelization. In fact, from the finest to the coarsest grid in Fig. \ref{fig:dmg_additive}, only down--sampling operations are performed. From the coarsest to the finest grid, only interpolations are performed. Therefore, these operations can be carried out in a single step and then the algorithm becomes analogous to a multi--rate system with multiple channels. The diagram of this algorithm is shown in Fig. \ref{fig:dmg_multichannel}. The computational complexity of a parallel implementation divides $\mathcal{O}(n \log n)$ by the number of channels and reaches $\mathcal{O}(\log n)$ in the best scenario.

%%%%%%%%%%%%%%%%%%%%%%%%%%%%%%%%%%%%%%%%%%%%%%%%%%%%
\section{Example}
\label{sec:example}

Consider a two--dimensional system in which $\finesys$ is the finite--difference discretization of Helmholtz's equation $-\nabla^2 u-k^2 u=0$ on a square domain with periodic boundary conditions and unit step size. The size of the problem is set to $n=32\times 32$ and the wavenumber is set to $k=\tfrac{\pi}{3}$. Thus, the stationary impulse response of $\finesys$ is $\left[\begin{smallmatrix} & -1 & \\-1 & \underline{4-k^2} & -1\\ & -1 & \end{smallmatrix}\right]$ (the underline denotes the diagonal element).

\begin{figure}[!t]
\centerline{
\subfloat[$\finegrid=\finegrid_{r}+\finegrid_{b}$]{\label{fig:e2_grid_1}
\frame{\includegraphics[width=1.5in]{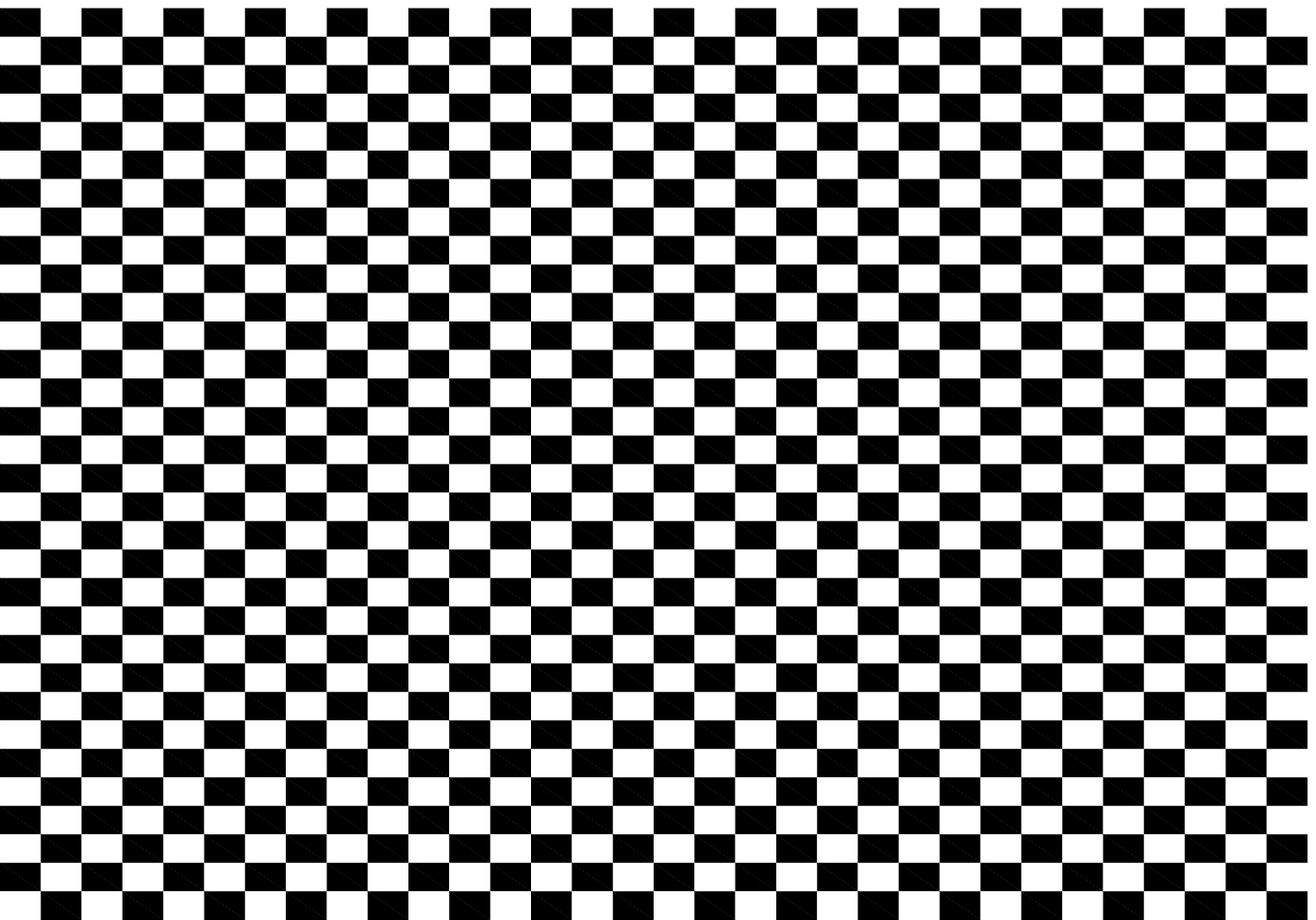}} \hfil
\frame{\includegraphics[width=1.5in]{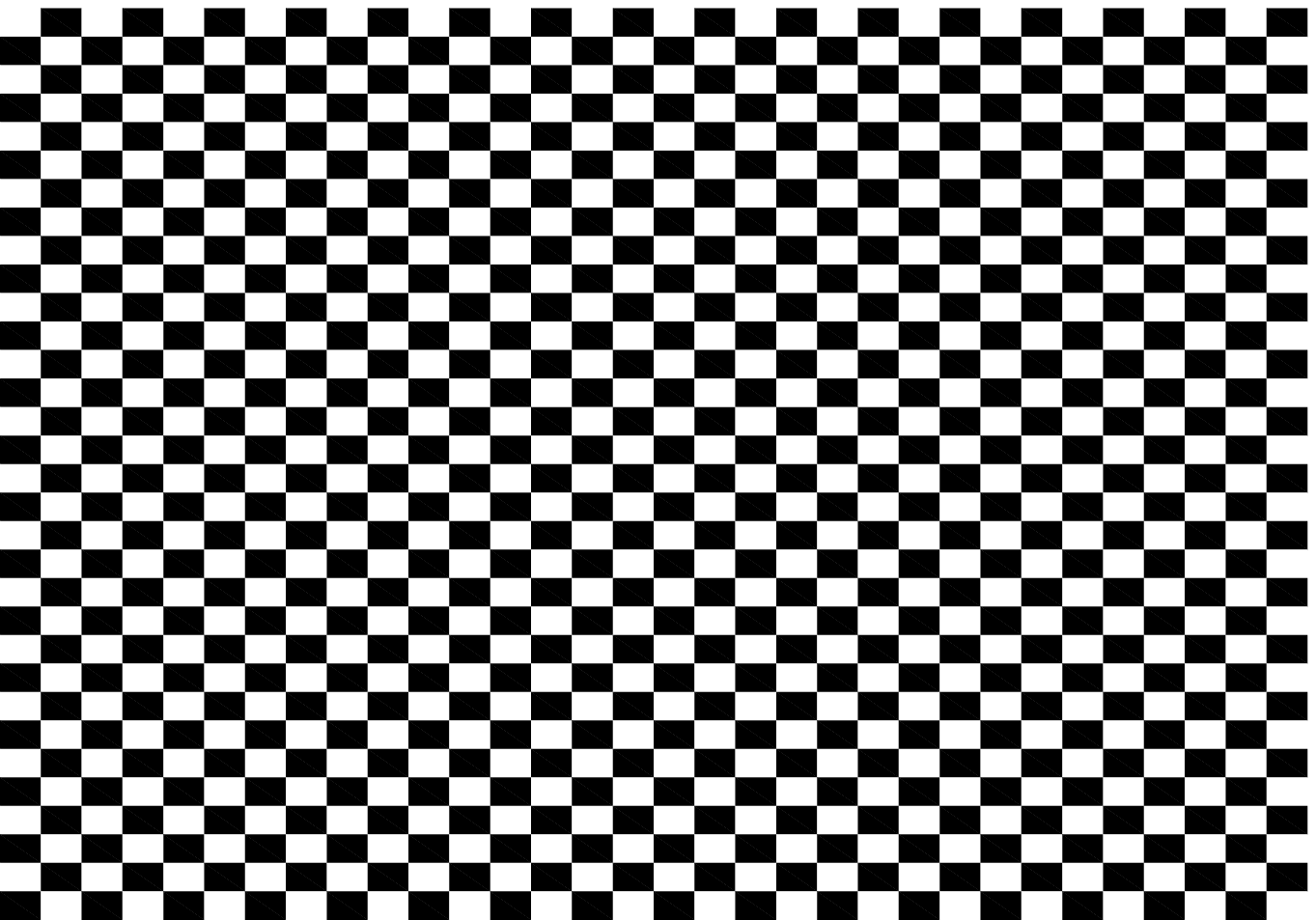}}}}
\centerline{
\subfloat[$\finegrid=\finegrid_{rr}+\finegrid_{rb}+\finegrid_{br}+\finegrid_{bb}$]{
\frame{\includegraphics[width=0.74in]{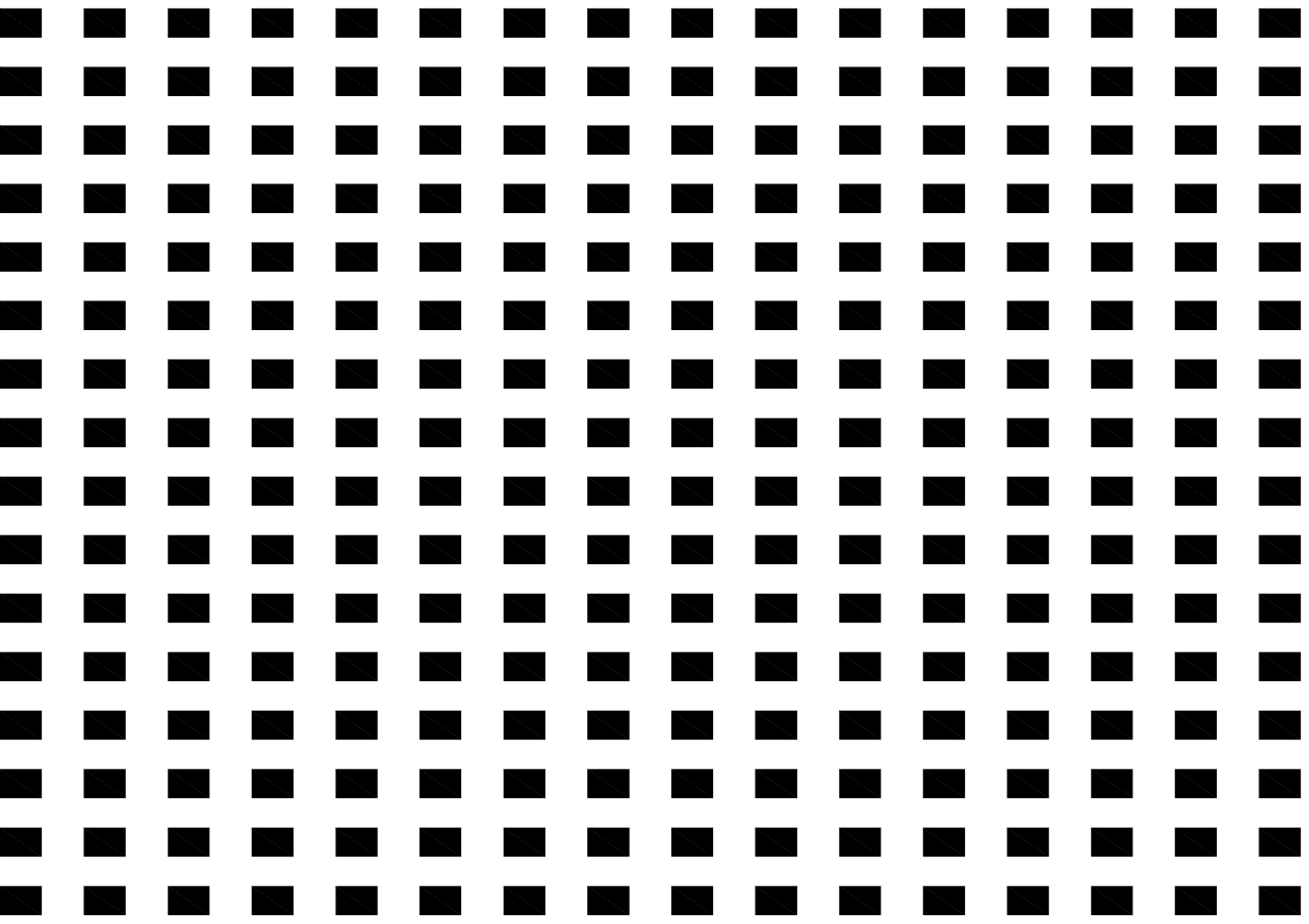}} \hfil
\frame{\includegraphics[width=0.74in]{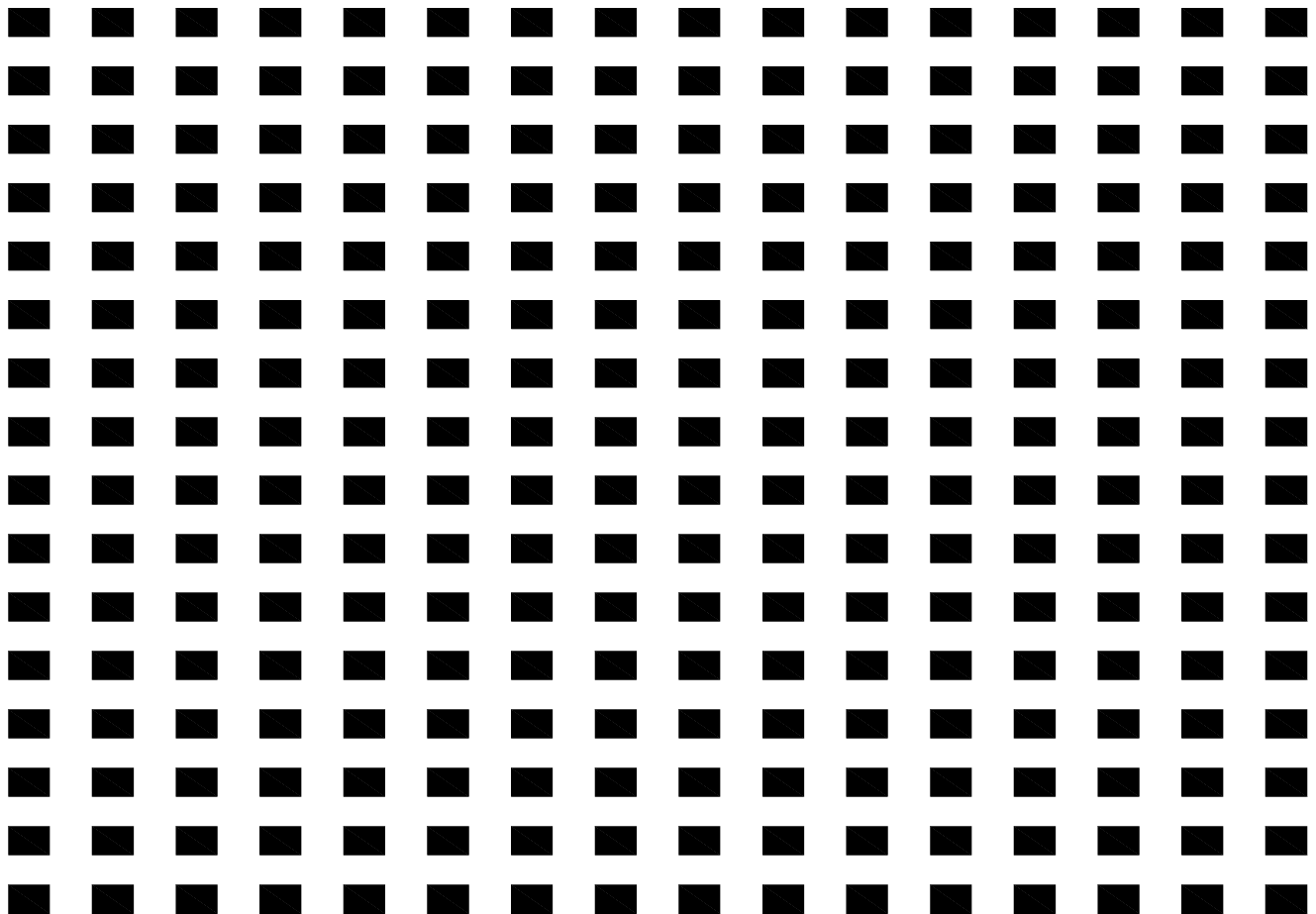}} \hfil
\frame{\includegraphics[width=0.74in]{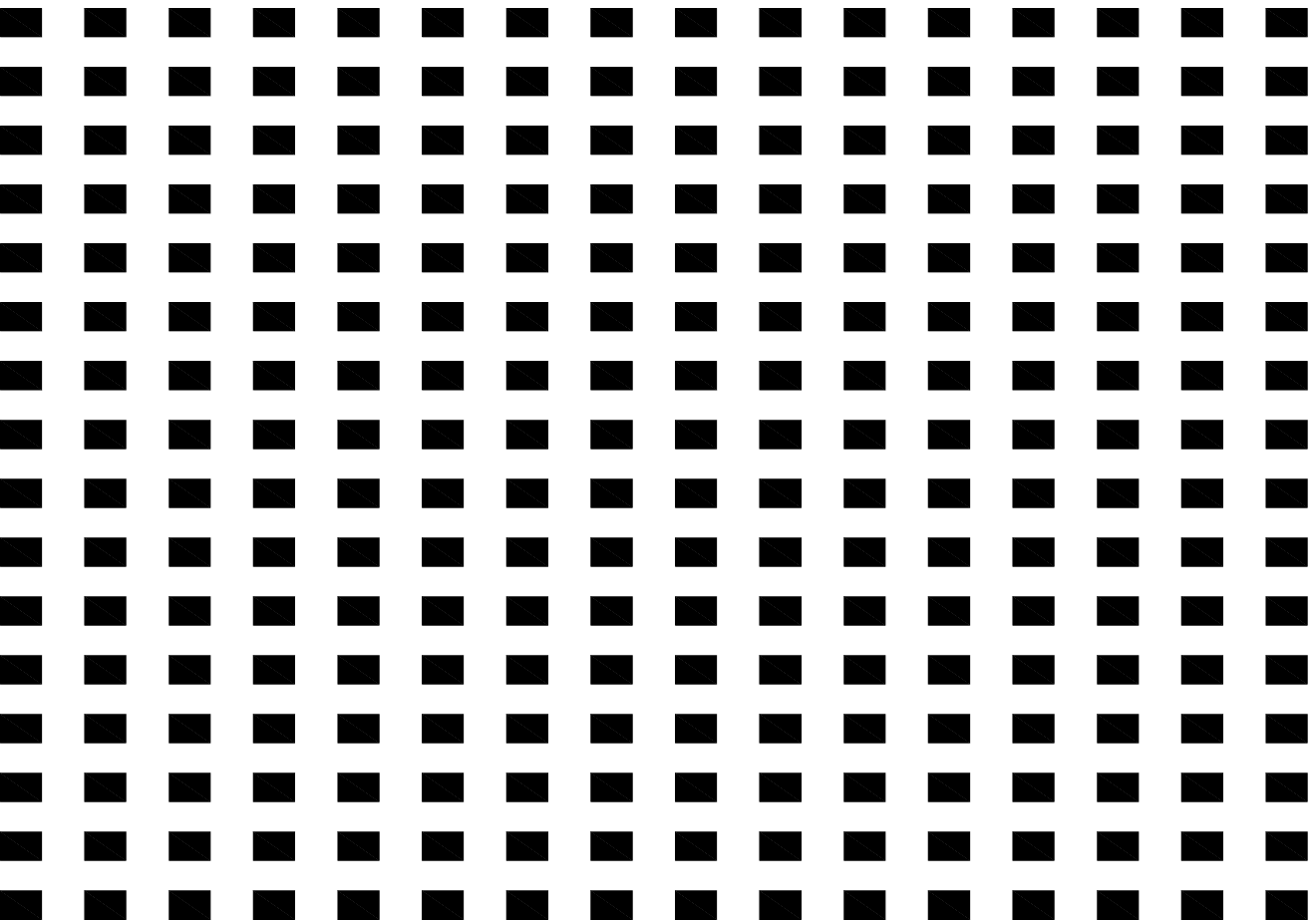}} \hfil
\frame{\includegraphics[width=0.74in]{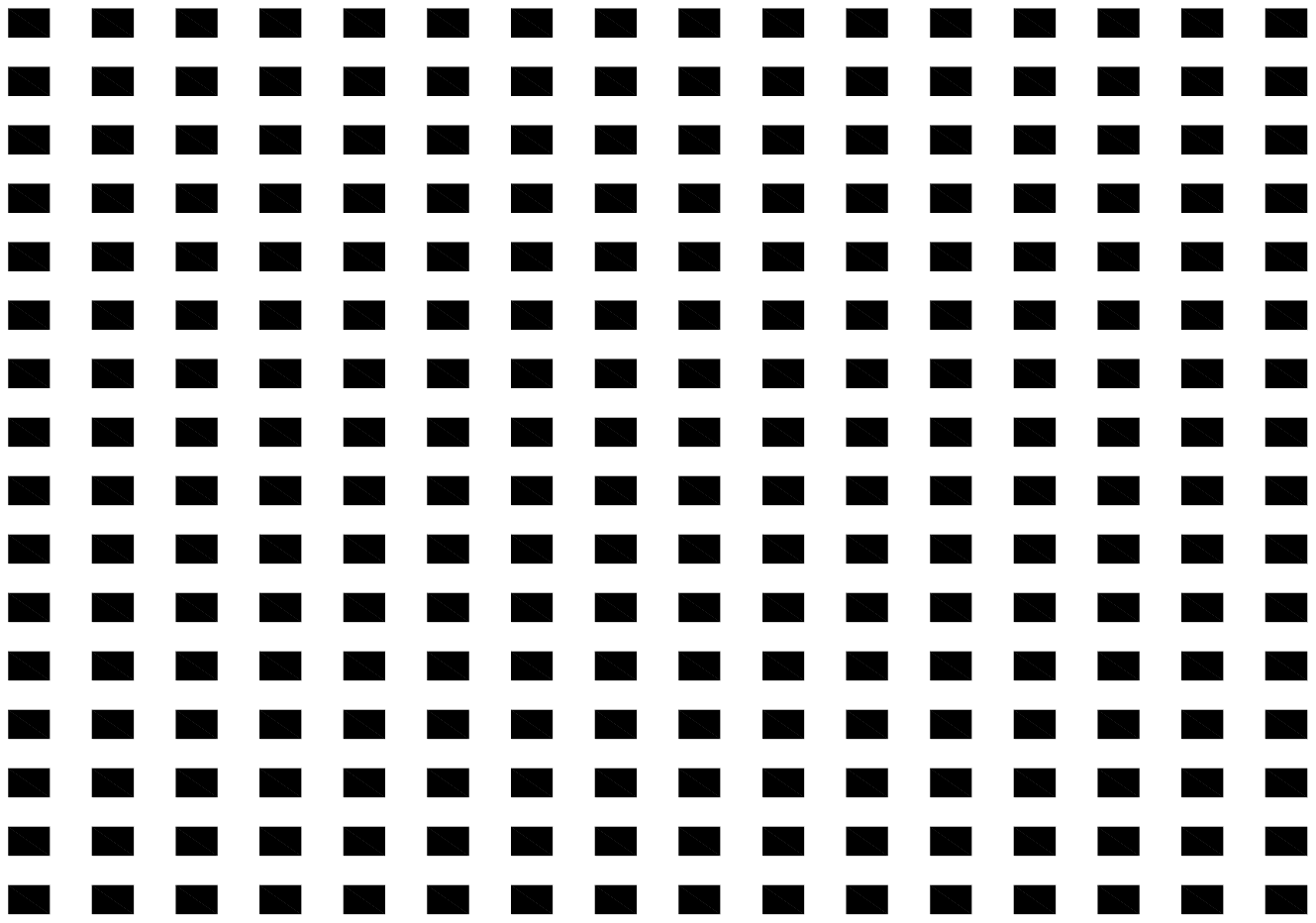}}}}
\centerline{
\subfloat[$\finegrid=\finegrid_{rrr}+\finegrid_{rrb}+\finegrid_{rbr}+\finegrid_{rbb}+\finegrid_{brr}+\finegrid_{brb}+\finegrid_{bbr}+\finegrid_{bbb}$]{
\frame{\includegraphics[width=0.34in]{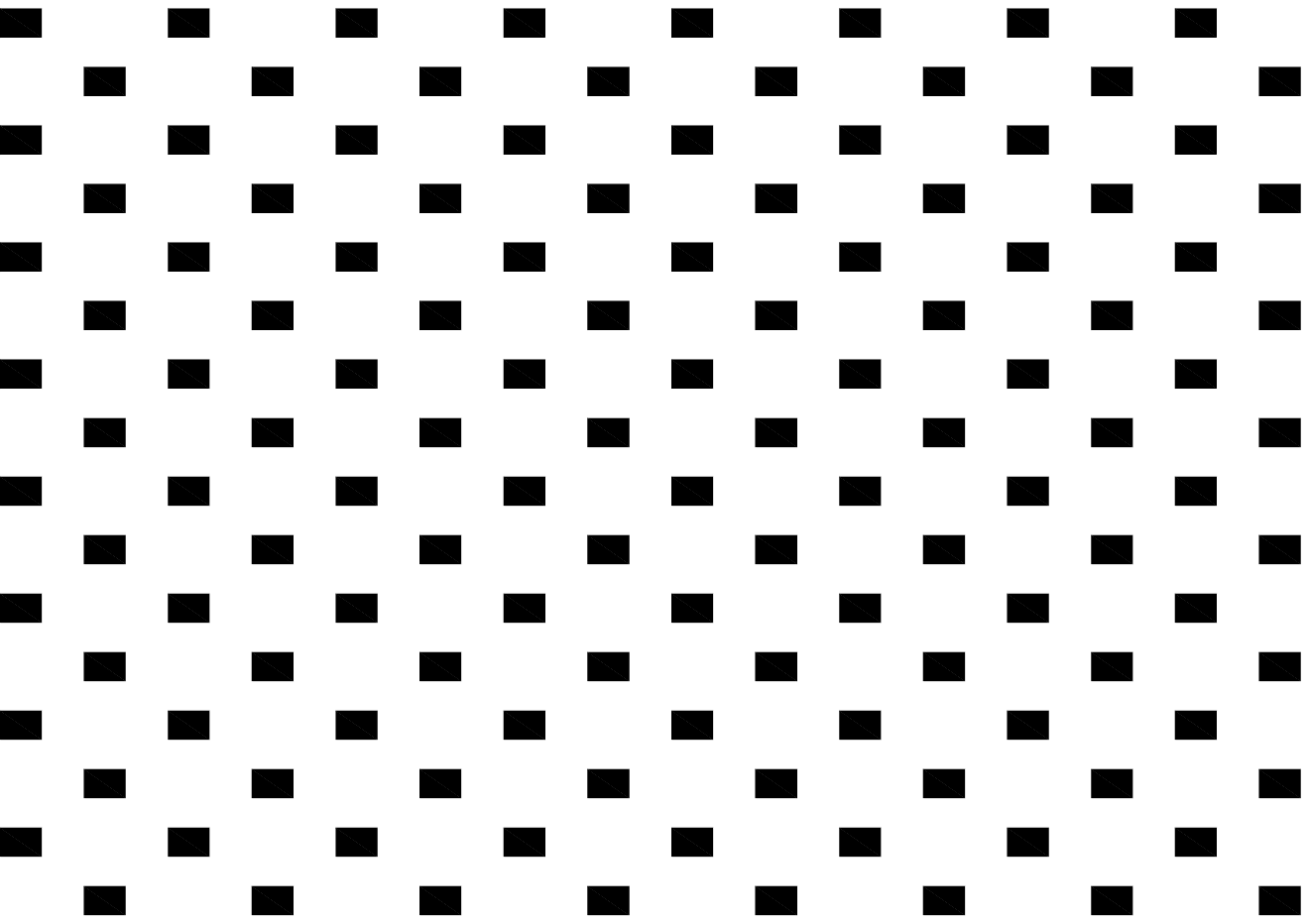}} \hfil
\frame{\includegraphics[width=0.34in]{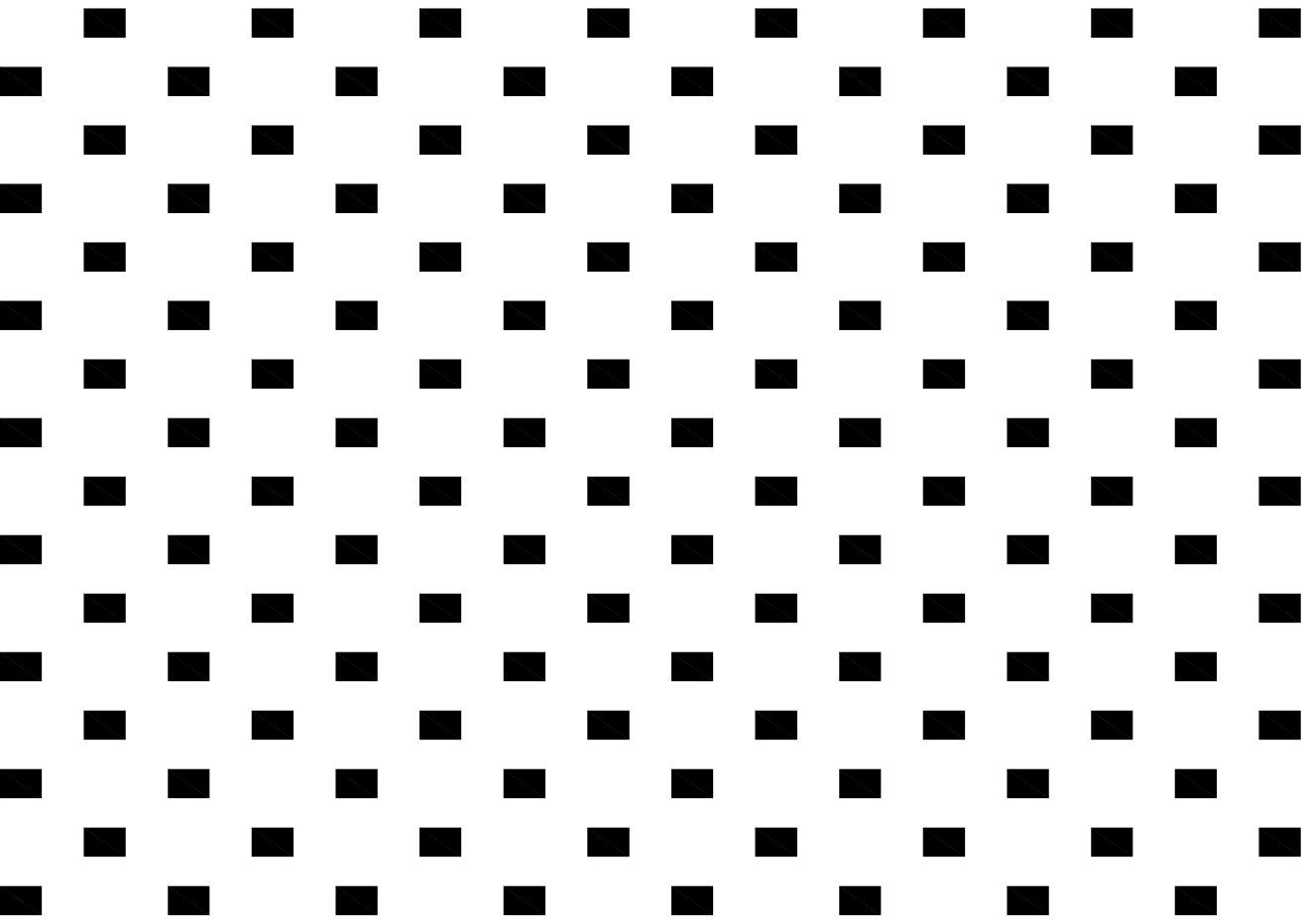}} \hfil
\frame{\includegraphics[width=0.34in]{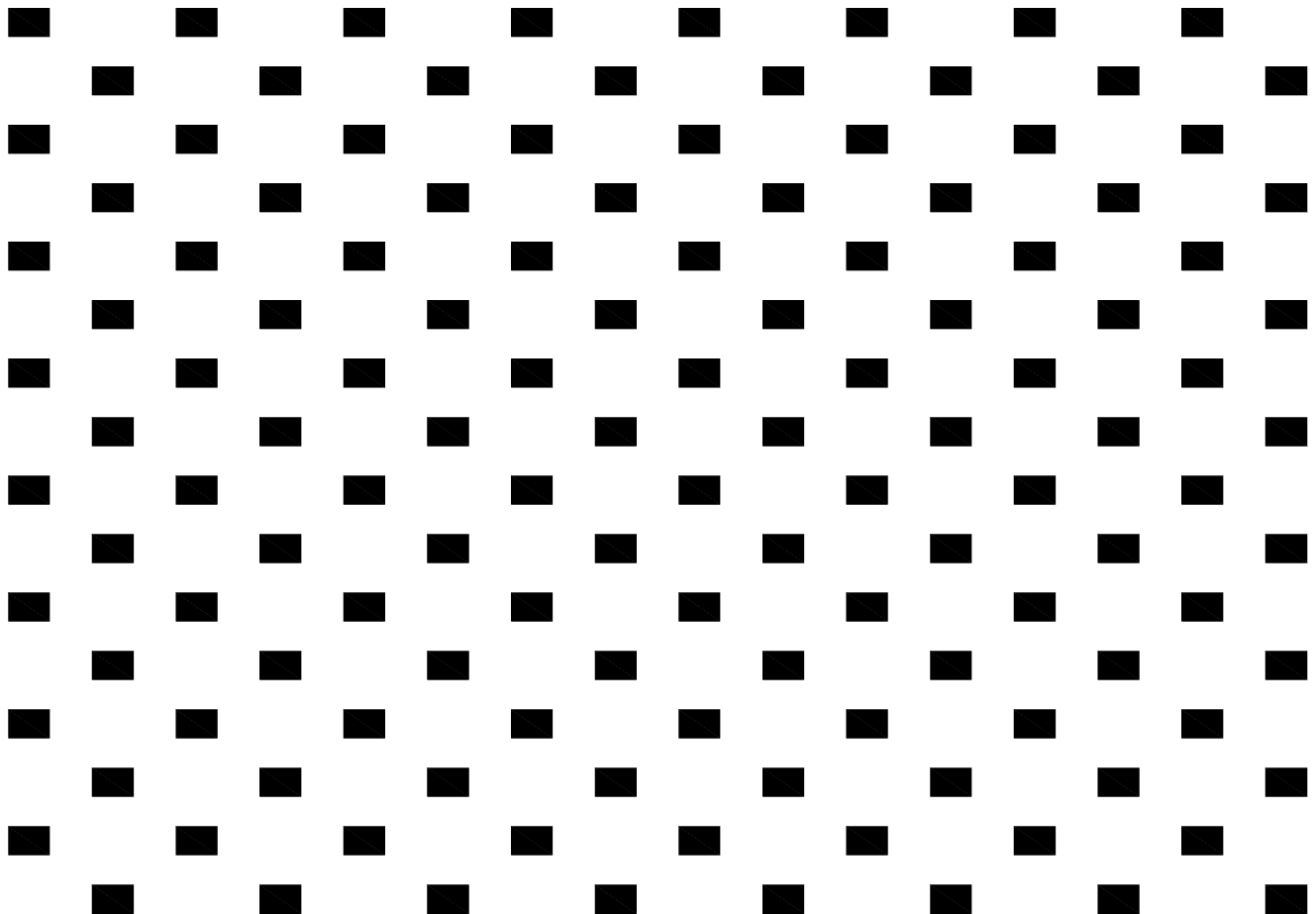}} \hfil
\frame{\includegraphics[width=0.34in]{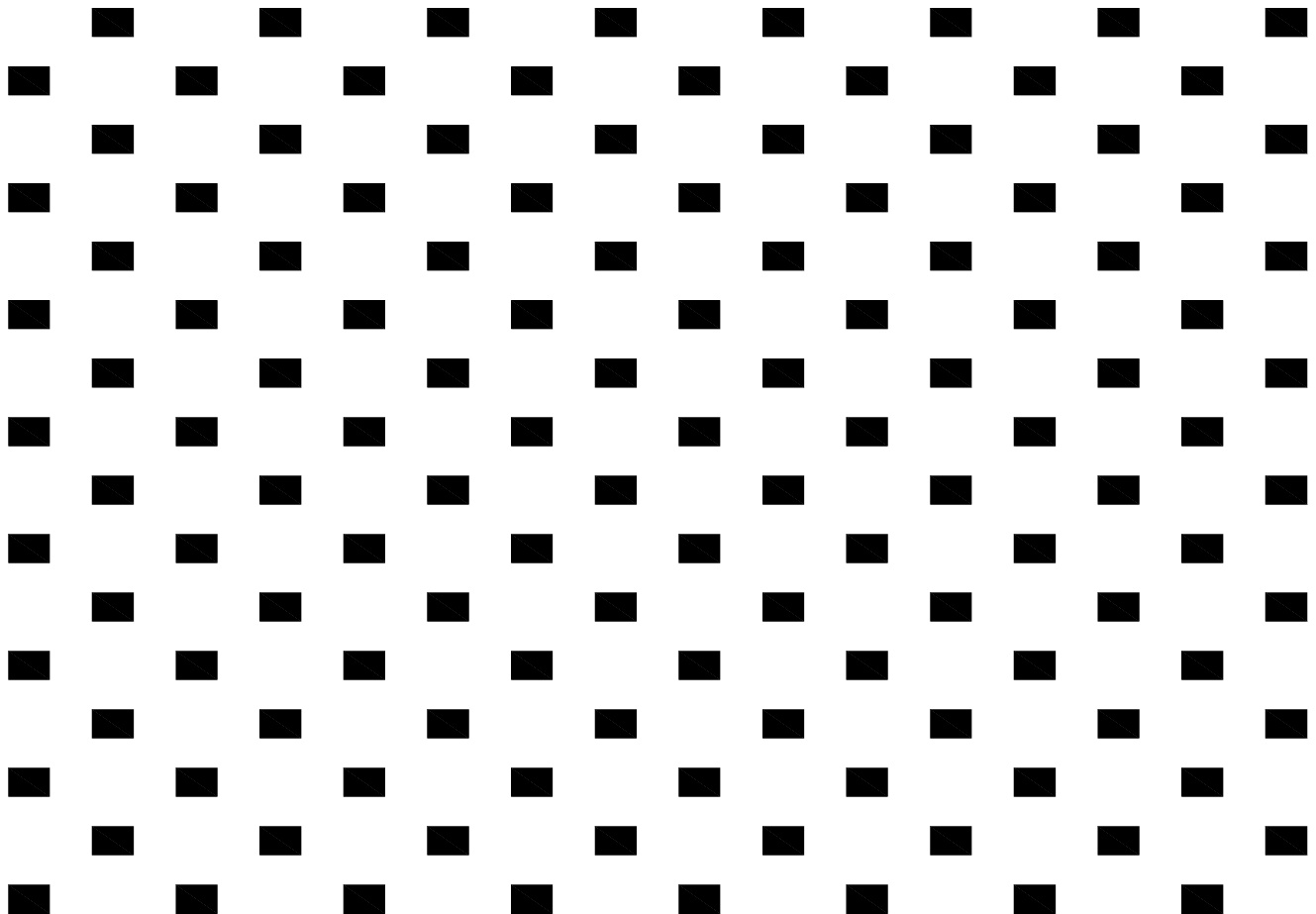}} \hfil
\frame{\includegraphics[width=0.34in]{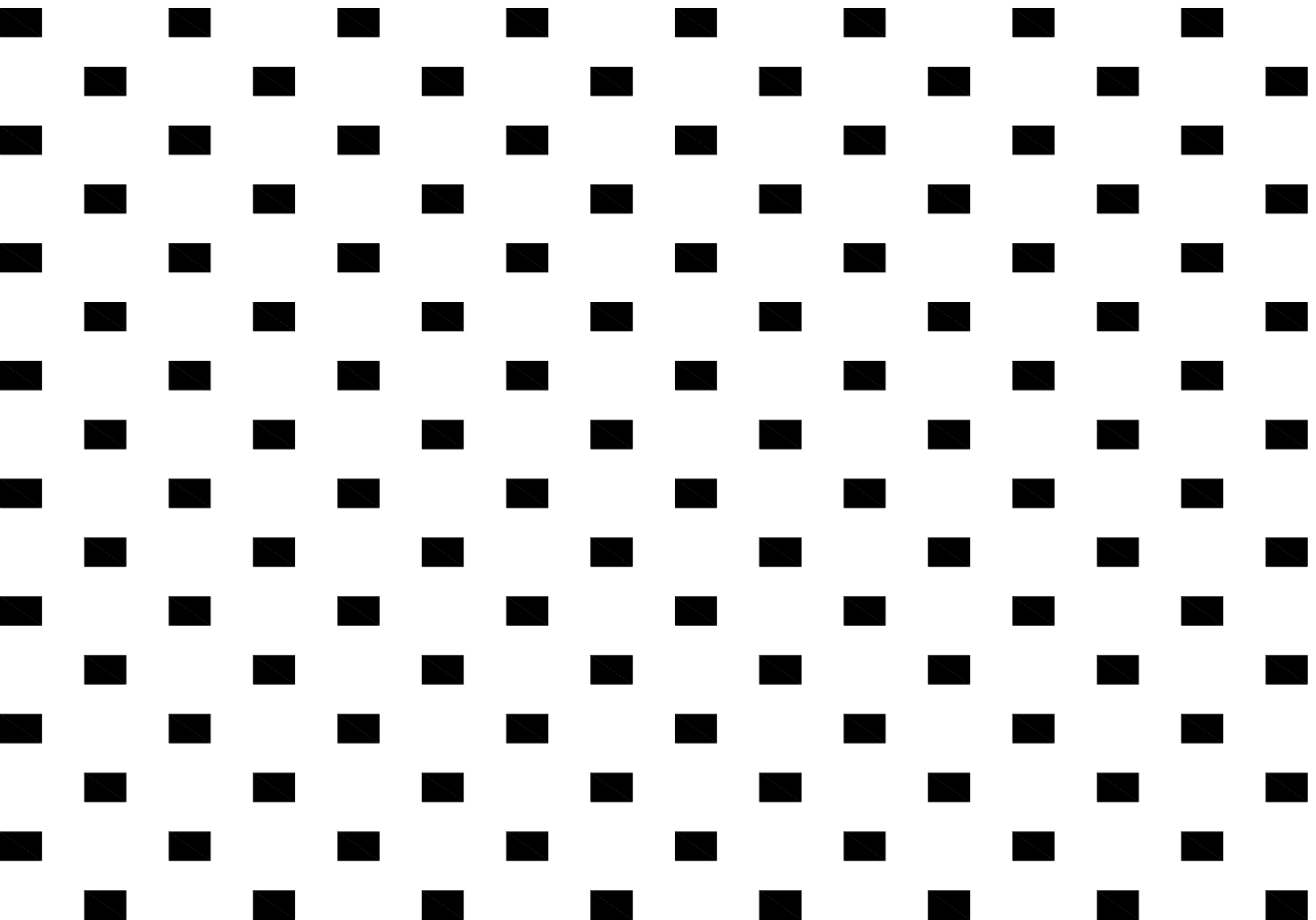}} \hfil
\frame{\includegraphics[width=0.34in]{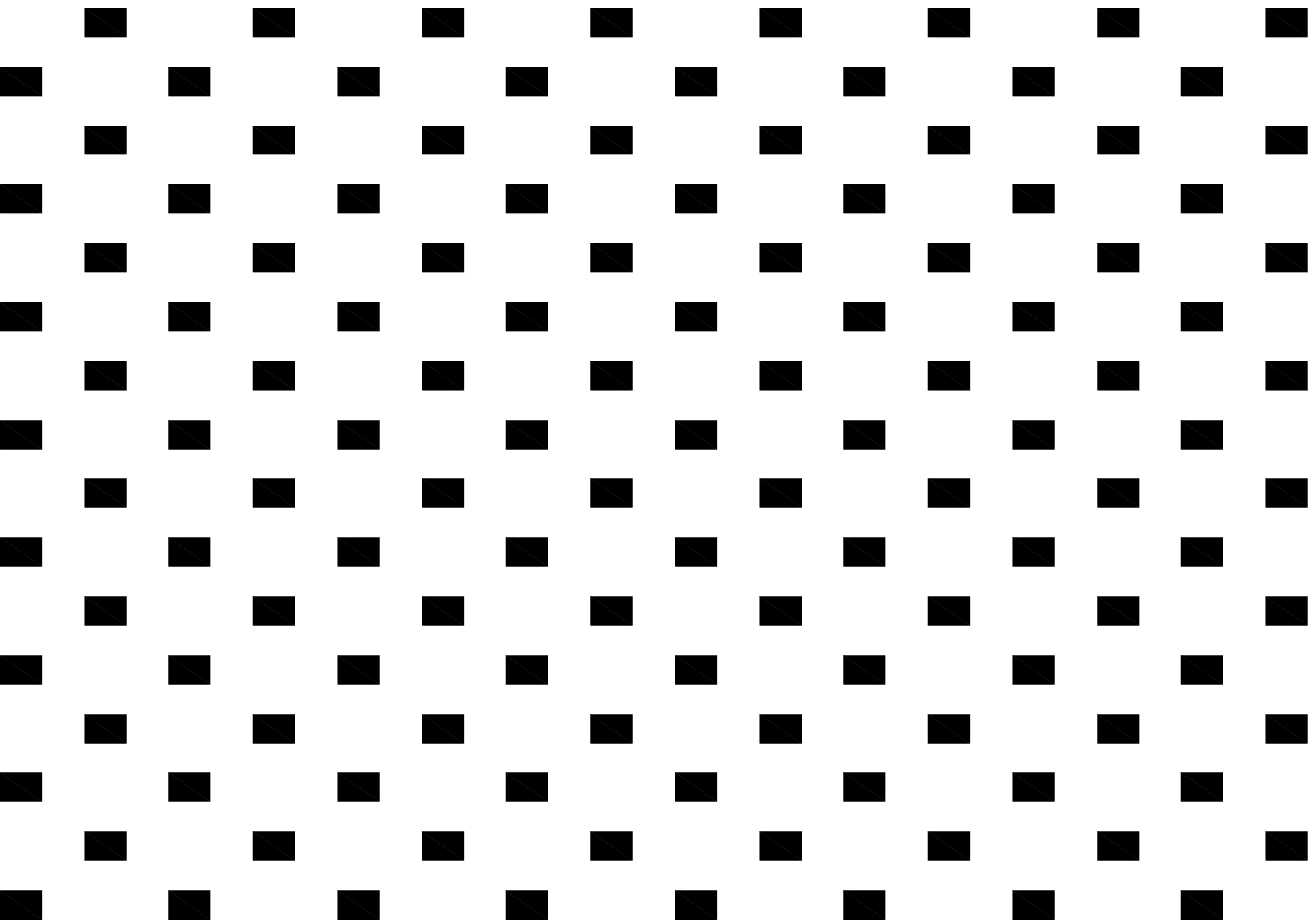}} \hfil
\frame{\includegraphics[width=0.34in]{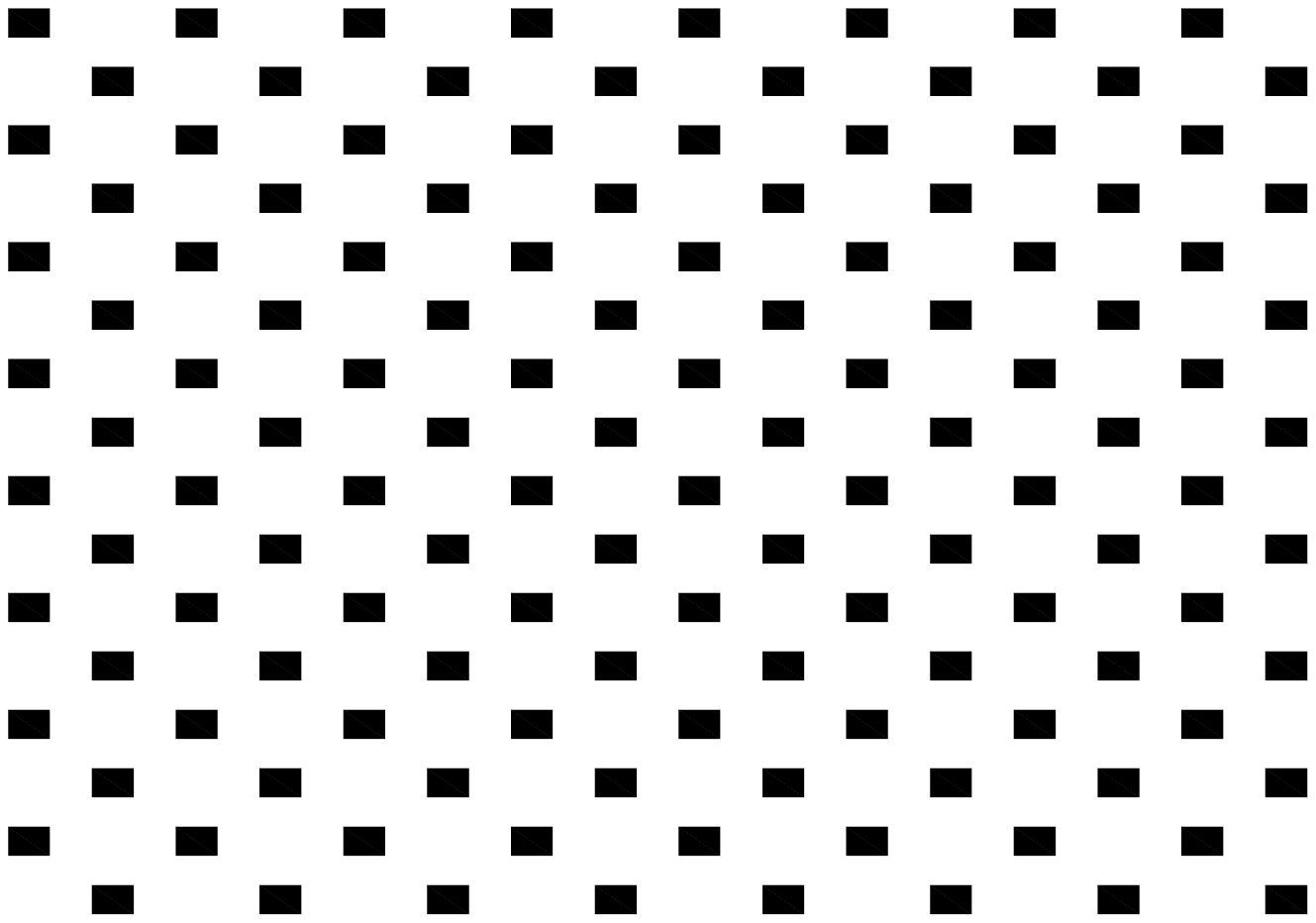}} \hfil
\frame{\includegraphics[width=0.34in]{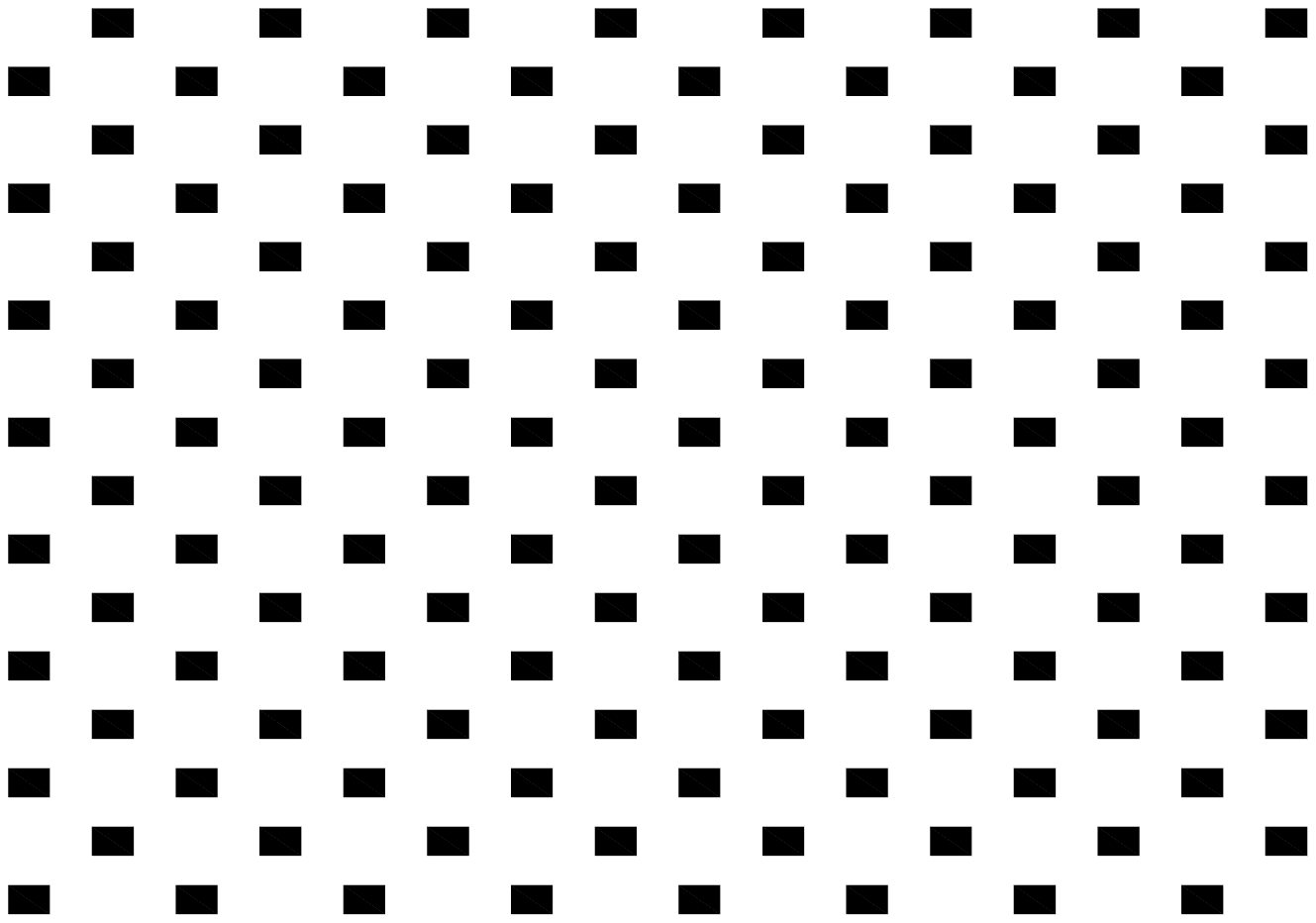}}}}
\caption{Red--black partition of nodes in a 2D square domain, $\finegrid$. Dark squares indicate the selected nodes.} \label{fig:e2_grids}
\end{figure}

The system is invariant under space shifts and therefore the eigenvectors are given by harmonic functions $(W)_{i,j}=\exp\left(\sqrt{-1}\;\tfrac{2\pi}{n}ij\right)$. After proper reordering, the basis has harmonic aliasing patterns for the down--sampling shown in Fig. \ref{fig:e2_grids}. The chessboard down--sampling shown in Fig. \ref{fig:e2_grid_1} gives a mirror matrix $\finesys^\star$ with impulse response $\left[\begin{smallmatrix} & 1 & \\1 & \underline{4-k^2} & 1\\ & 1 & \end{smallmatrix}\right]$.

\begin{figure*}[!t]
\centerline{
\subfloat[Impulse response in $\finegrid$]{\includegraphics[width=1.5in]{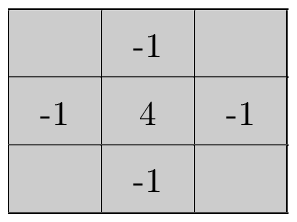}} \hfil
\subfloat[Impulse response in $\finegrid_r$]{\includegraphics[width=1.5in]{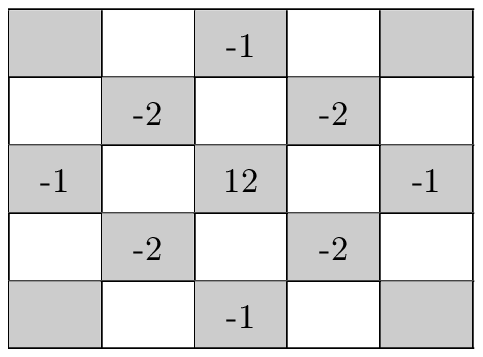}} \hfil
\subfloat[Impulse response in $\finegrid_{rr}$]{\includegraphics[width=1.5in]{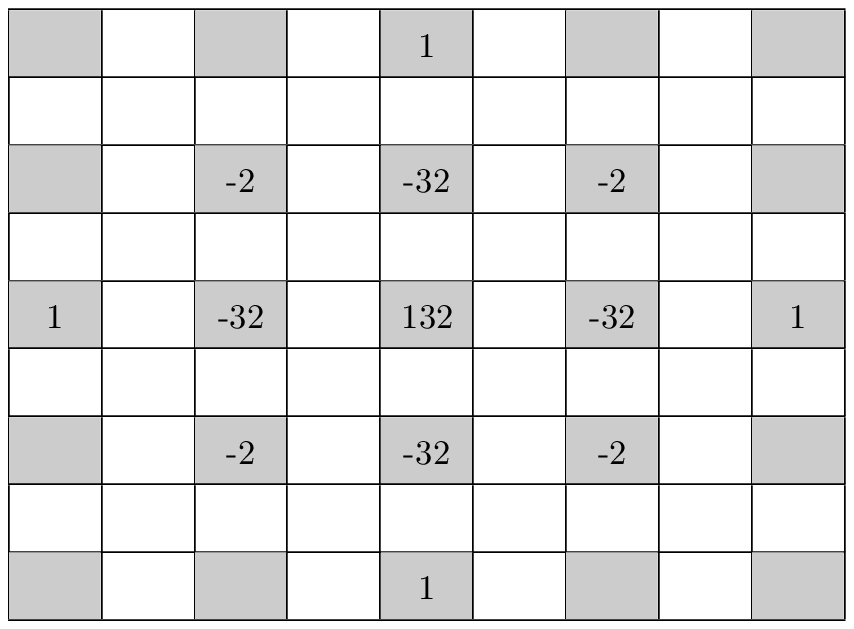}}}
\caption{Impulse response coefficients from a finite--difference discretization of the Laplacian operator: $-\nabla^2$. The coarse grid system matrices obtained with mirror filters give impulse responses that grow in coarse grids.} \label{fig:e2_stencils}
\end{figure*}

An important problem arises naturally for systems in two or more dimensions. The down/up--sampling by a factor of $2$ is not enough to reduce the impulse response of the product $\finesys\finesys^\star$. Then, the impulse response of system matrices grows larger and larger in coarse grids. An example for the impulse response of a Laplacian operator is shown in Fig. \ref{fig:e2_stencils}. The same problem appears in other direct multi--grid solvers like total and partial (or cyclic) reduction methods \cite{JSchroder_UTrottenberg_1973a,JSchroder_UTrottenberg_HReutersberg_1976a,BLBuzbee_GHGolub_CWNielson_1970a}.

\begin{figure*}[!t]
\centerline{
\subfloat[Source vector: $f$.]{\frame{\includegraphics[width=1.5in]{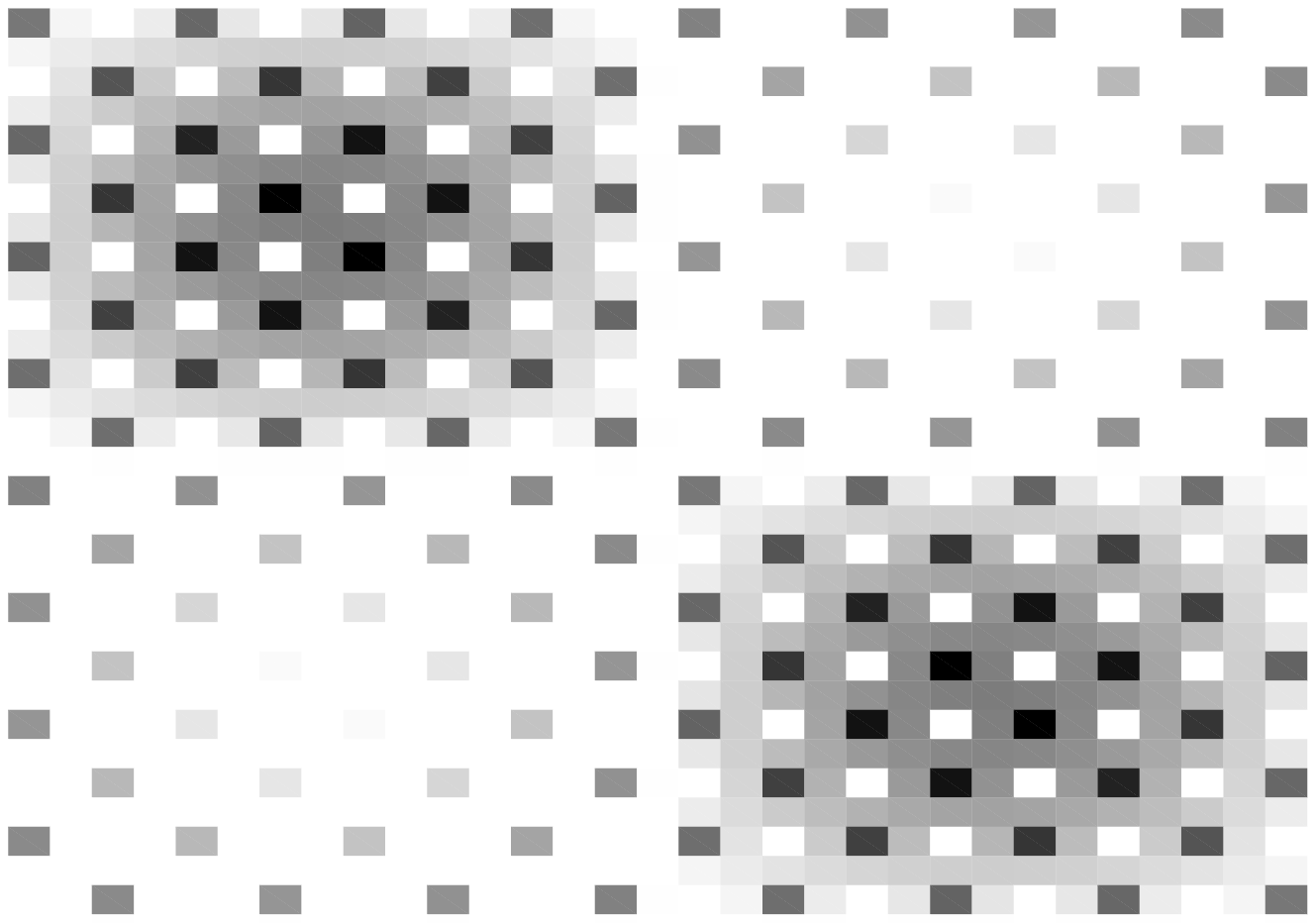}}\label{fig:e2_source1}} \hfil
\subfloat[Exact solution: $u$.]{\frame{\includegraphics[width=1.5in]{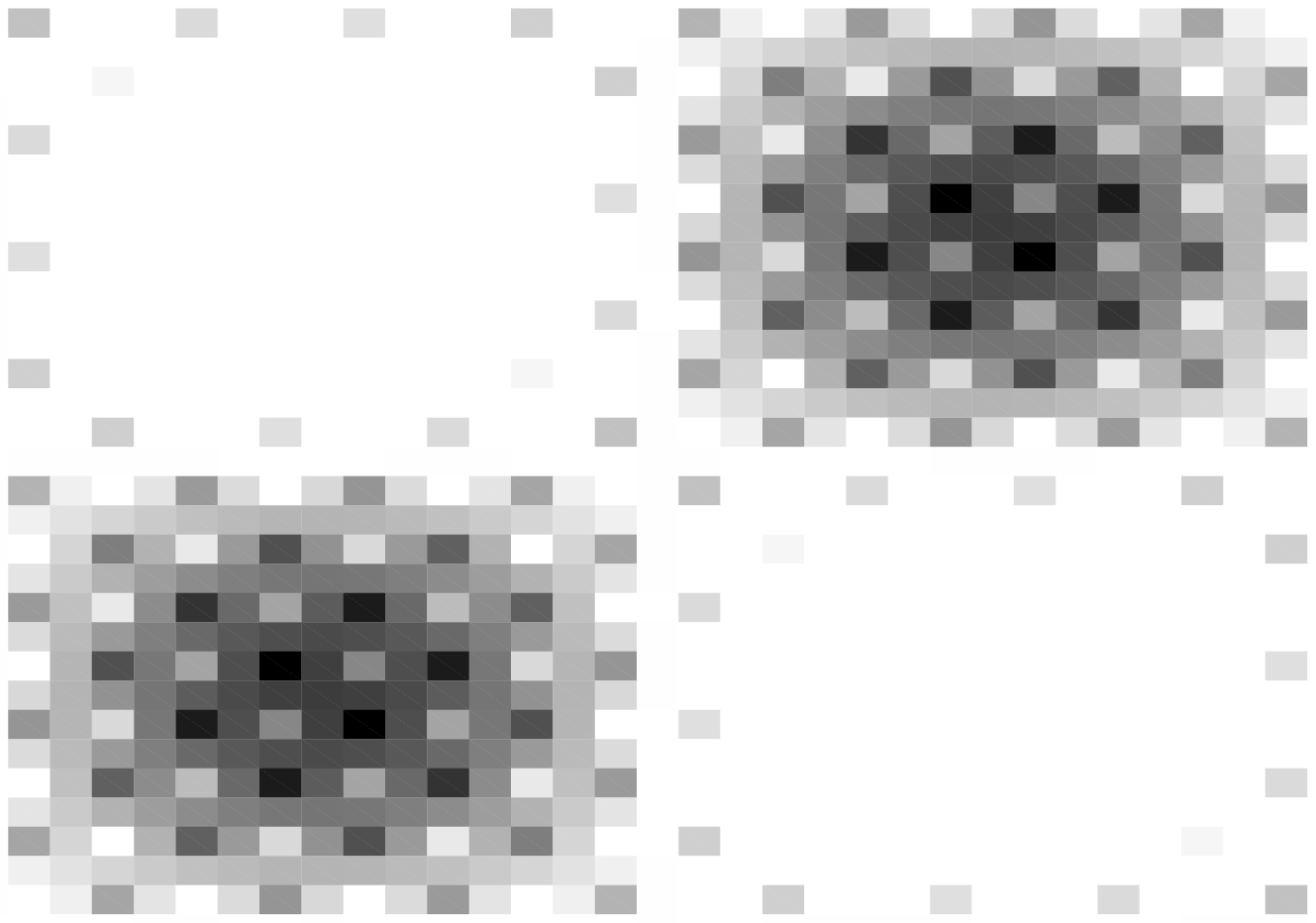}}\label{fig:e2_solution1}} \hfil
\subfloat[Source vector: $f$.]{\frame{\includegraphics[width=1.5in]{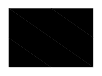}}\label{fig:e2_source2}} \hfil
\subfloat[Exact solution: $u$.]{\frame{\includegraphics[width=1.5in]{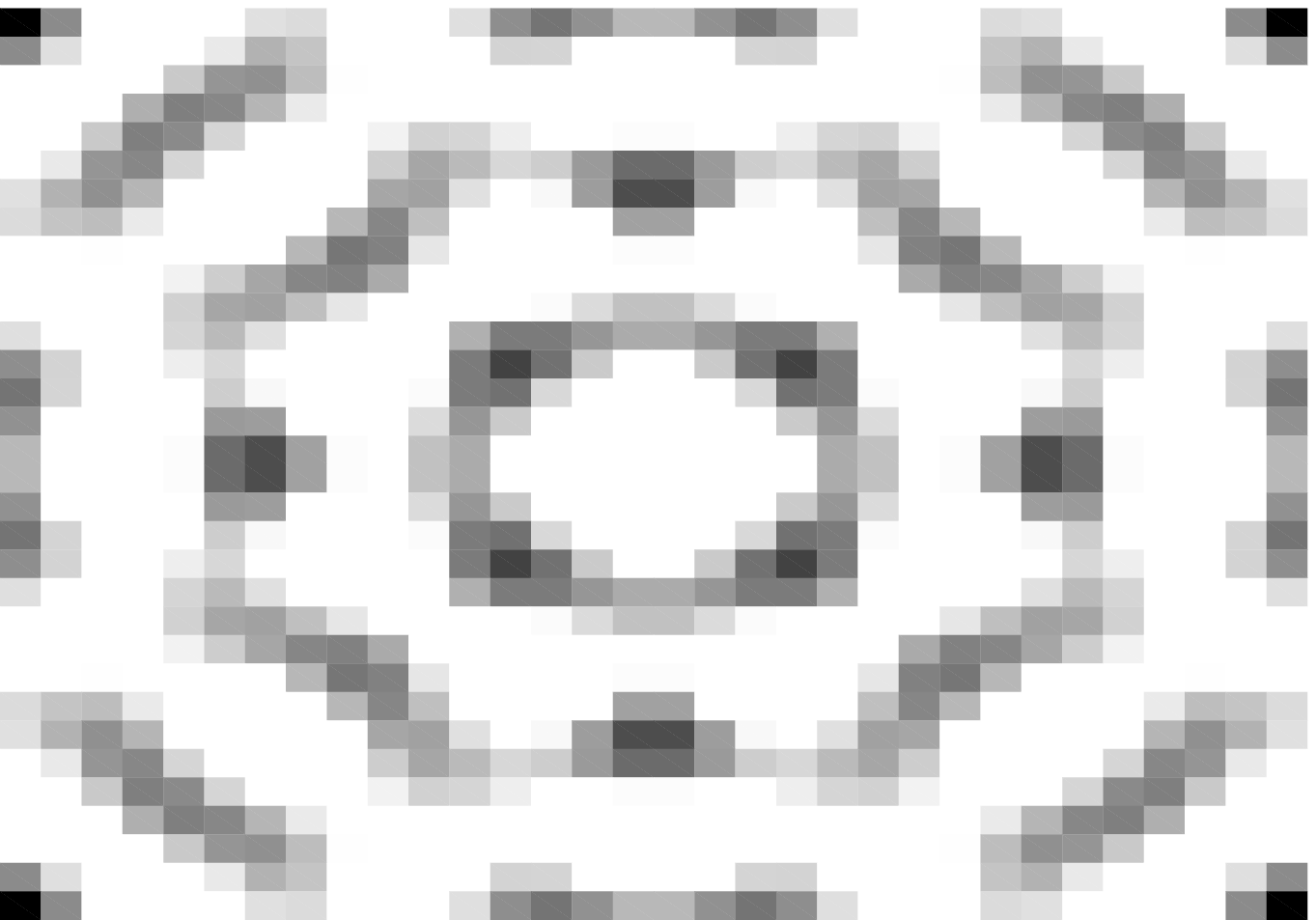}}\label{fig:e2_solution2}}}
\caption{Helmholtz's equation, $-\nabla^2 u-k^2 u=0$, with $k=\tfrac{\pi}{3}$, is solved in a square domain with periodic boundary conditions for two source vectors. In Fig. \ref{fig:e2_source1} the source vector is set to $f(i,j)=\sin(i\pi/16)\sin(j\pi/16)+\sin(i\pi/2)\sin(j\pi/2)$, mixing low-- and high--frequency sources. The exact solution is shown in Fig. \ref{fig:e2_solution1}. In Fig. \ref{fig:e2_source2} the source vector is equal to $1$ in four neighboring nodes at the center of the figure and zero elsewhere. The exact solution is shown in \ref{fig:e2_solution2}.}
\end{figure*}

\begin{figure*}[!t]
\centerline{
\subfloat[$v_0$ for $f$ in Fig. \ref{fig:e2_source1}]{\frame{\includegraphics[width=1.5in]{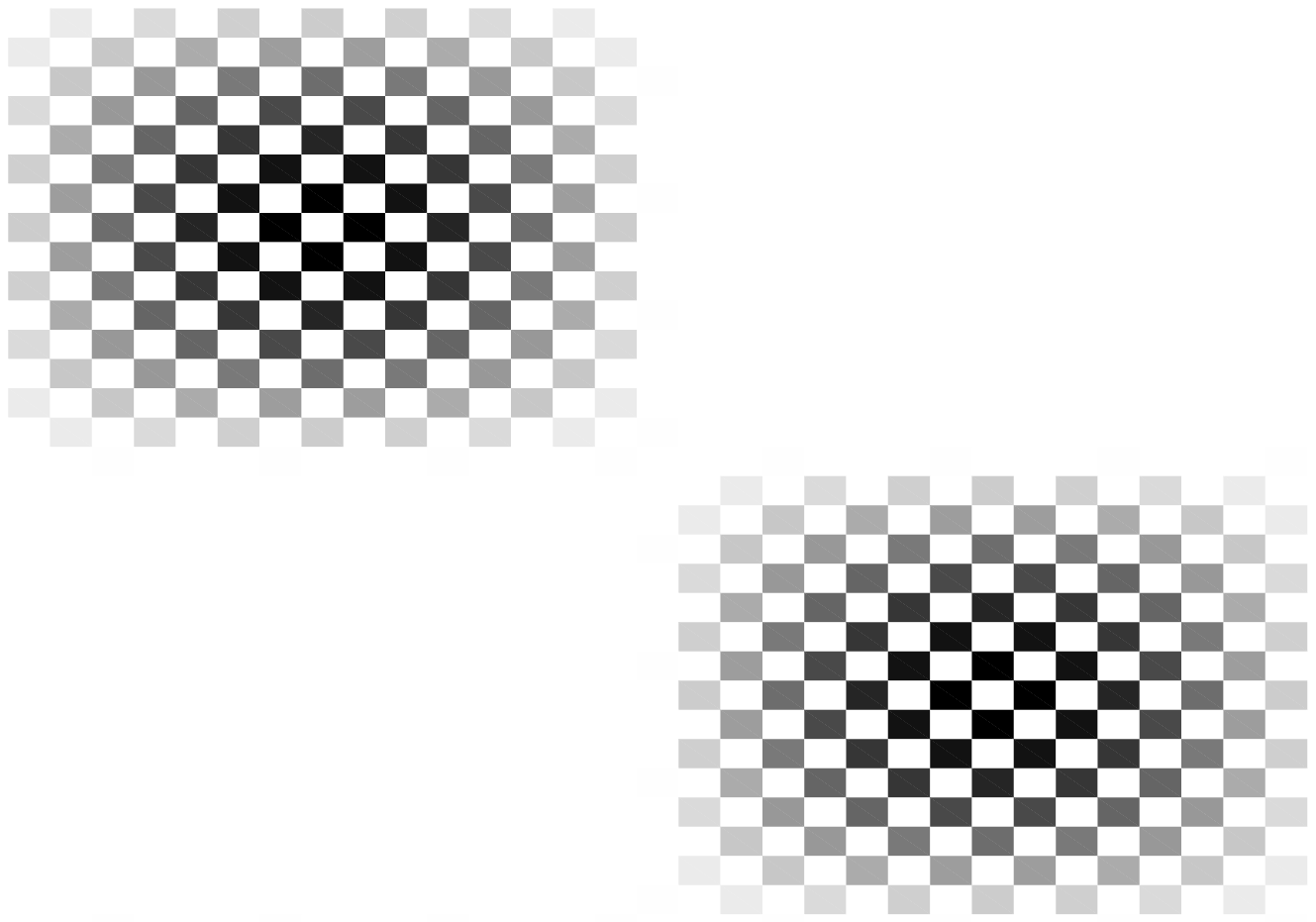}}\label{fig:e2_multiplicative_v0_1}} \hfil
\subfloat[$e_0$ for $f$ in Fig. \ref{fig:e2_source1}]{\frame{\includegraphics[width=1.5in]{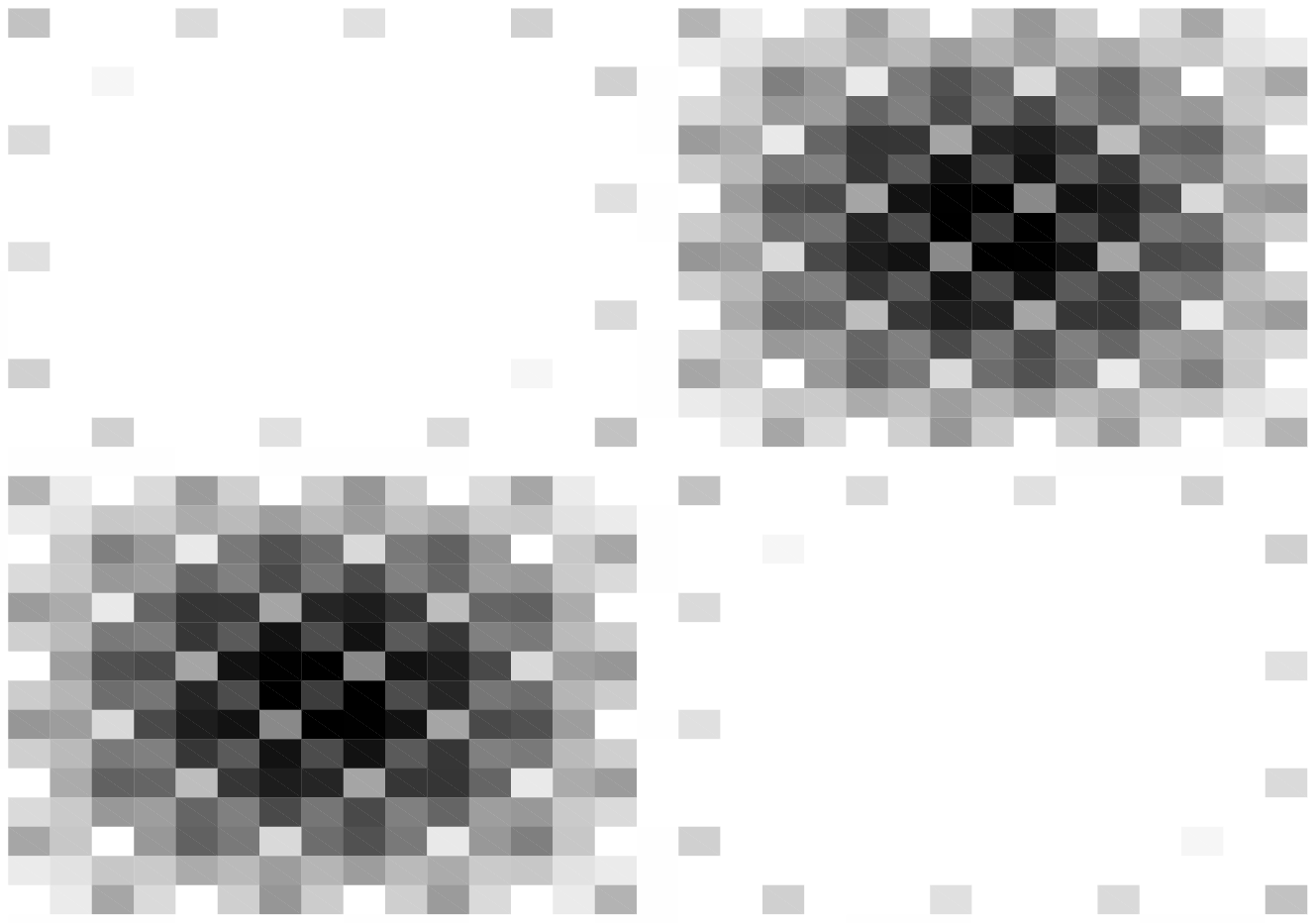}}\label{fig:e2_multiplicative_e0_1}} \hfil
\subfloat[$v_0$ for $f$ in Fig. \ref{fig:e2_source2}]{\frame{\includegraphics[width=1.5in]{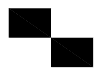}}\label{fig:e2_multiplicative_v0_2}} \hfil
\subfloat[$e_0$ for $f$ in Fig. \ref{fig:e2_source2}]{\frame{\includegraphics[width=1.5in]{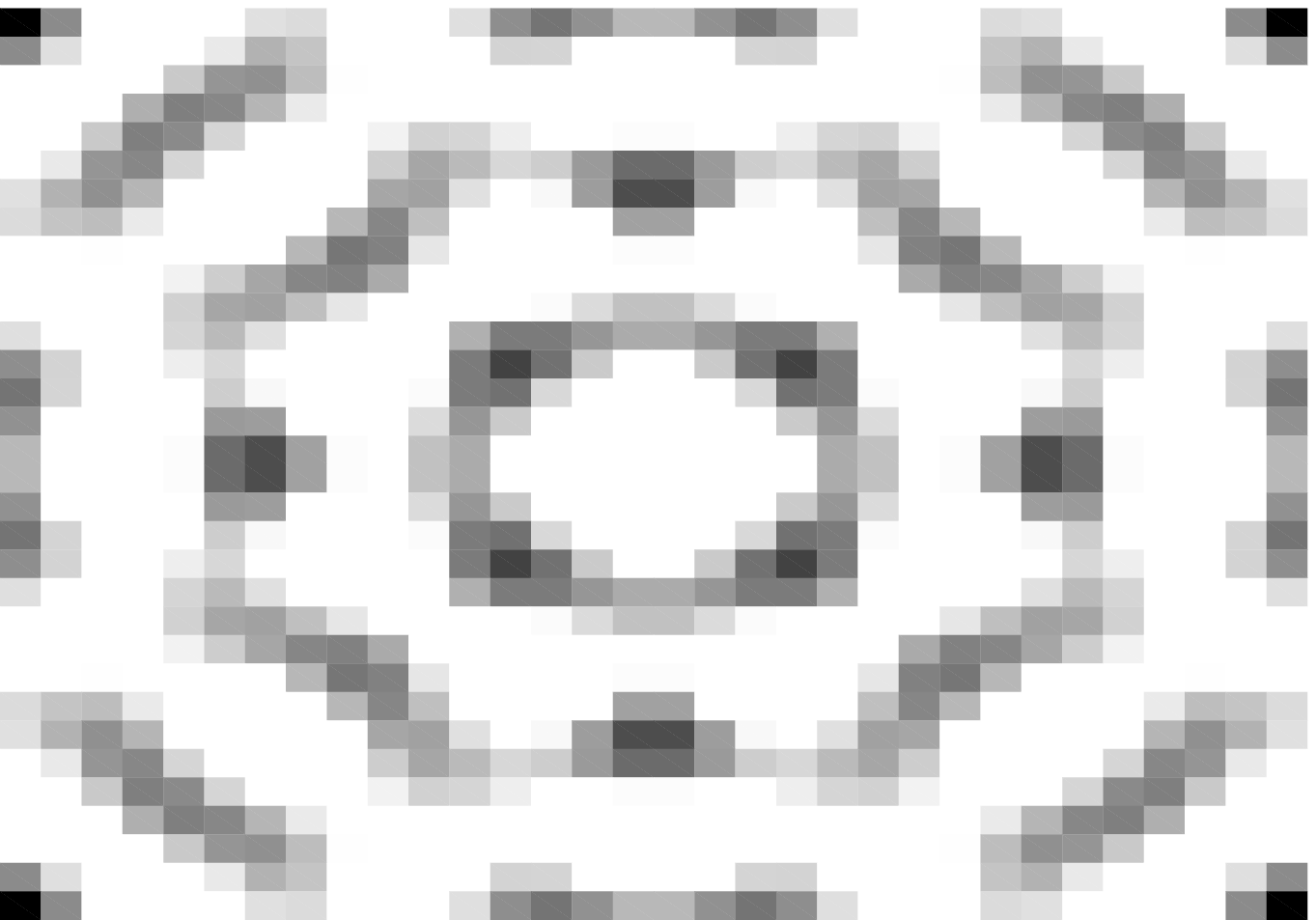}}\label{fig:e2_multiplicative_e0_2}}}
\caption{Intermediate solutions of Helmholtz's equation using the direct two--grid algorithm with a multiplicative approach.} \label{fig:e2_multiplicative_solutions}
\end{figure*}

\begin{figure*}[!t]
\centerline{
\subfloat[$u_{r}$ and $u_{b}$]{ \label{fig:e2_additive_l1_1}
\frame{\includegraphics[width=1.5in]{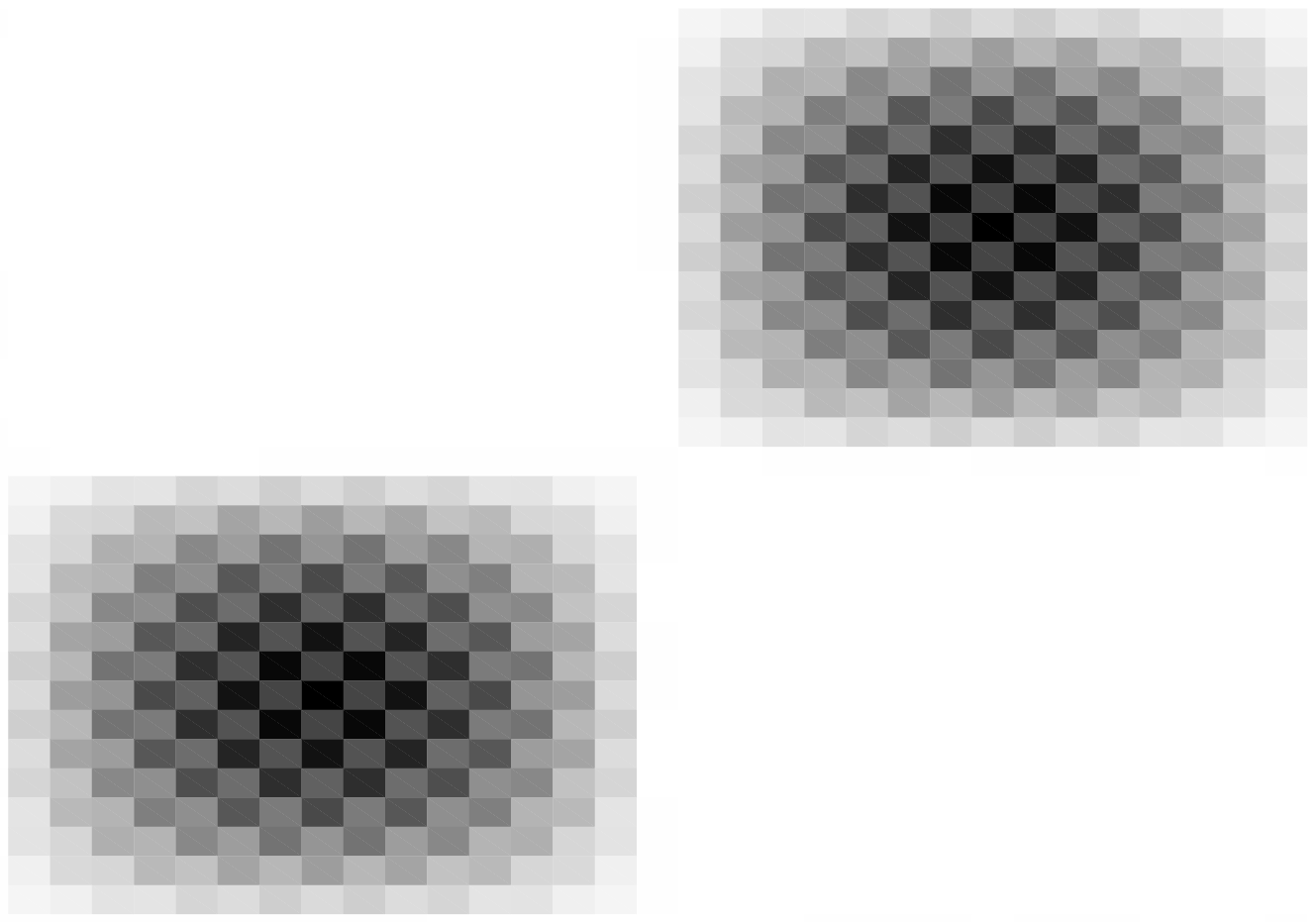}} \hfil
\frame{\includegraphics[width=1.5in]{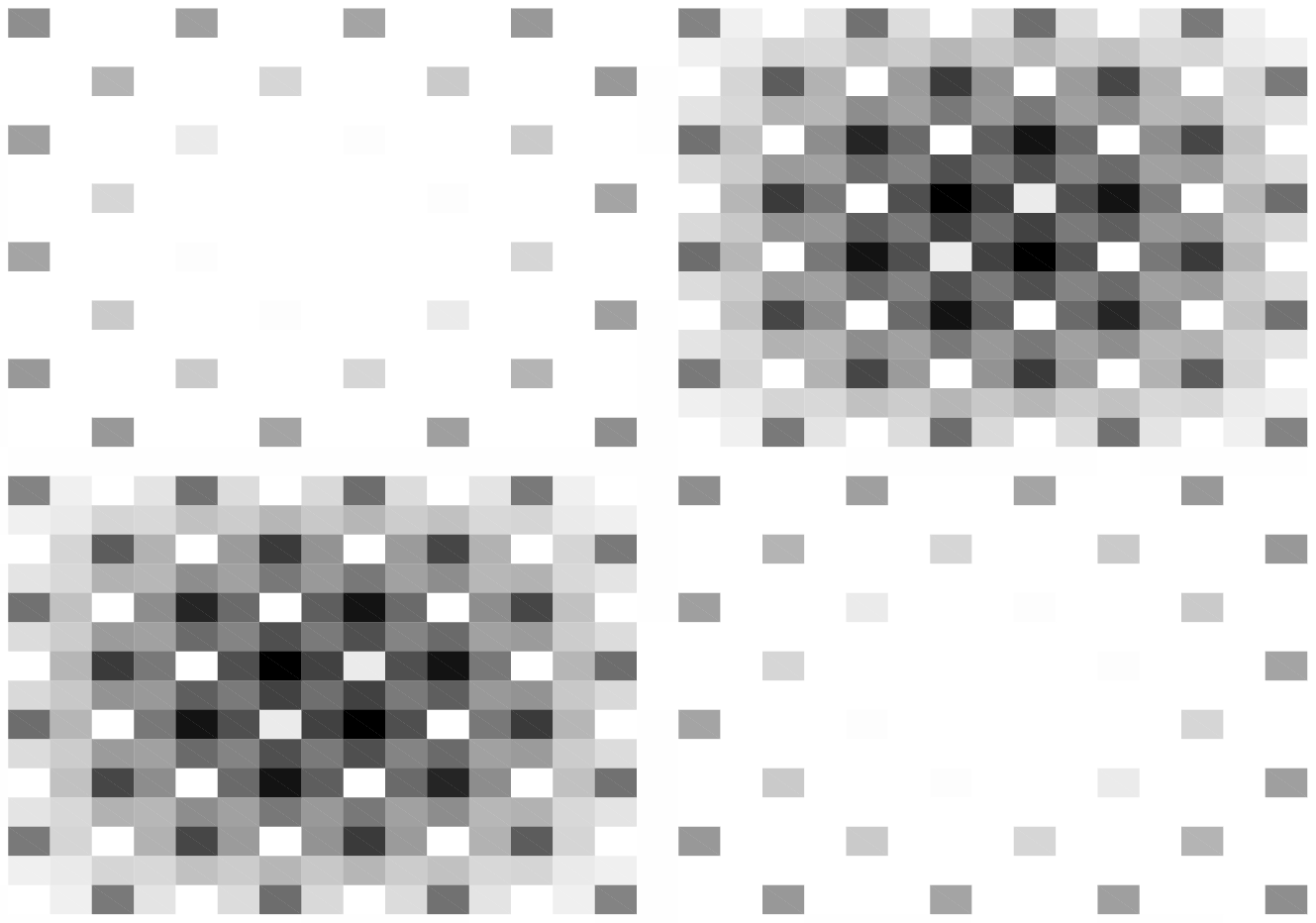}}} \hfil
\subfloat[$u_{r}$ and $u_{b}$]{ \label{fig:e2_additive_l1_2}
\frame{\includegraphics[width=1.5in]{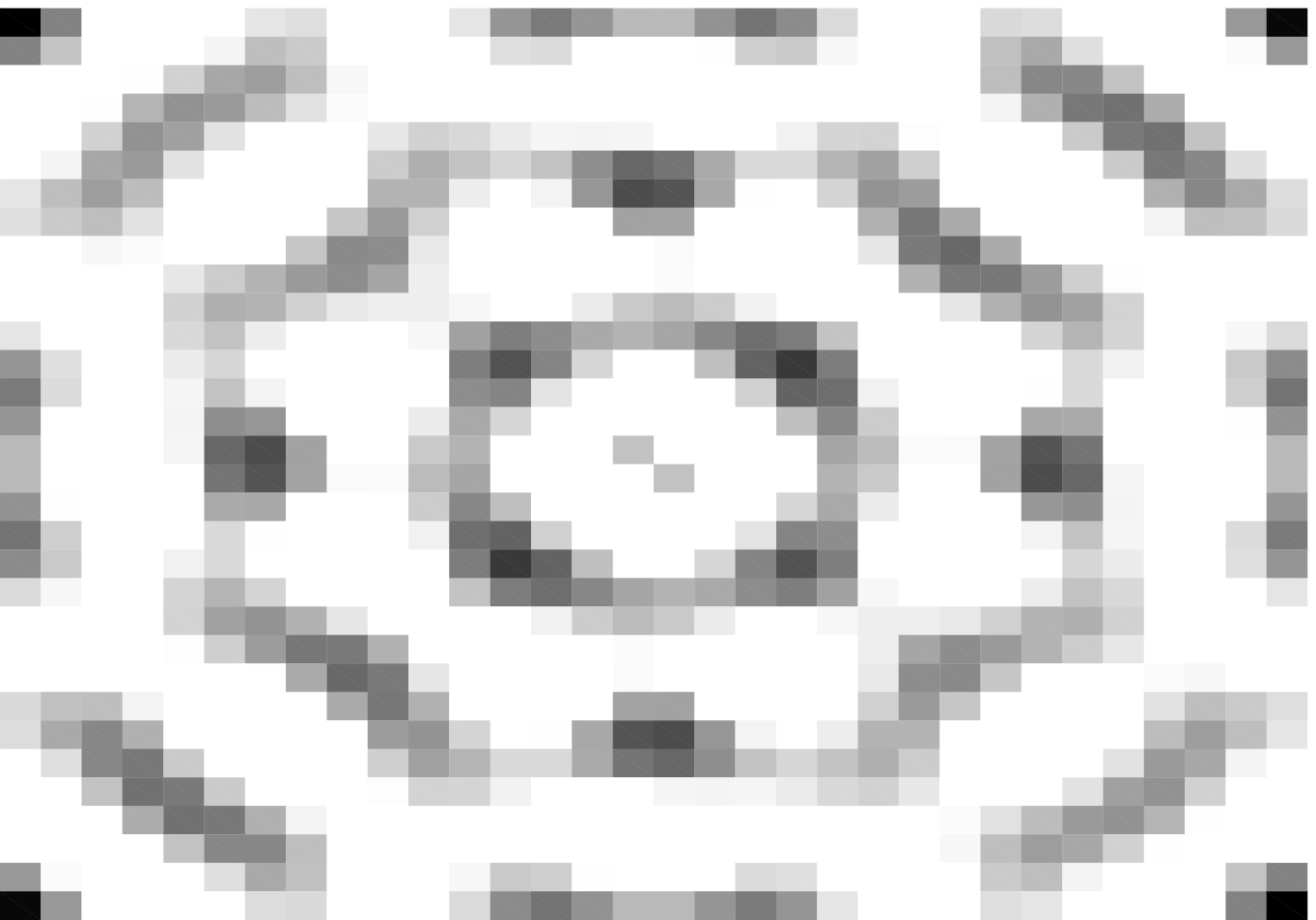}} \hfil
\frame{\includegraphics[width=1.5in]{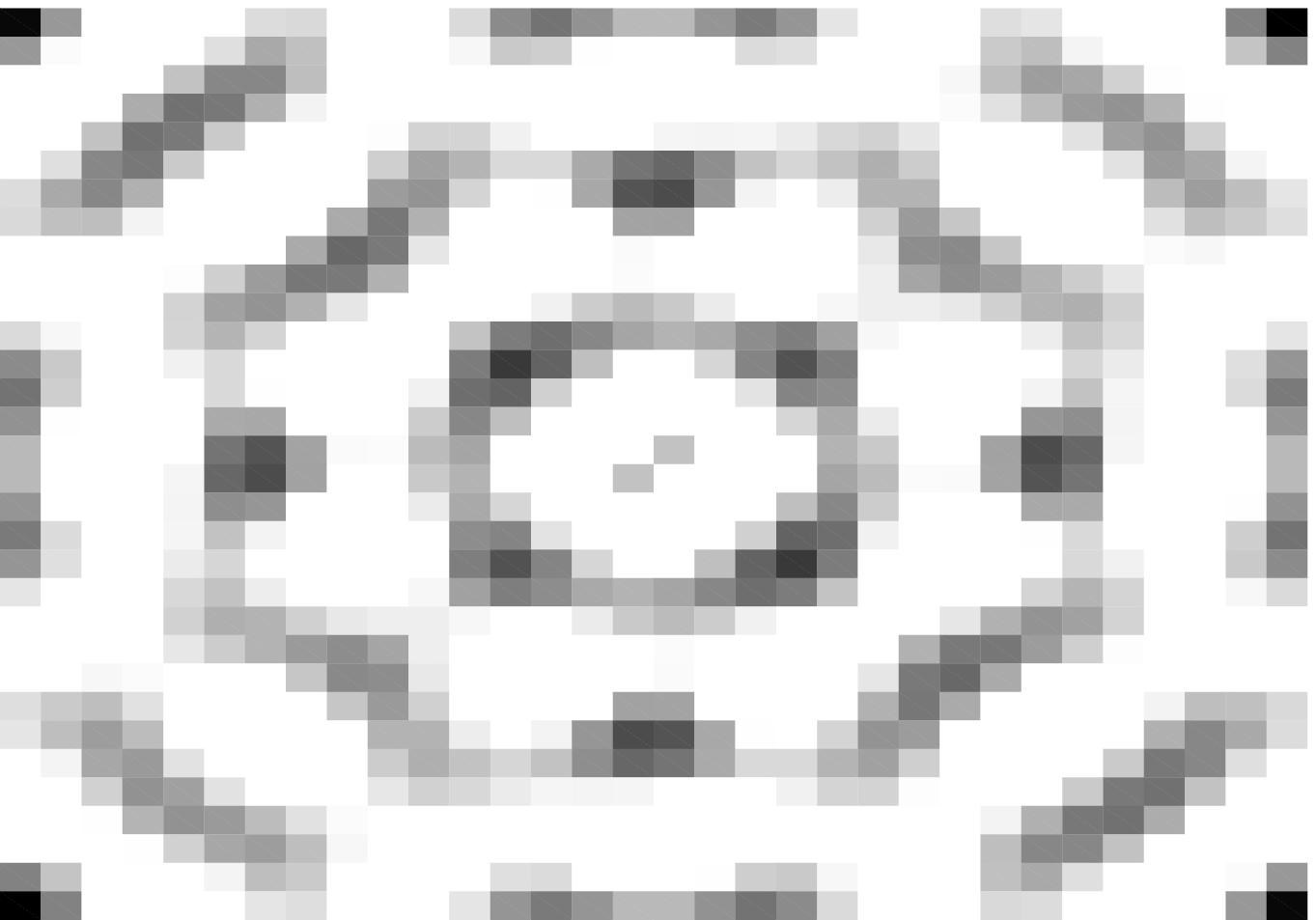}}}}
\centerline{
\subfloat[$u_{rr}$, $u_{rb}$, $u_{br}$ and $u_{bb}$]{ \label{fig:e2_additive_l2_1}
\frame{\includegraphics[width=0.74in]{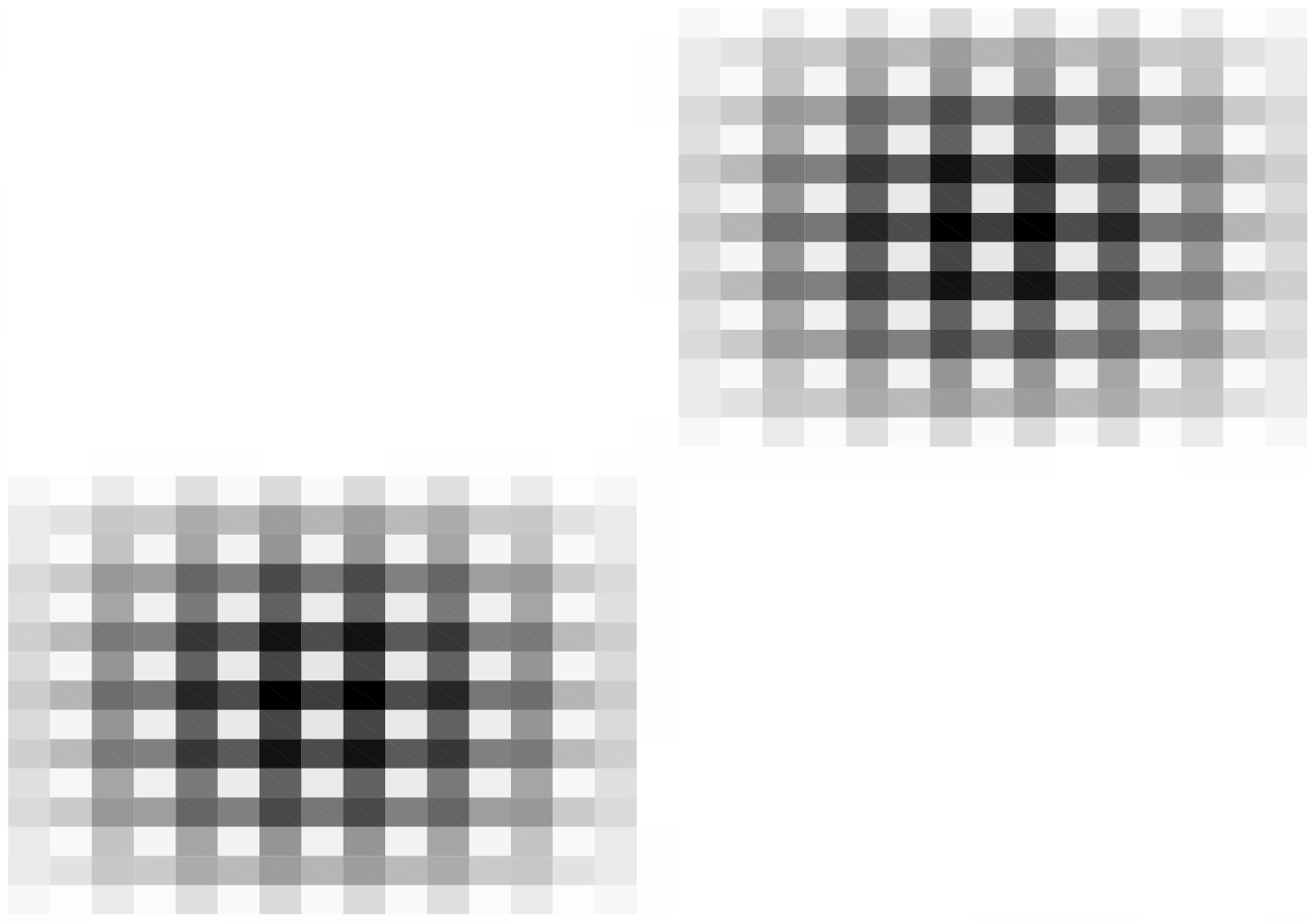}} \hfil
\frame{\includegraphics[width=0.74in]{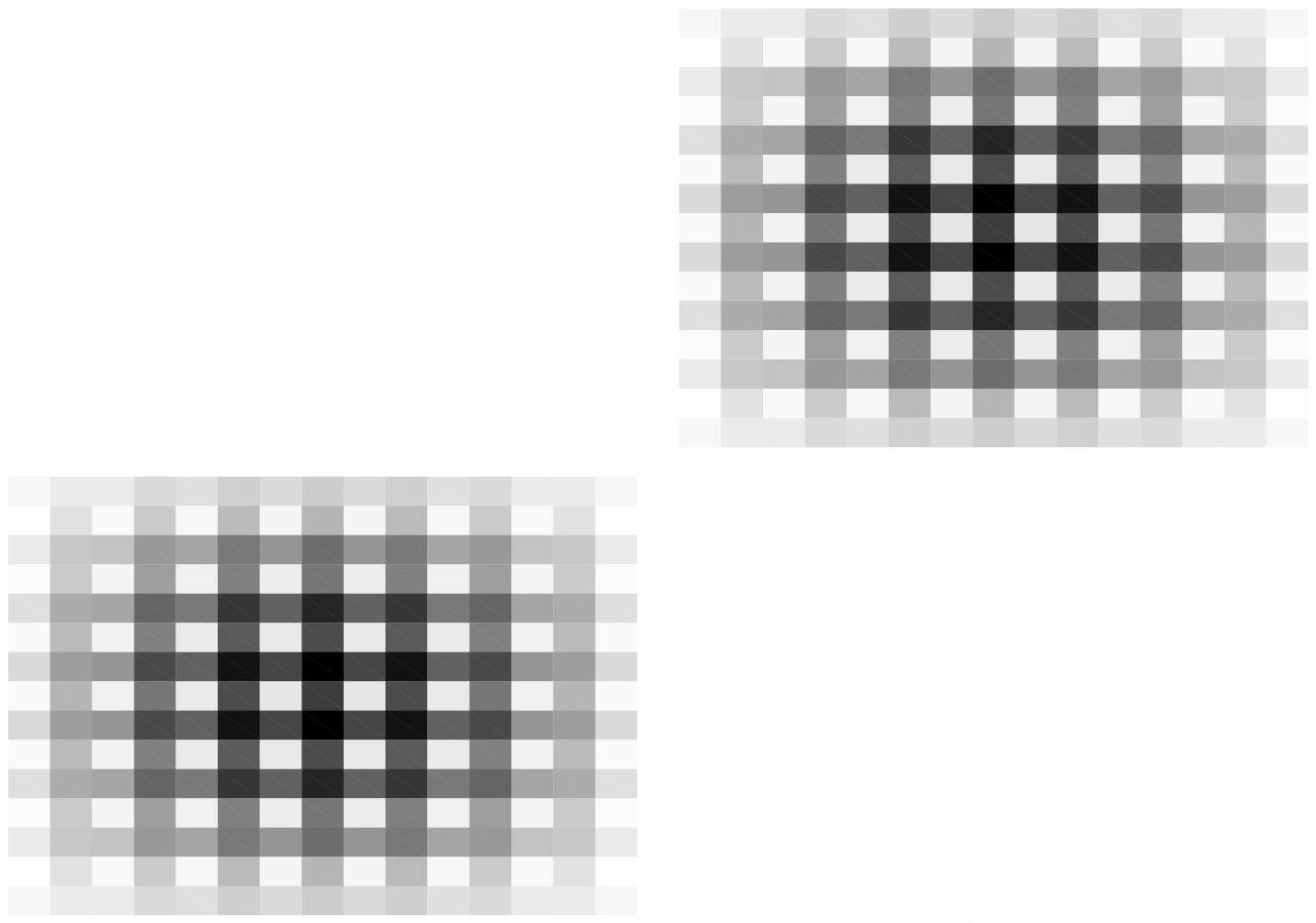}} \hfil
\frame{\includegraphics[width=0.74in]{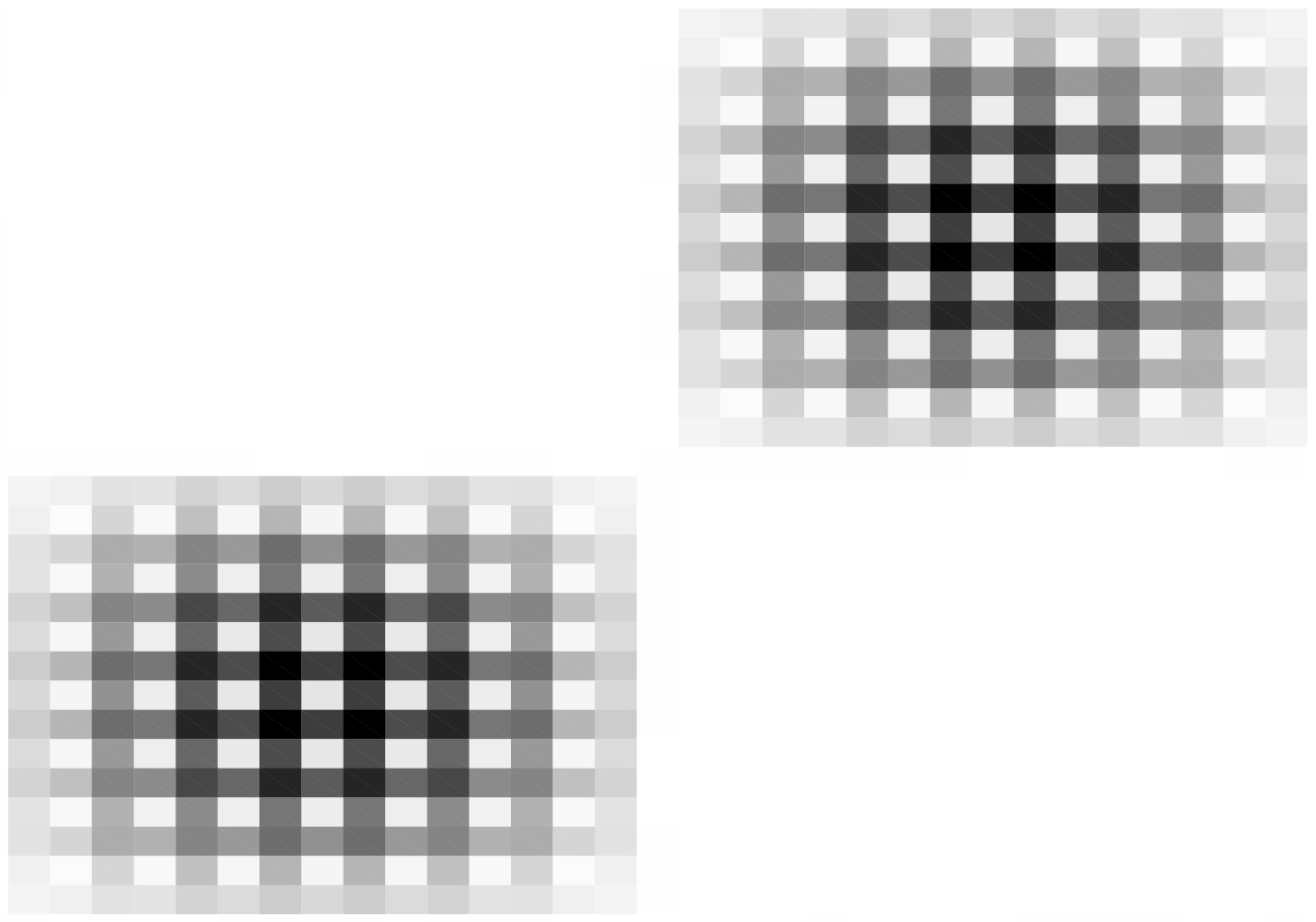}} \hfil
\frame{\includegraphics[width=0.74in]{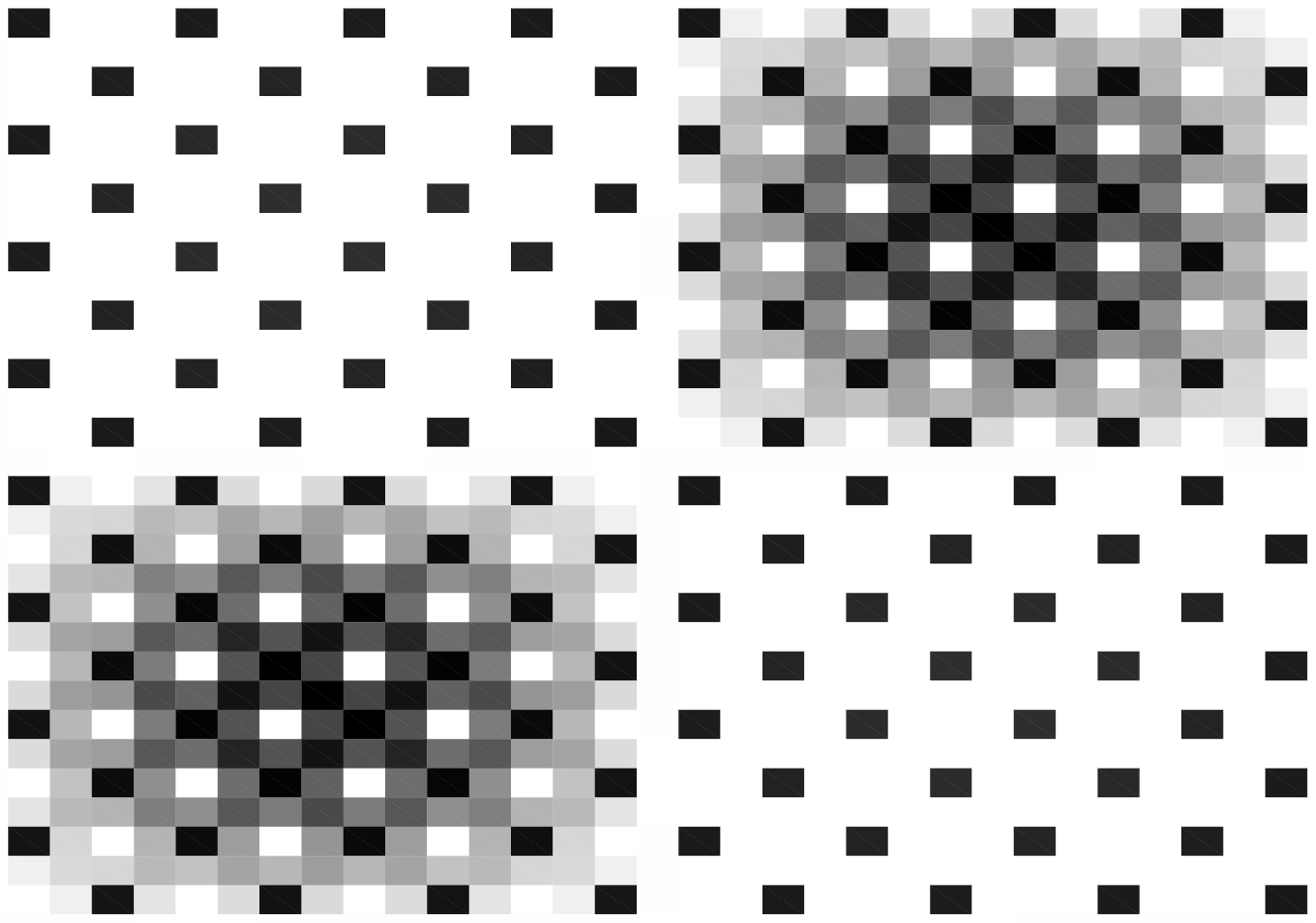}}} \hfil
\subfloat[$u_{rr}$, $u_{rb}$, $u_{br}$ and $u_{bb}$]{ \label{fig:e2_additive_l2_2}
\frame{\includegraphics[width=0.74in]{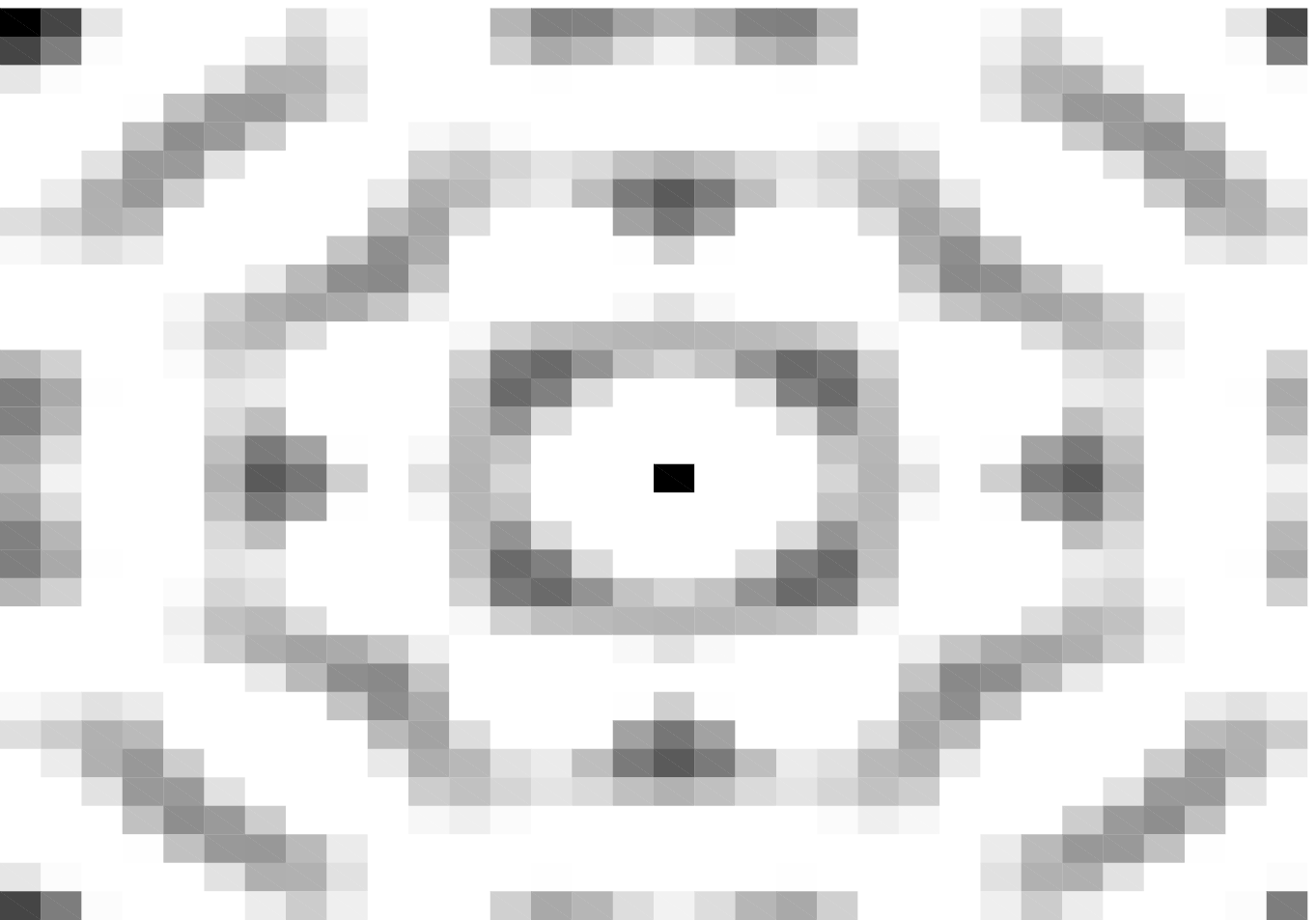}} \hfil
\frame{\includegraphics[width=0.74in]{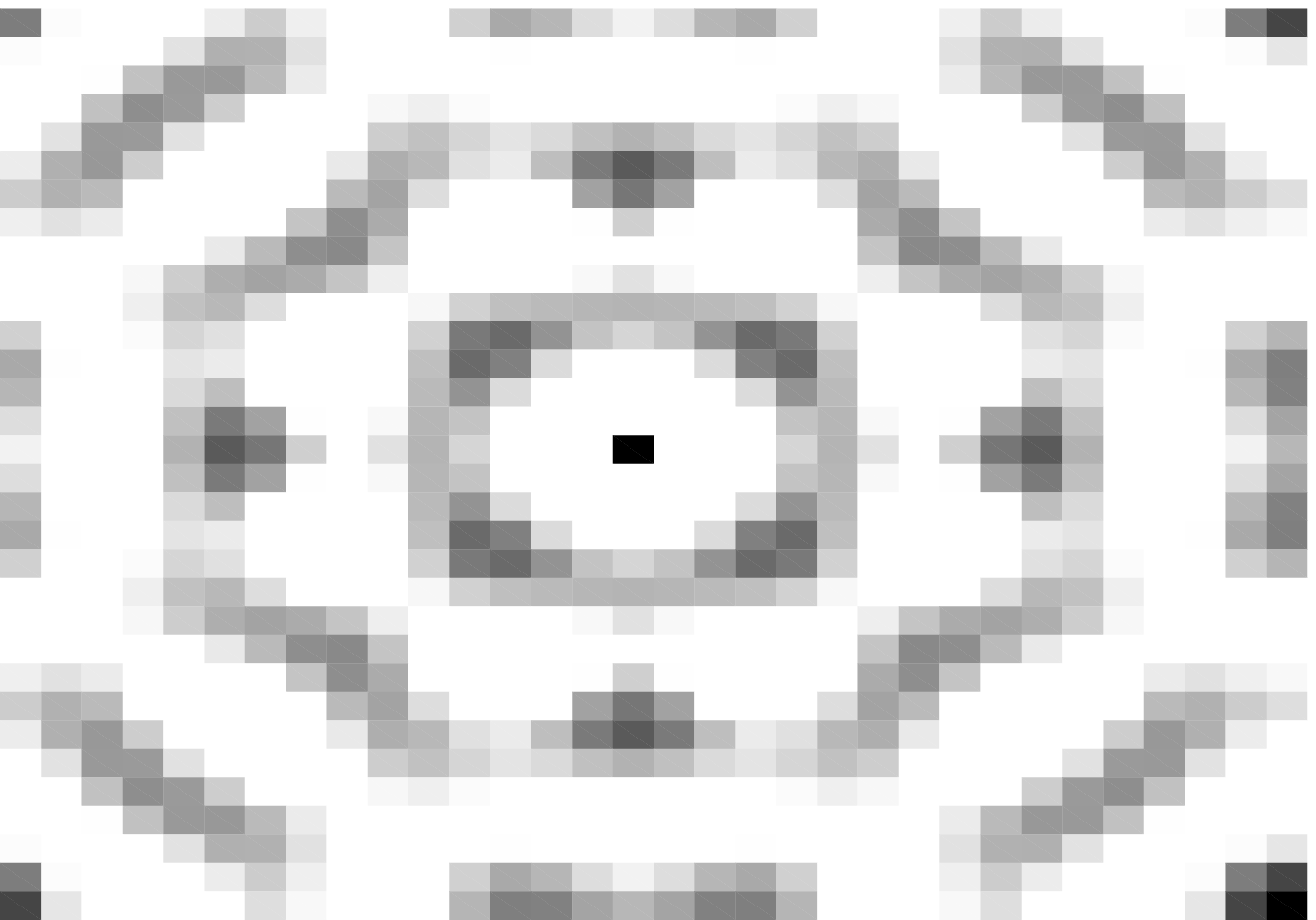}} \hfil
\frame{\includegraphics[width=0.74in]{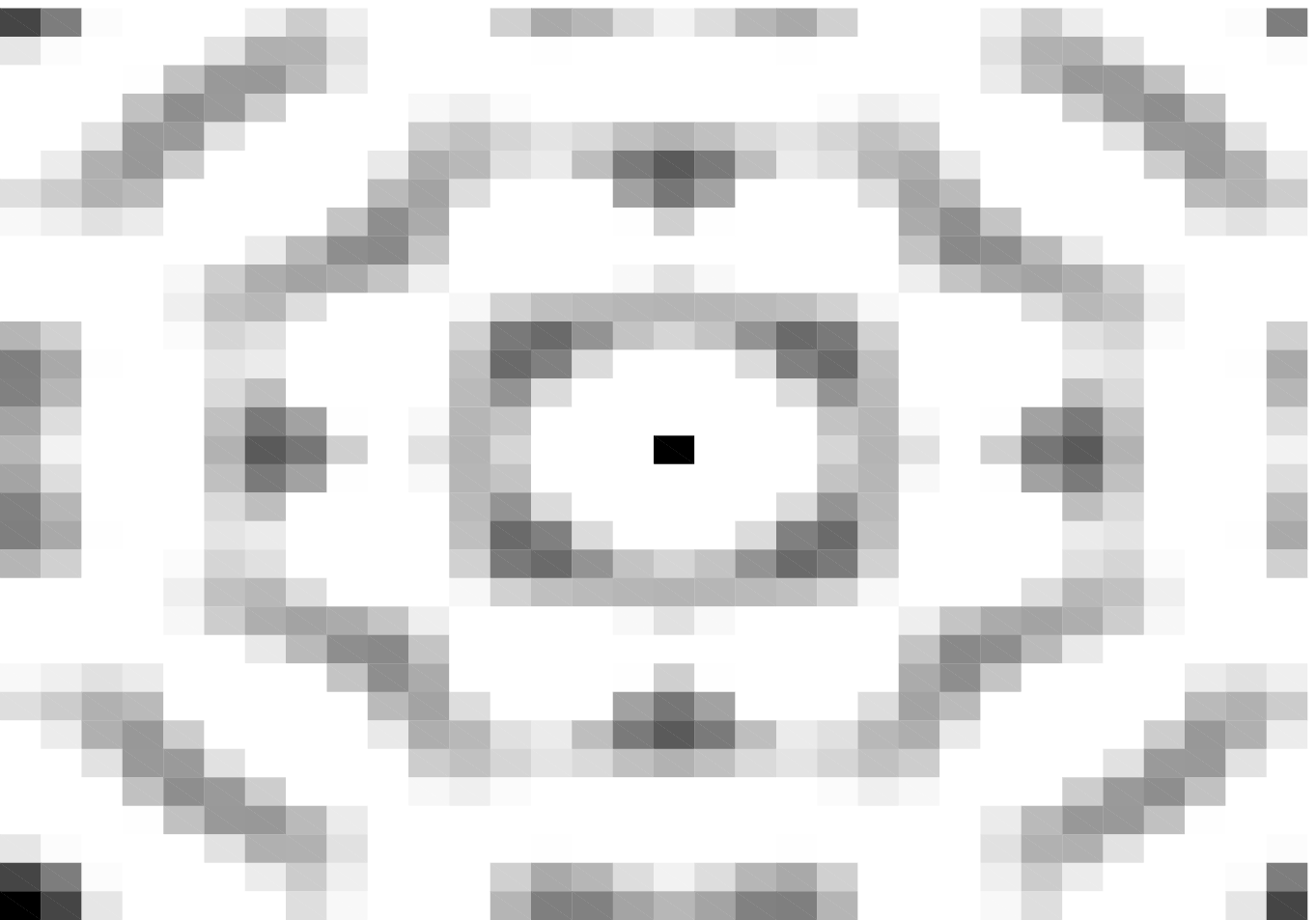}} \hfil
\frame{\includegraphics[width=0.74in]{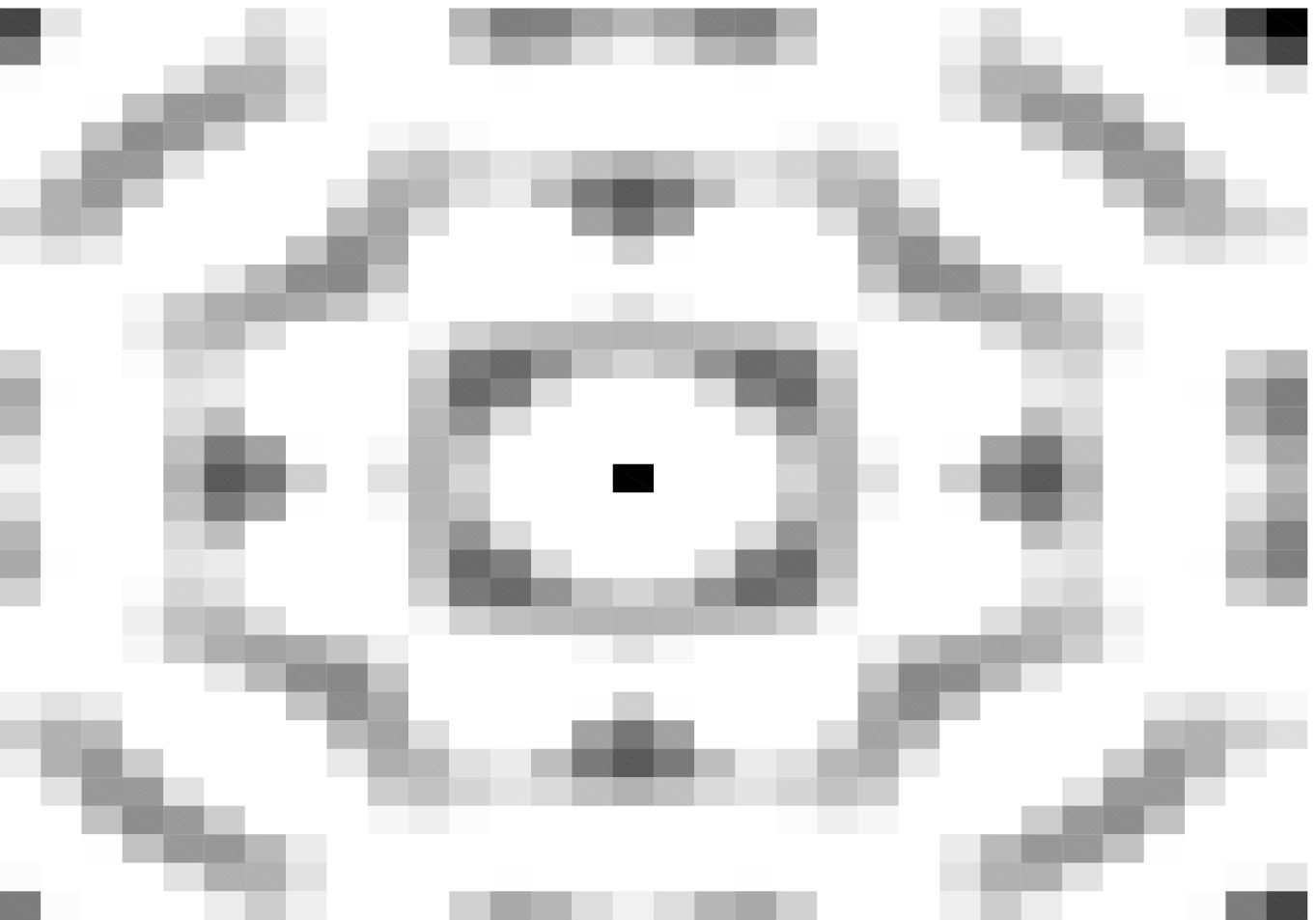}}}}
\centerline{
\subfloat[$u_{rrr}$, $u_{rrb}$, $u_{rbr}$, $u_{rbb}$, $u_{brr}$, $u_{brb}$, $u_{bbr}$ and $u_{bbb}$]{ \label{fig:e2_additive_l3_1}
\frame{\includegraphics[width=0.34in]{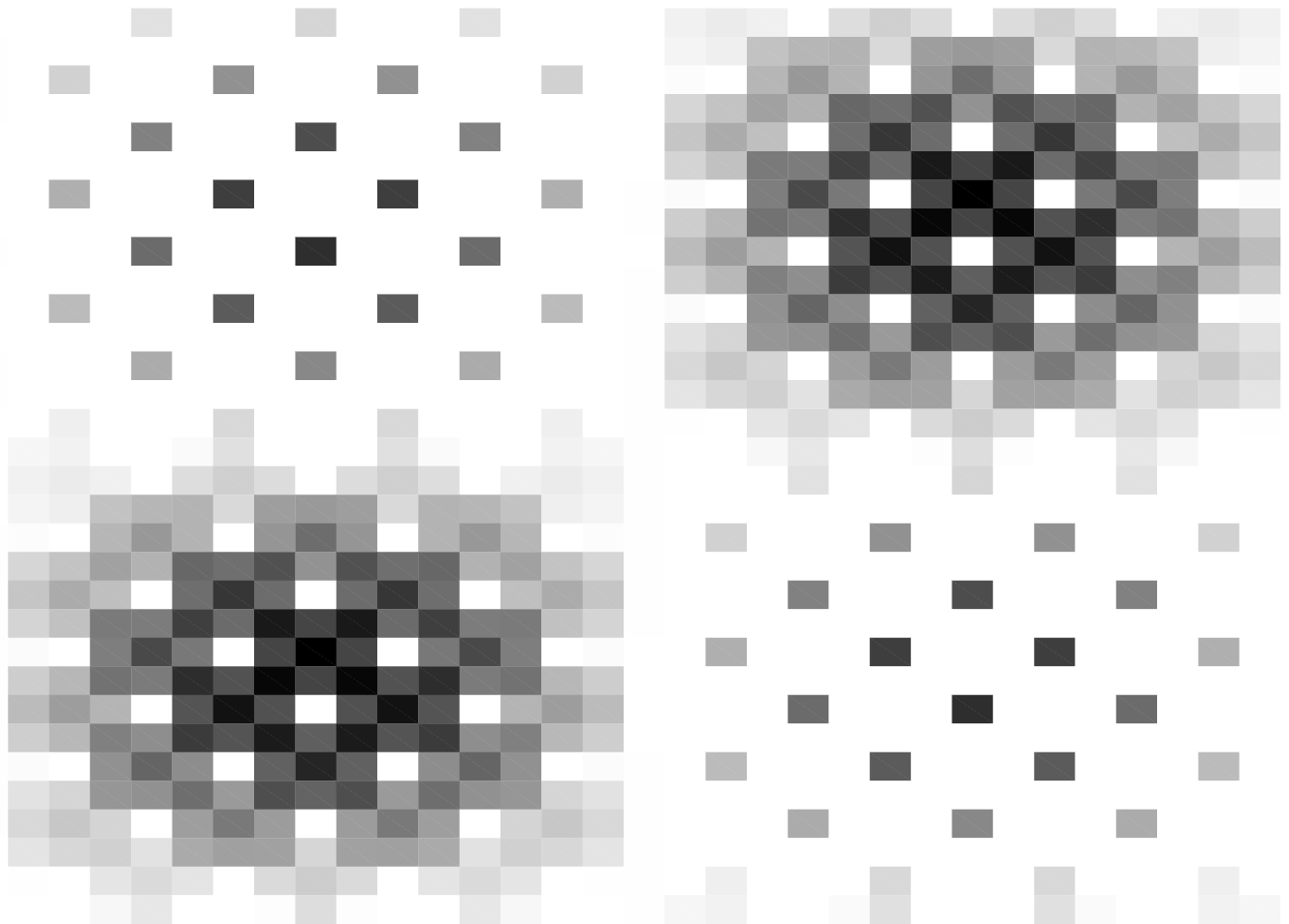}} \hfil
\frame{\includegraphics[width=0.34in]{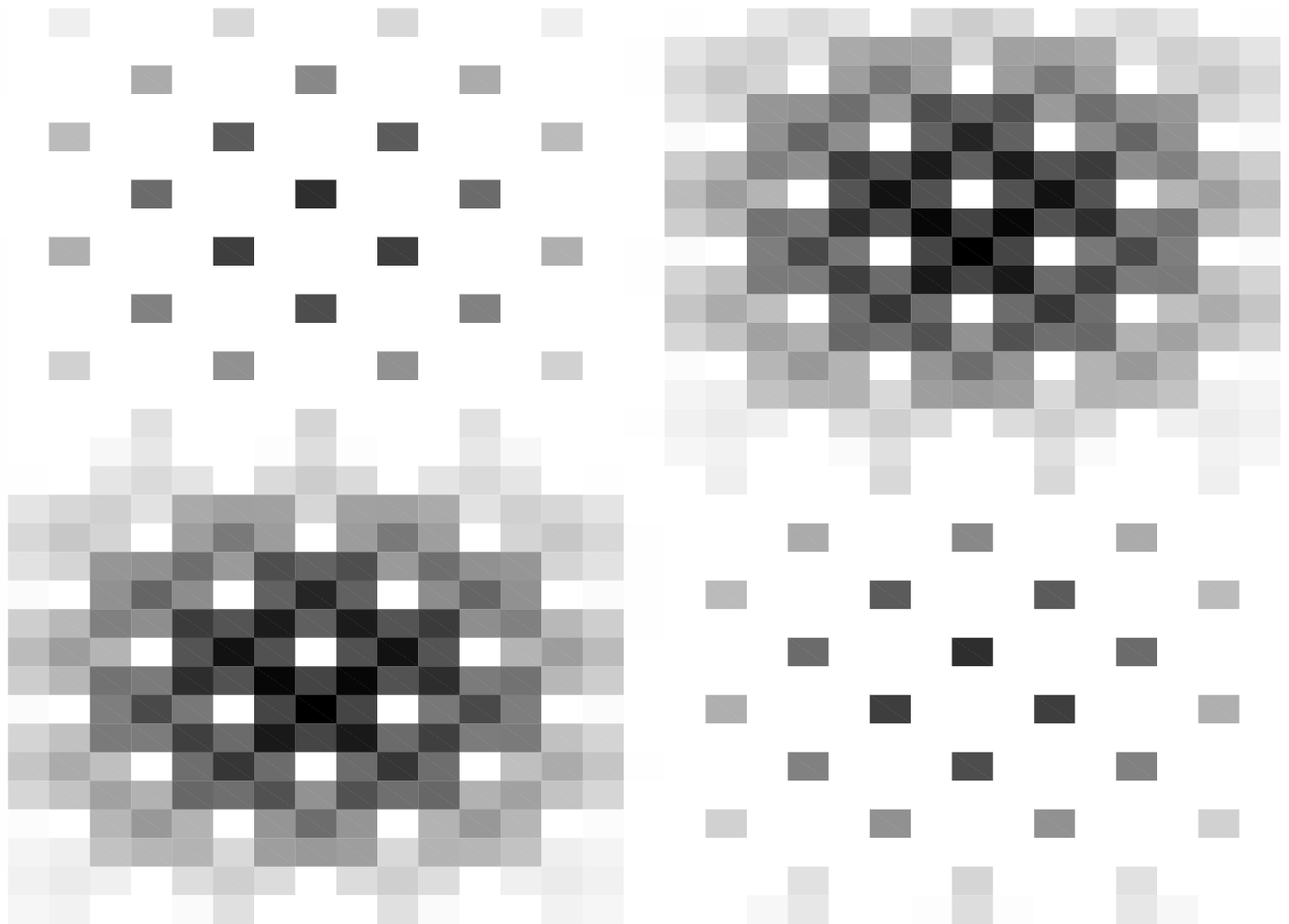}} \hfil
\frame{\includegraphics[width=0.34in]{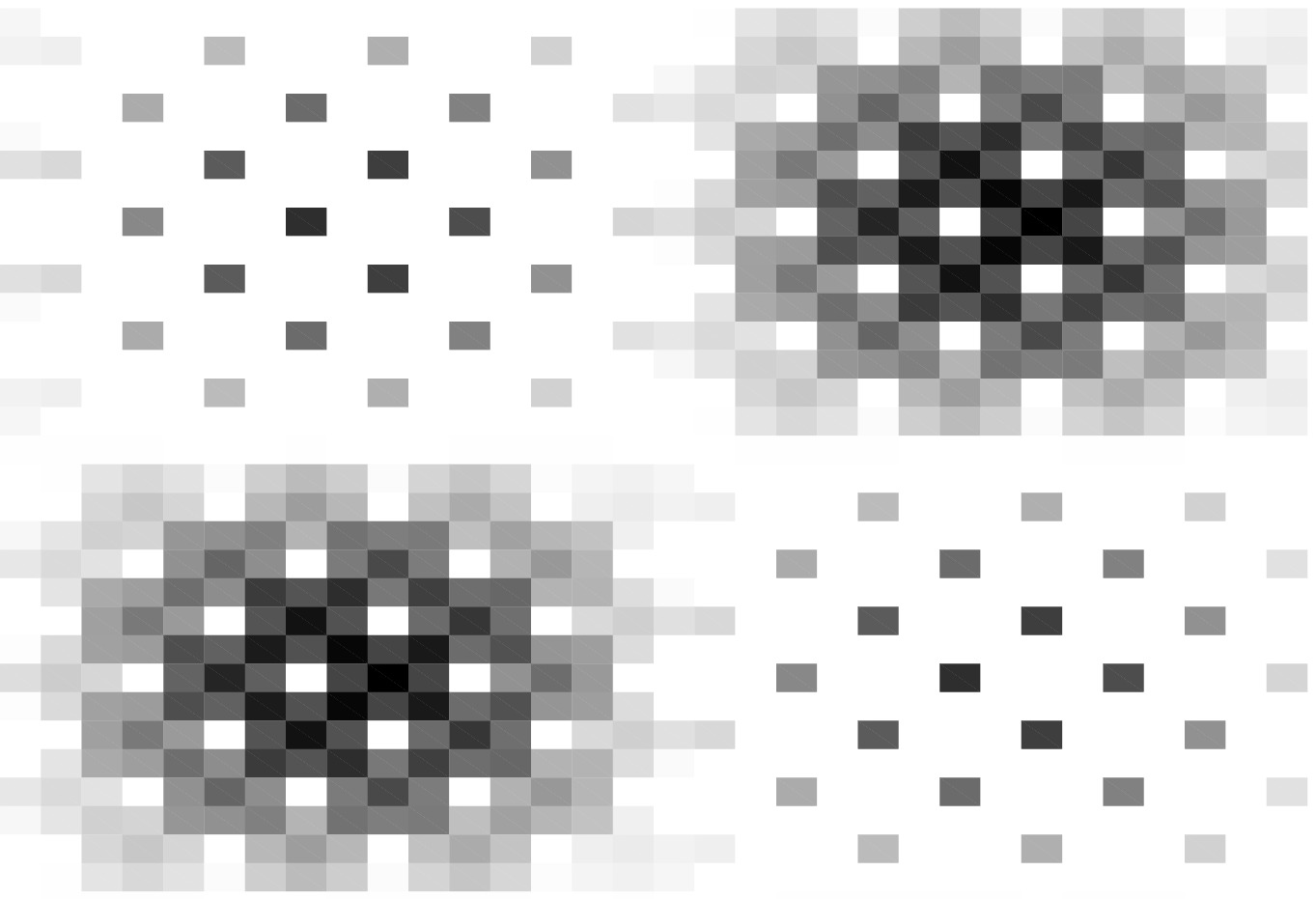}} \hfil
\frame{\includegraphics[width=0.34in]{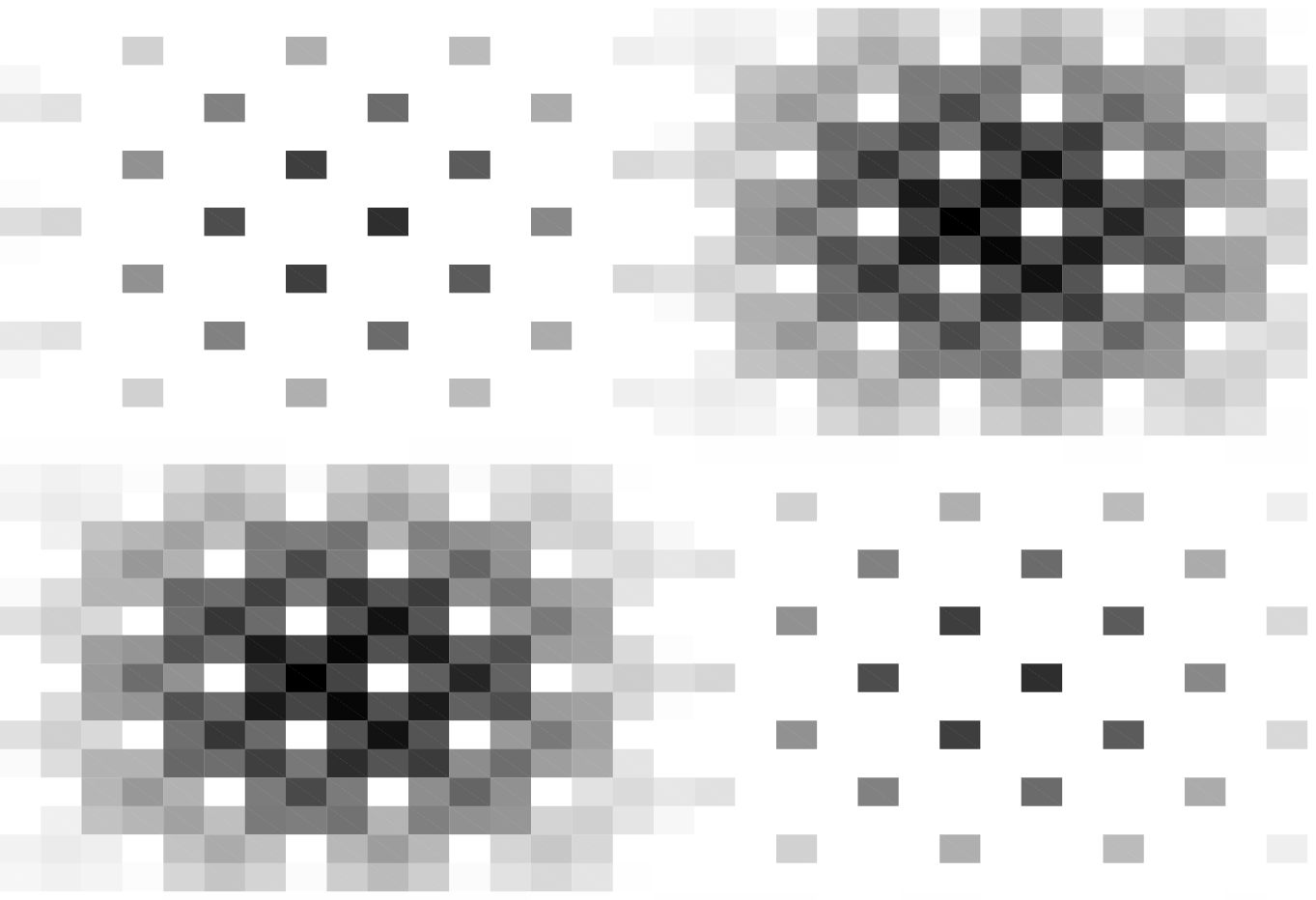}} \hfil
\frame{\includegraphics[width=0.34in]{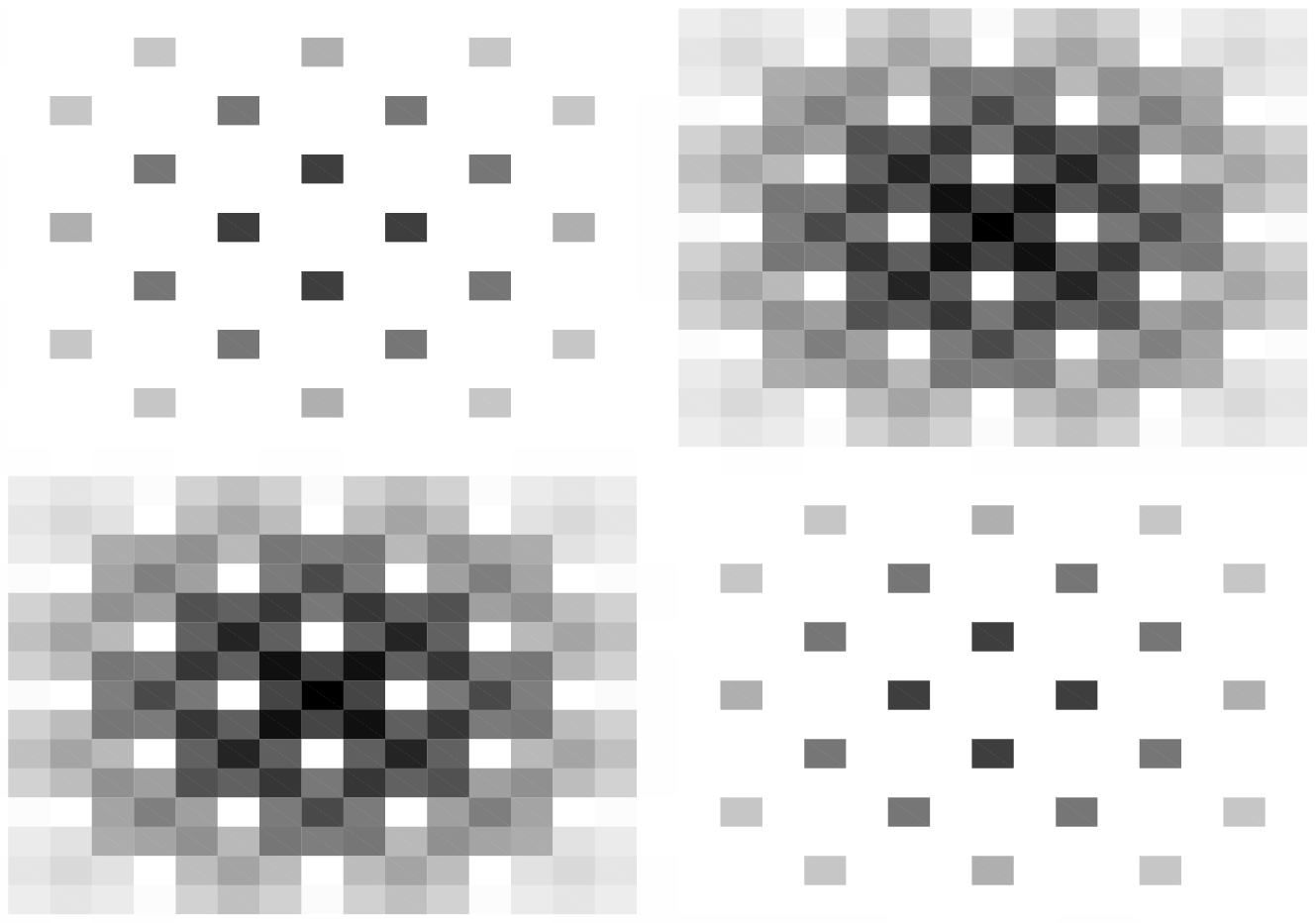}} \hfil
\frame{\includegraphics[width=0.34in]{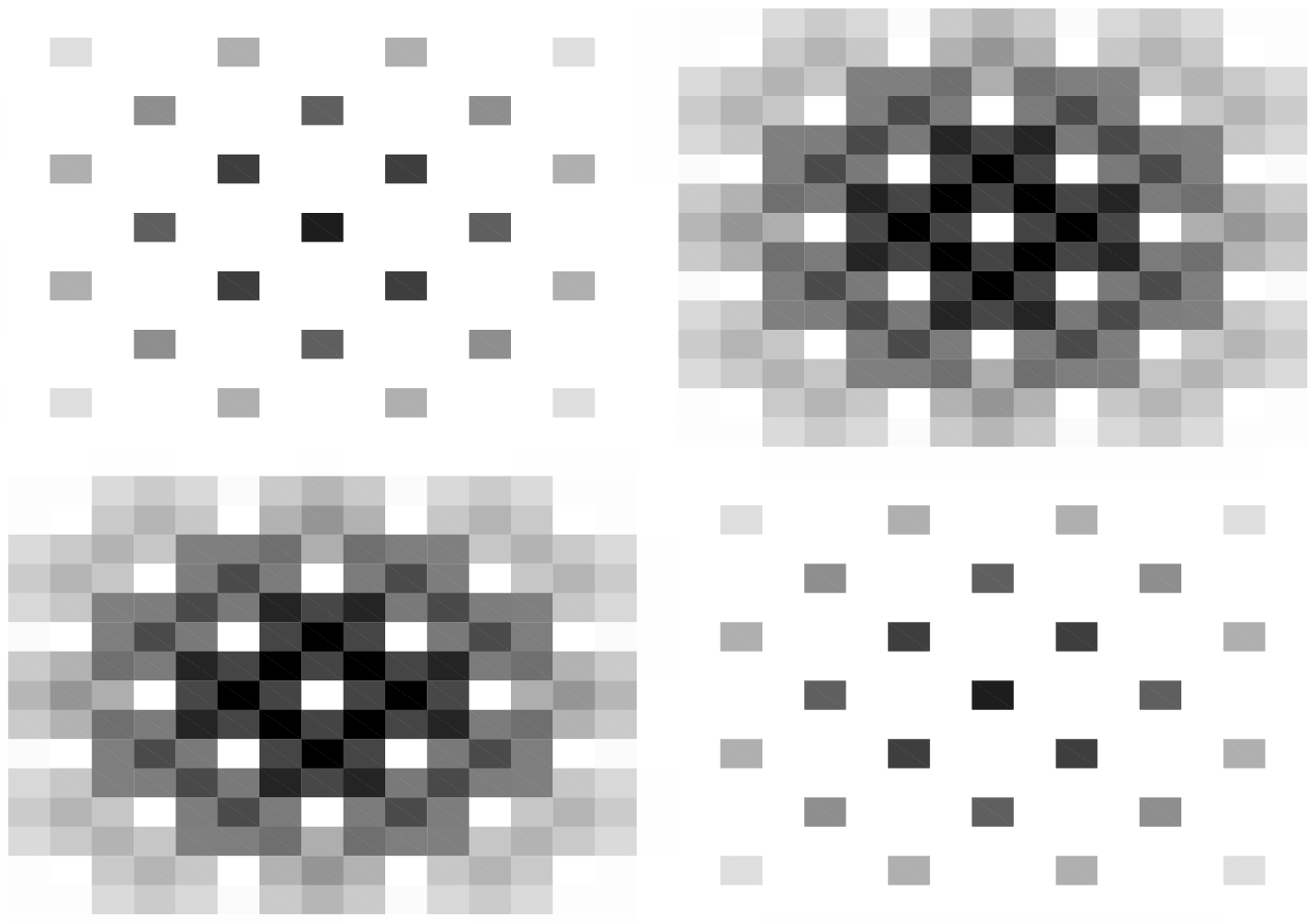}} \hfil
\frame{\includegraphics[width=0.34in]{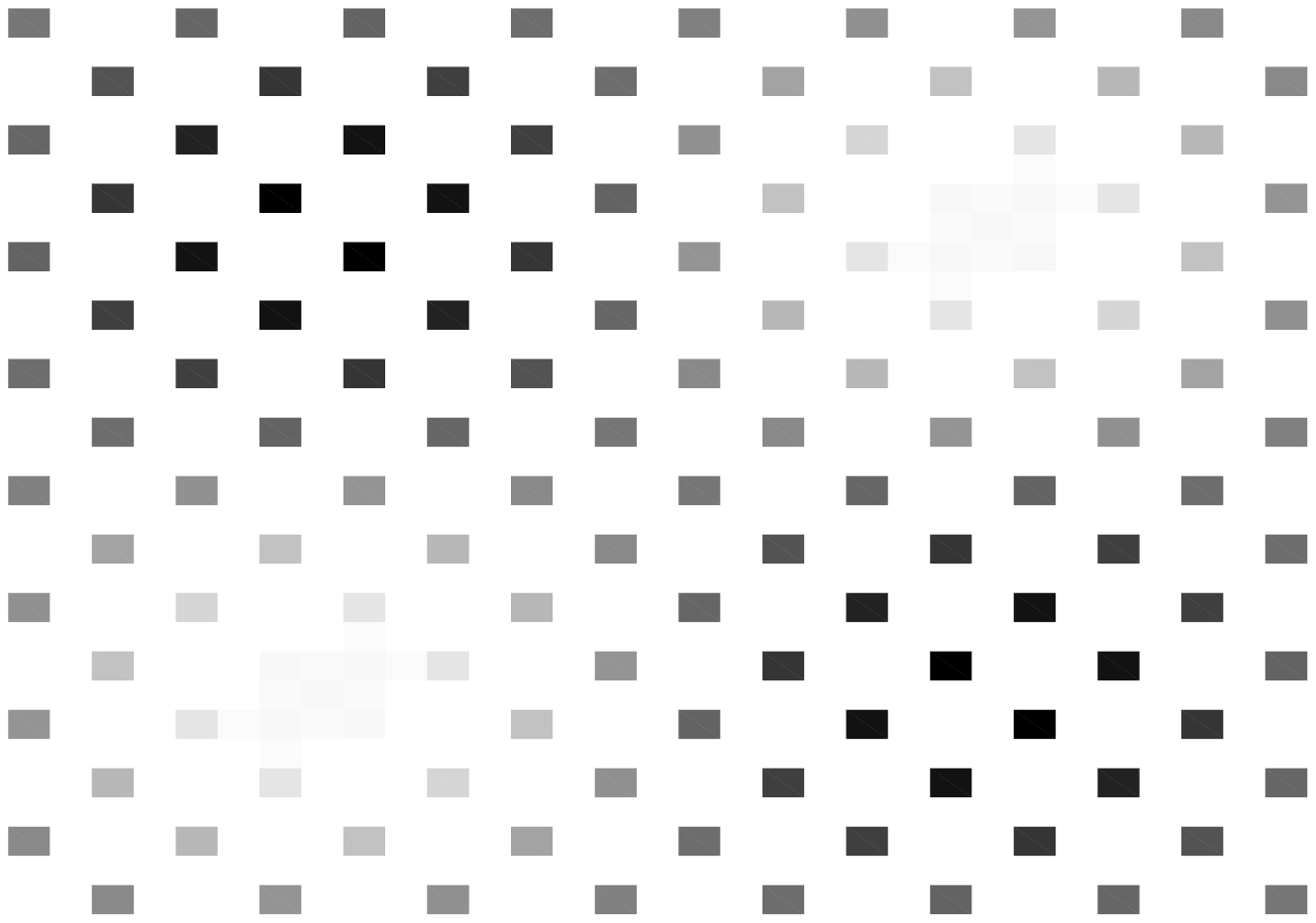}} \hfil
\frame{\includegraphics[width=0.34in]{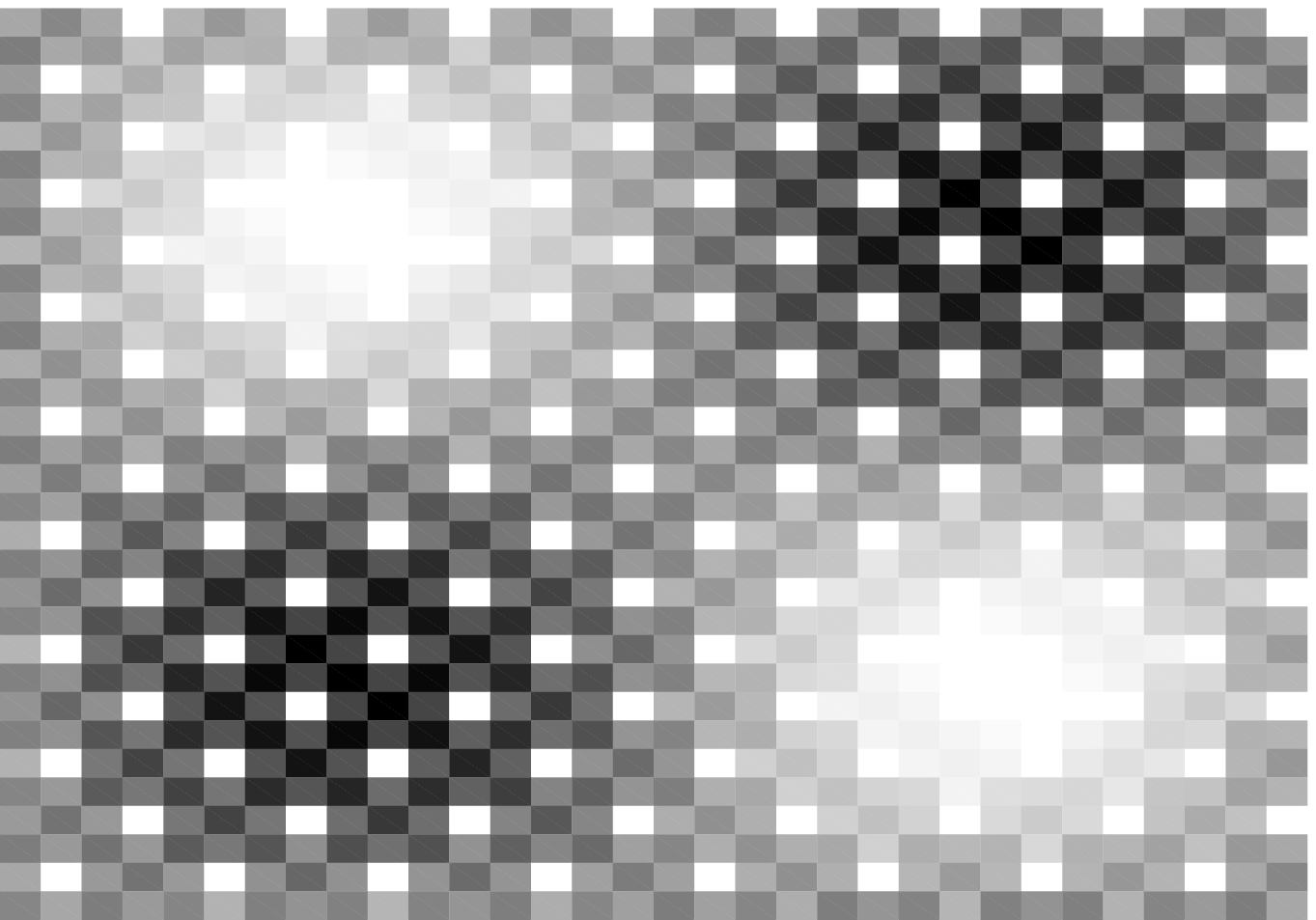}}} \hfil
\subfloat[$u_{rrr}$, $u_{rrb}$, $u_{rbr}$, $u_{rbb}$, $u_{brr}$, $u_{brb}$, $u_{bbr}$ and $u_{bbb}$]{ \label{fig:e2_additive_l3_2}
\frame{\includegraphics[width=0.34in]{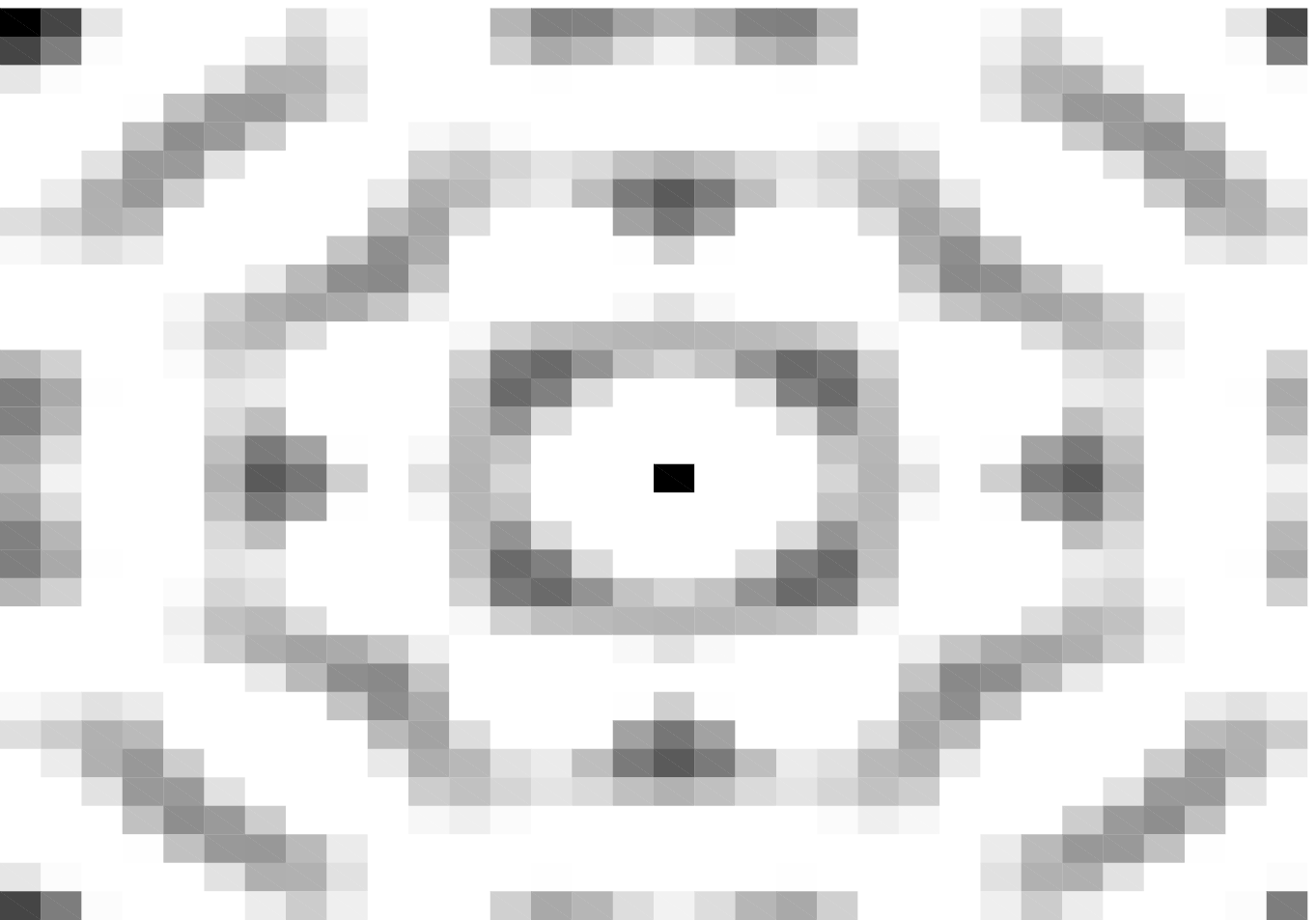}} \hfil
\frame{\includegraphics[width=0.34in]{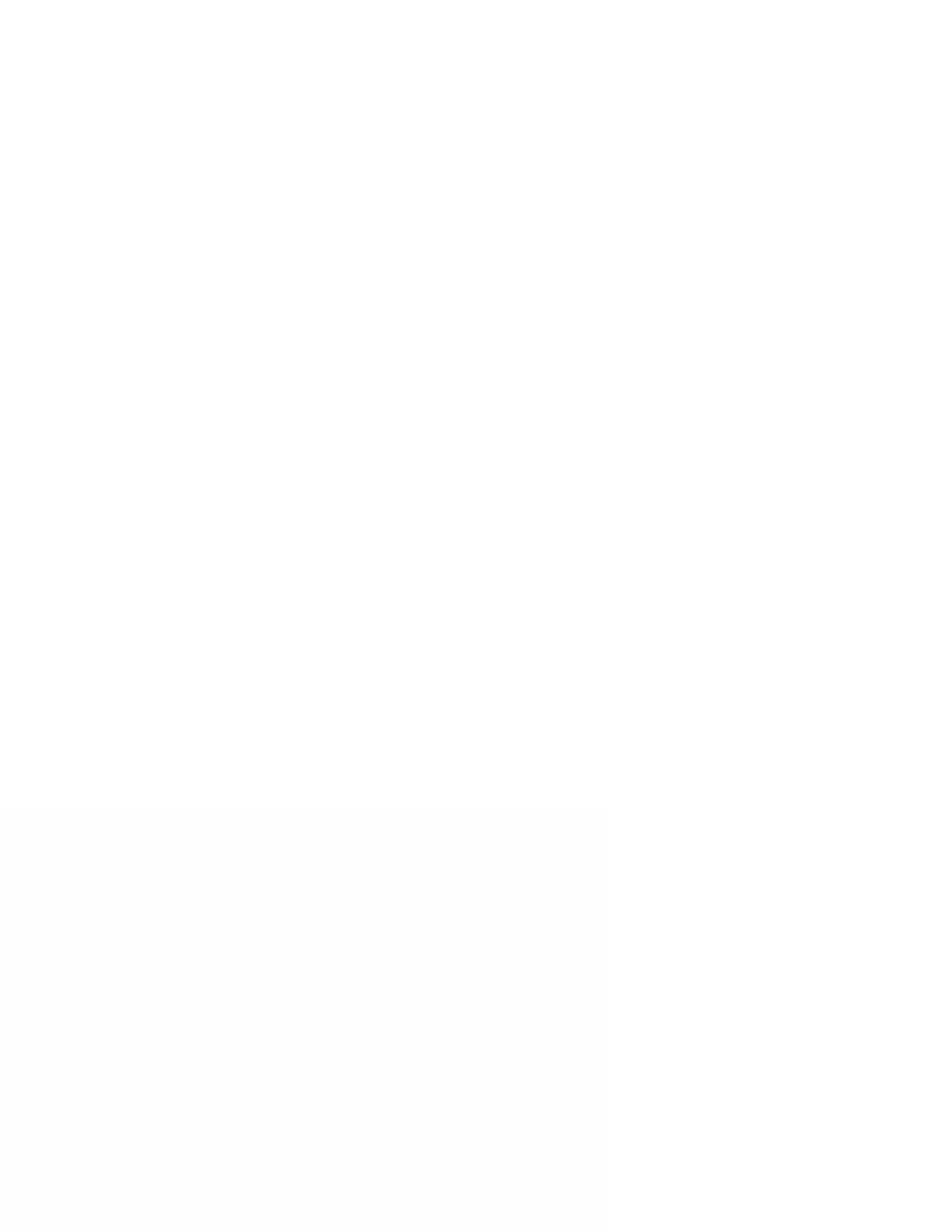}} \hfil
\frame{\includegraphics[width=0.34in]{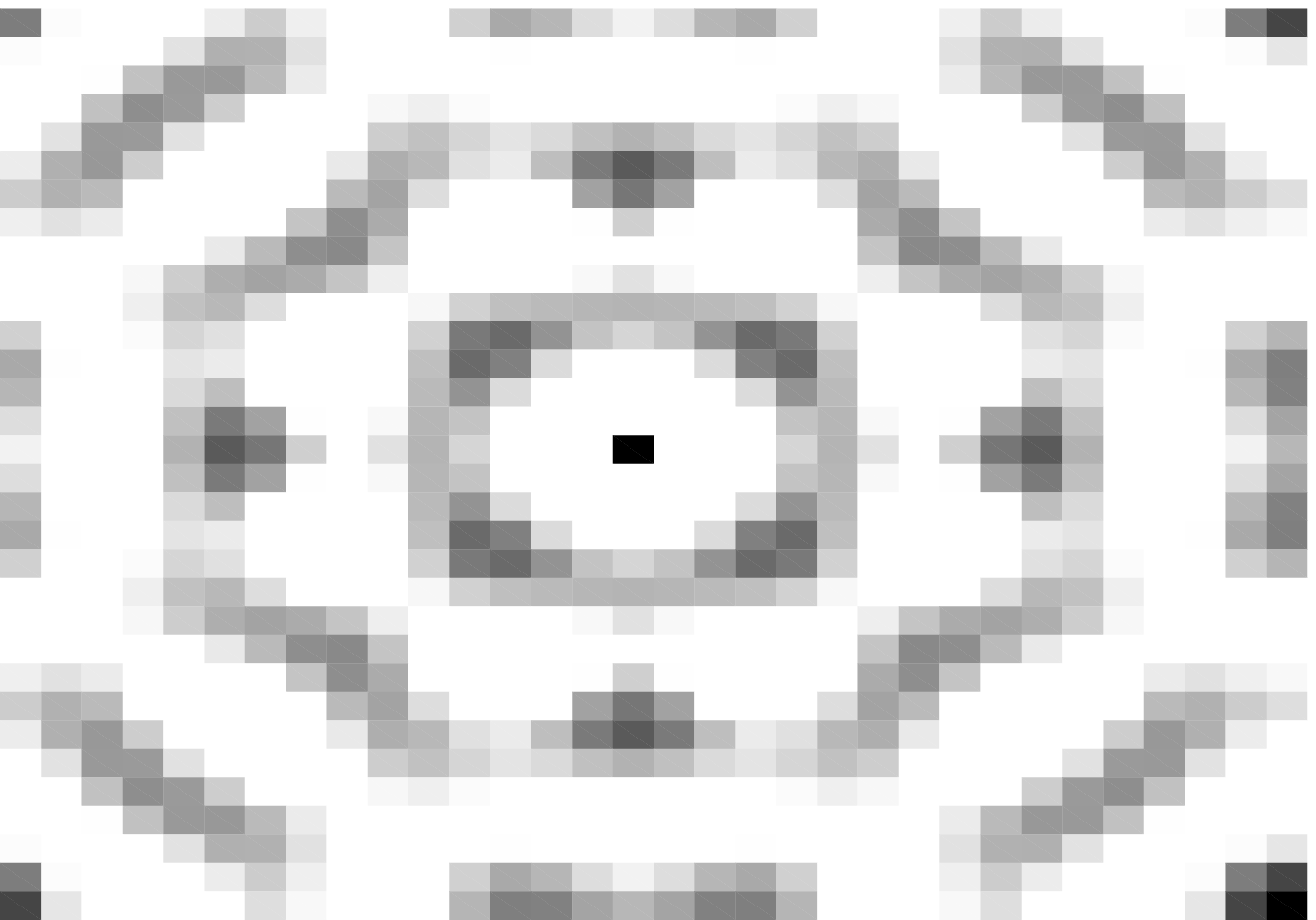}} \hfil
\frame{\includegraphics[width=0.34in]{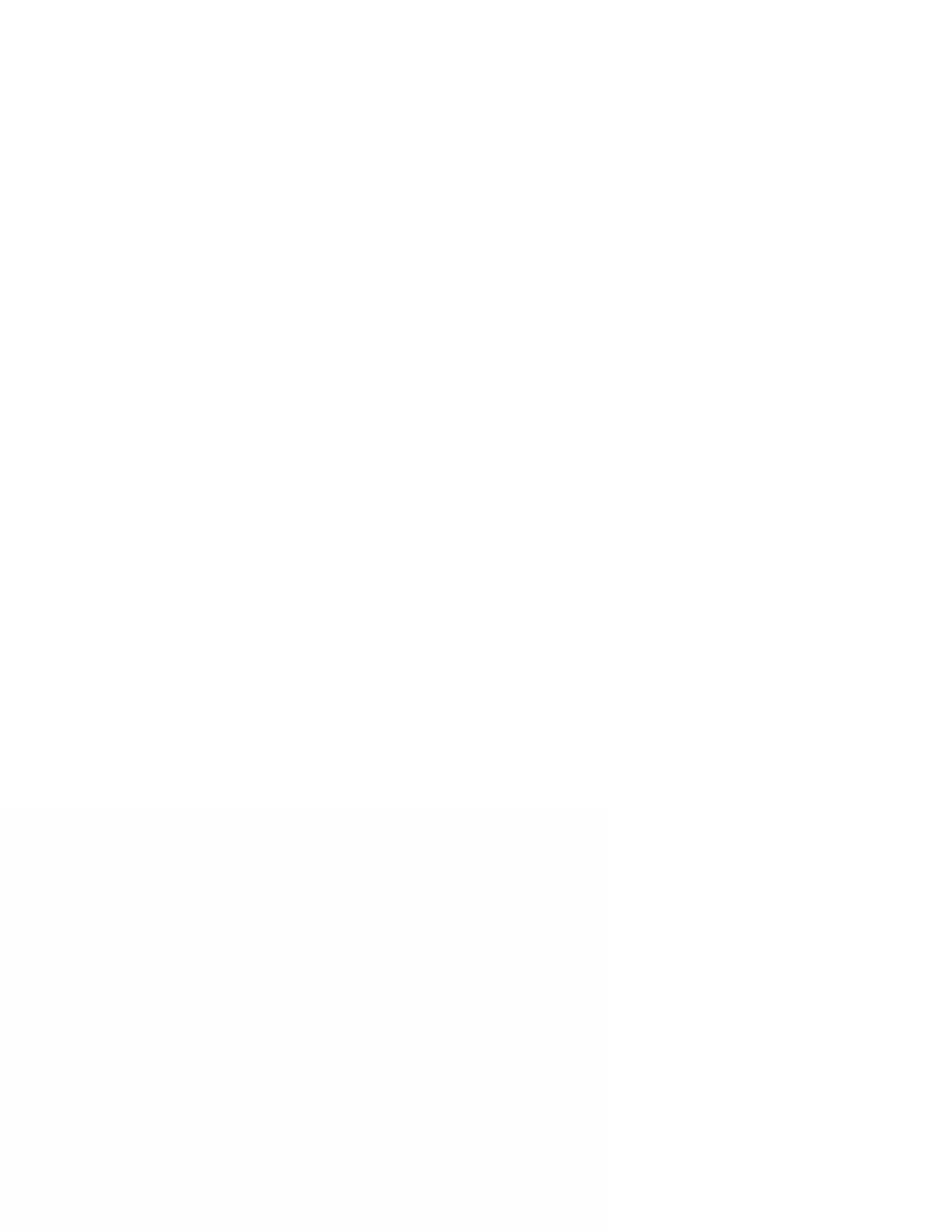}} \hfil
\frame{\includegraphics[width=0.34in]{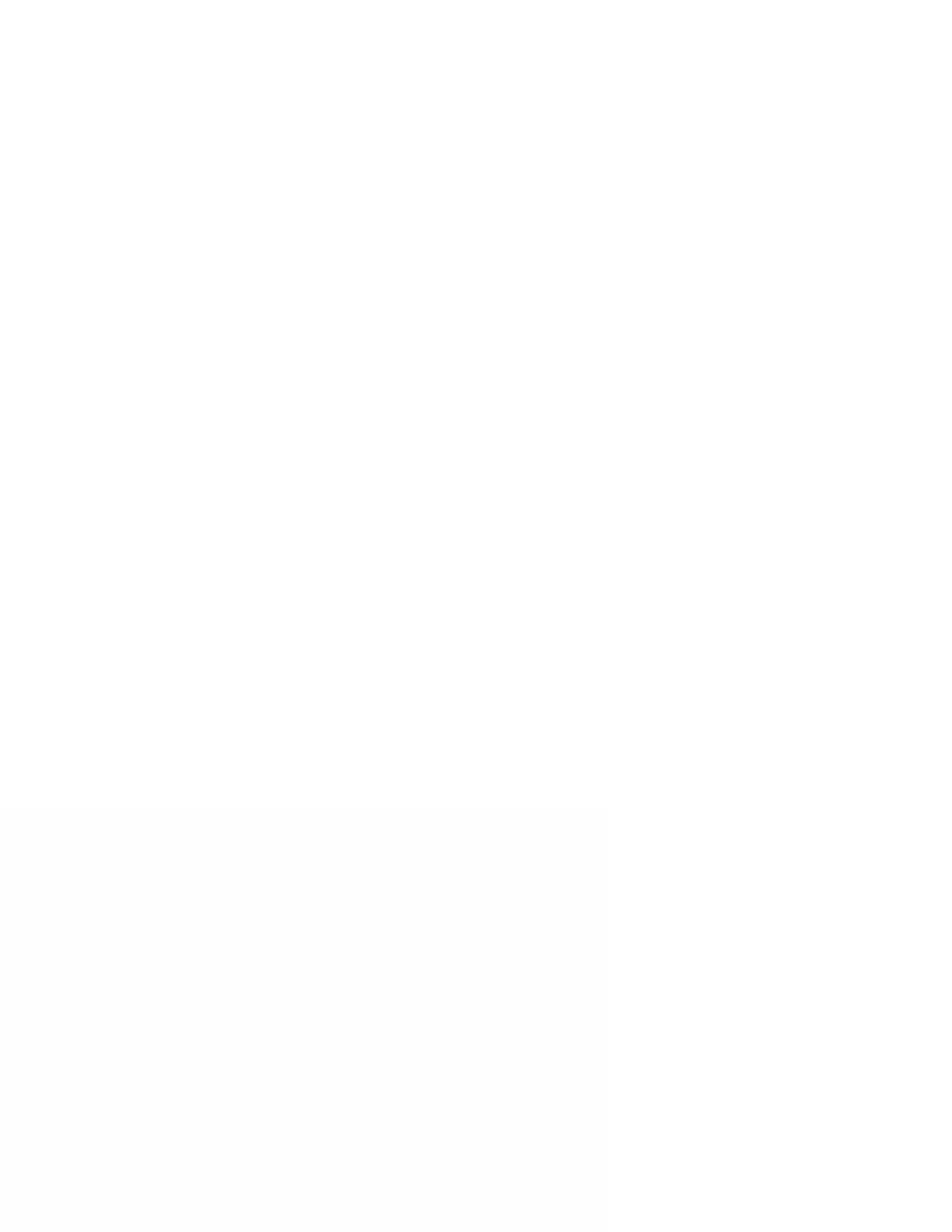}} \hfil
\frame{\includegraphics[width=0.34in]{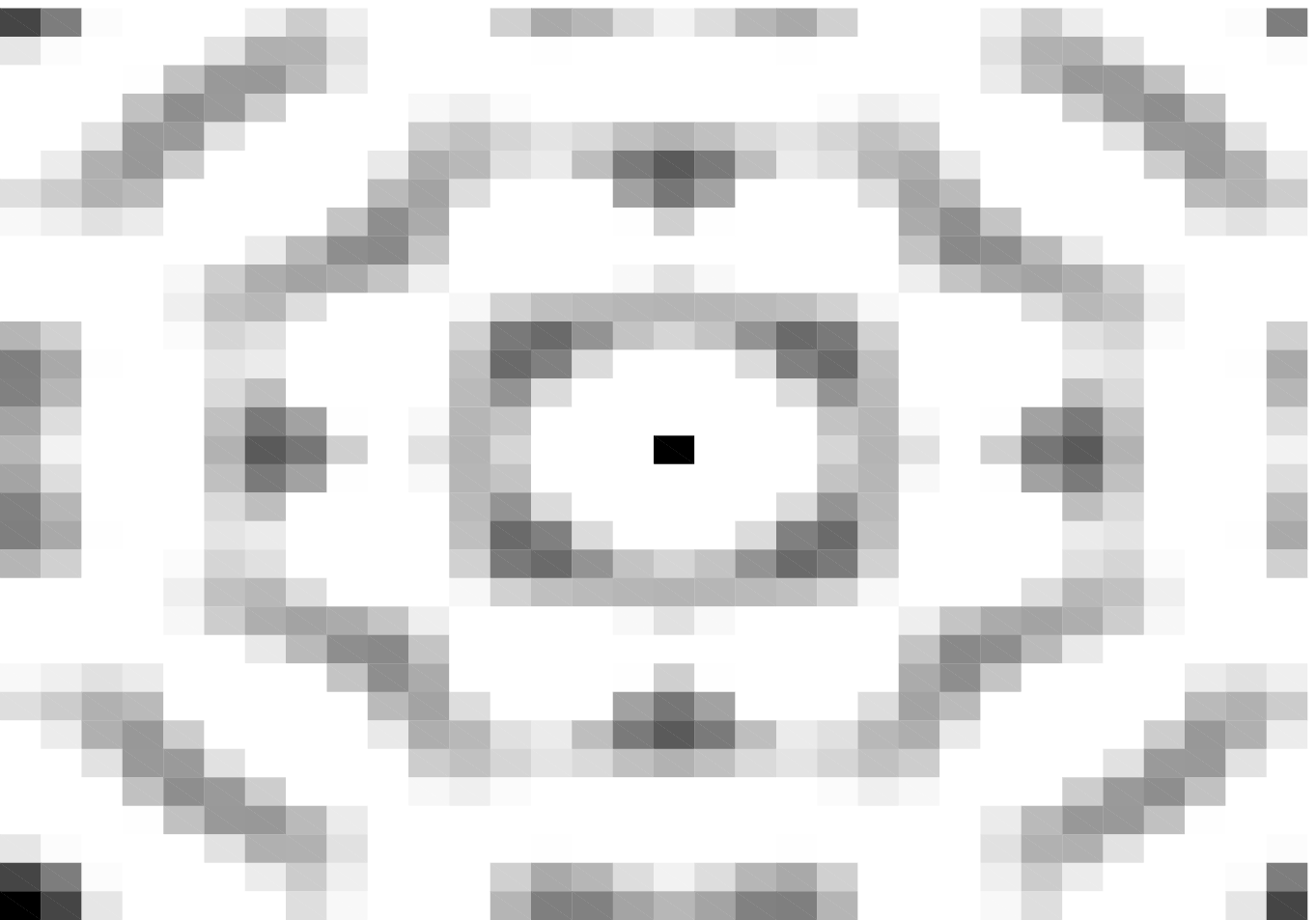}} \hfil
\frame{\includegraphics[width=0.34in]{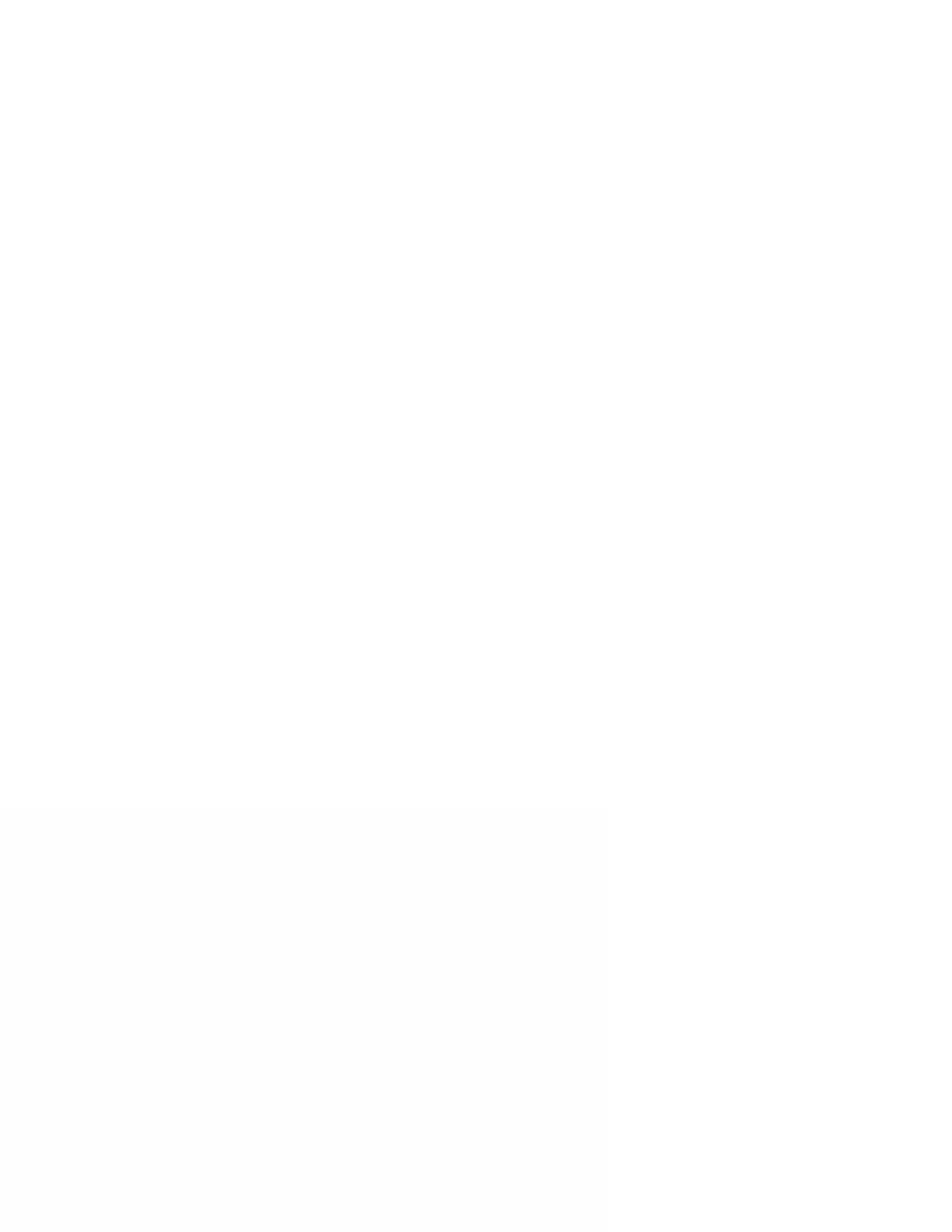}} \hfil
\frame{\includegraphics[width=0.34in]{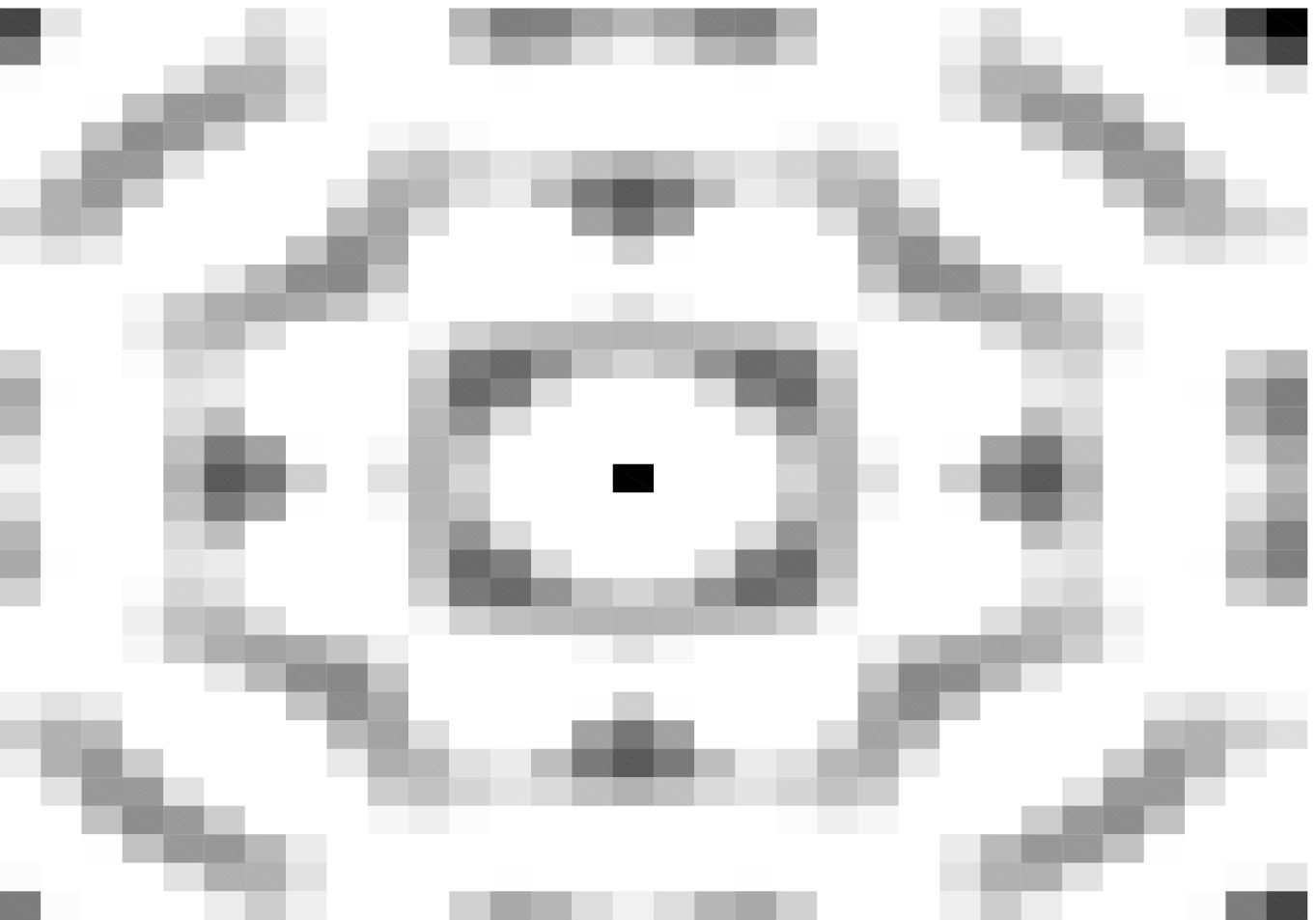}}}}
\caption{Intermediate solutions of Helmholtz's equation using the direct multi--grid algorithm with an additive approach. Considering the two sources in Fig. \ref{fig:e2_source1} and \ref{fig:e2_source2}, the exact solution is given by: $u=u_{r}+u_{b}$ at the first coarse level in Fig. \ref{fig:e2_additive_l1_1} and \ref{fig:e2_additive_l1_2}, respectively; $u=u_{rr}+u_{rb}+u_{br}+u_{bb}$ at the second coarse level in Fig. \ref{fig:e2_additive_l2_1} and \ref{fig:e2_additive_l2_2}, respectively; and $u=u_{rrr}+u_{rrb}+u_{rbr}+u_{rbb}+u_{brr}+u_{brb}+u_{bbr}+u_{bbb}$ at the third coarse level in Fig. \ref{fig:e2_additive_l3_1} and \ref{fig:e2_additive_l3_2}, respectively.} \label{fig:e2_additive_solutions}
\end{figure*}

In Fig. \ref{fig:e2_source1} and \ref{fig:e2_solution1} a solution of Helmholtz's equation is shown when the source vector is the superposition of two frequencies, low and high. In Fig. \ref{fig:e2_multiplicative_v0_1} and \ref{fig:e2_multiplicative_e0_1} the intermediate solutions of a two--grid multiplicative approach are shown. In this case the red coarse system matrix is diagonal and the high--frequency in the forcing function is not seen at the red coarse grid. Therefore, the solution from the red coarse system is trivial and only considers the low--frequency component. The black coarse system adds a correction with the high--frequency and part of the low--frequency component of the solution. In Fig. \ref{fig:e2_additive_l1_1}, \ref{fig:e2_additive_l2_1} and \ref{fig:e2_additive_l3_1} the intermediate solutions of an additive multi--grid approach are shown for three coarse levels. The high--frequency component is not seen in some of the coarse grids.

In Fig. \ref{fig:e2_source2} and \ref{fig:e2_solution2} a solution of Helmholtz's equation is shown for a sparse source vector. In Fig. \ref{fig:e2_multiplicative_v0_2} and \ref{fig:e2_multiplicative_e0_2} the intermediate solutions of a two--grid multiplicative approach are shown. Since the red coarse system matrix is diagonal, the solution of nested iteration is trivial, but not very significant as the correction scheme gives most of the solution. In Fig. \ref{fig:e2_additive_l1_2}, \ref{fig:e2_additive_l2_2} and \ref{fig:e2_additive_l3_2} the solutions at coarse grids are shown by using an additive approach. One of the non--zero values of the source vector appears in each of the coarse grids $\Omega_{rr}$, $\Omega_{rb}$, $\Omega_{br}$ and $\Omega_{bb}$. At the third coarse level the source vector is zero in four of the eight grids. This is a consequence of the down--sampling interpolation and shows the advantage of this approach when the source vector is sparse.

\section{Conclusions}
\label{sec:conclusions}

Numerical methods to solve linear systems of equations were obtained based on the similarities of the full two--grid algorithm and perfect reconstruction filter banks. The two alternatives, multiplicative and additive, correspond to direct Schwartz domain decomposition methods based on a partition of the original domain. The additive approach can be used to parallelize the problem among all the available processors, whereas the multiplicative approach is more efficient in a single processor.

Future research will focus on the application of these algorithms in systems that are not LSI. On one hand, the algorithms are ready to work on systems that are known to have harmonic aliasing patterns but more numerical studies are necessary. And, on the other hand, the most challenging problem is to understand the physical and geometrical implications of harmonic aliasing patterns. This is essential to construct practical methods to check aliasing patterns and be able to use these methods in more challenging problems.

% if have a single appendix:
%\appendix[Proof of the Zonklar Equations]
% or
%\appendix  % for no appendix heading
% do not use \section anymore after \appendix, only \section*
% is possibly needed

% use appendices with more than one appendix
% then use \section to start each appendix
% you must declare a \section before using any
% \subsection or using \label (\appendices by itself
% starts a section numbered zero.)
%

\appendices
%---------------------------------------------------
\section{Proof of Theorem \ref{thm:red_black_surjective}} \label{app:proof_red_black_surjective}
\begin{proof}
The proof that the red harmonic aliasing pattern is equivalent to (\ref{eq:red_surjective}) is due to Navarrete and Coyle \cite{PNavarrete_2008a}. The proof that the black harmonic aliasing pattern is equivalent to (\ref{eq:black_surjective}) follows the same reasoning. Given the partition of eigenvectors $W=\left[W_L W_H\right]$ and $V=\left[V_L V_H\right]$, the black harmonic aliasing pattern can be written as the following set of biorthogonal relationships
\begin{align}
(\bdown V_L)^H (\bdown W_L) & = \tfrac{1}{2} \; I \;, \label{eq:coarse_biorthogonal_LL} \\
(\bdown V_L)^H (\bdown W_H) & = - \tfrac{1}{2} \; I \;, \label{eq:coarse_biorthogonal_LH} \\
(\bdown V_H)^H (\bdown W_L) & = - \tfrac{1}{2} \; I \;\;\mbox{and} \label{eq:coarse_biorthogonal_HL} \\
(\bdown V_H)^H (\bdown W_H) & = \tfrac{1}{2} \; I \;. \label{eq:coarse_biorthogonal_HH}
\end{align}
Since $W$ and $V$ form a biorthogonal basis, $W V^H = W_L V^H_L + W_H V^H_H = I$. Pre-multiplying by $\bdown$ and post-multiplying by $\bup$ gives
\begin{equation} \label{eq:downsample_biorthogonal}
(\bdown W) (\bdown V)^H = (\bdown W_L) (\bdown V_L)^H + (\bdown W_H) (\bdown V_H)^H = I \;.
\end{equation}

First, (\ref{eq:black_surjective}) is assumed. Then, equation (\ref{eq:downsample_biorthogonal}) immediately implies the set of biorthogonal relationships above, and the black harmonic aliasing pattern is fulfilled. Second, the black harmonic aliasing pattern is assumed. Pre--multiplying (\ref{eq:downsample_biorthogonal}) by $(\bdown V_L)^H$ and using equations (\ref{eq:coarse_biorthogonal_LL}) and (\ref{eq:coarse_biorthogonal_LH}) gives $\bdown V_L = - \bdown V_H$. Similarly, post--multiplying (\ref{eq:downsample_biorthogonal}) by $\bdown W_H$ and using equations (\ref{eq:coarse_biorthogonal_LH}) and (\ref{eq:coarse_biorthogonal_HH}) gives $\bdown W_L = - \bdown W_H$. Therefore, the black harmonic aliasing pattern implies (\ref{eq:black_surjective}).
\end{proof}

%---------------------------------------------------
\section{Proof of Lemma \ref{lem:galerkin_inverses}} \label{app:proof_galerkin_inverses}
\begin{proof}
The proof of (\ref{eq:red_galerkin_inverse}) is due to Navarrete and Coyle \cite{PNavarrete_2008a}. The proof of (\ref{eq:black_galerkin_inverse}) follows from
\begin{align}
\bcoarsesys^{-1} & = \left\{ \bdown\brfilter\finesys\bifilter\bup \right\}^{-1} \label{eq:proof_galerkin_1} \\
 & = \left\{ (\bdown W) \brfilterev \Lambda \bifilterev (\bdown V)^H \right\}^{-1} \label{eq:proof_galerkin_2} \\
 & = \left\{ (\bdown W_L) \bdelta (\bdown V_L)^H \right\}^{-1}  \label{eq:proof_galerkin_3} \\
 & = 4\; (\bdown W_L) \bdelta^{-1} (\bdown V_L)^H \;, \label{eq:proof_galerkin_4}
\end{align}
where (\ref{eq:igrid_decomposition}) is used in (\ref{eq:proof_galerkin_1}), the eigen--decompositions of filters is used in (\ref{eq:proof_galerkin_2}), Theorem \ref{thm:red_black_surjective} is used in (\ref{eq:proof_galerkin_3}) and the biorthogonal relationships (\ref{eq:coarse_biorthogonal_LL}) to (\ref{eq:coarse_biorthogonal_HH}) are used in (\ref{eq:proof_galerkin_4}).
\end{proof}

%---------------------------------------------------
\section{Proof of Theorem \ref{thm:cgc_decomposition}} \label{app:proof_cgc_decomposition}
\begin{proof}
The proof of (\ref{eq:rcgc_decomposition}) is due to Navarrete and Coyle \cite{PNavarrete_2008a}. The proof of (\ref{eq:bcgc_decomposition}) follows from
\begin{align}
\bcgc & = I - \bifilter\bup\bcoarsesys^{-1}\bdown\brfilter \finesys \label{eq:proof_bcgc_1} \\
 & = I - 4\; W \bifilterev ( V^H \bup \bdown W_L ) \bdelta^{-1} ( V_L^H \bup \bdown W ) \brfilterev \Lambda V^H \label{eq:proof_bcgc_2} \\
 & = W V^H - 4\; W \bifilterev \left(
\tfrac{1}{2} \left[ \begin{smallmatrix}
I \\
-I
\end{smallmatrix} \right]
\right)
\bdelta^{-1}
\left(
\tfrac{1}{2} \left[ \begin{smallmatrix}
I & -I
\end{smallmatrix} \right]
\right) \brfilterev \Lambda V^H \label{eq:proof_bcgc_3} \\
& = W
\left[\begin{smallmatrix}
I - \bifilterevL \bdelta^{-1} \brfilterevL \Lambda_L
& \bifilterevL \bdelta^{-1} \brfilterevH \Lambda_H \\
\bifilterevH \bdelta^{-1} \brfilterevL \Lambda_{L}
& I - \bifilterevH \bdelta^{-1} \brfilterevH \Lambda_{H} \\
\end{smallmatrix}\right]
V^H \;. \label{eq:proof_bcgc_4} 
\end{align}
where (\ref{eq:igrid_decomposition}) is used in (\ref{eq:proof_bcgc_1}), the eigen--decomposition of filters and Lemma \ref{lem:galerkin_inverses} are used in (\ref{eq:proof_bcgc_2}), definition \ref{def:rbHAP} is used in (\ref{eq:proof_bcgc_3}) and the partition of eigenvalues is used in (\ref{eq:proof_bcgc_4}).
\end{proof}

% use section* for acknowledgement
%\section*{Acknowledgment}
%The authors would like to thank...

% Can use something like this to put references on a page
% by themselves when using endfloat and the captionsoff option.
\ifCLASSOPTIONcaptionsoff
  \newpage
\fi

% trigger a \newpage just before the given reference
% number - used to balance the columns on the last page
% adjust value as needed - may need to be readjusted if
% the document is modified later
%\IEEEtriggeratref{8}
% The "triggered" command can be changed if desired:
%\IEEEtriggercmd{\enlargethispage{-5in}}

% references section

% can use a bibliography generated by BibTeX as a .bbl file
% BibTeX documentation can be easily obtained at:
% http://www.ctan.org/tex-archive/biblio/bibtex/contrib/doc/
% The IEEEtran BibTeX style support page is at:
% http://www.michaelshell.org/tex/ieeetran/bibtex/
\bibliographystyle{IEEEtran}
% argument is your BibTeX string definitions and bibliography database(s)
\bibliography{IEEEabrv,manuscript}
%
% <OR> manually copy in the resultant .bbl file
% set second argument of \begin to the number of references
% (used to reserve space for the reference number labels box)
%\begin{thebibliography}{1}
%
%\bibitem{IEEEhowto:kopka}
%H.~Kopka and P.~W. Daly, \emph{A Guide to \LaTeX}, 3rd~ed.\hskip 1em plus
%  0.5em minus 0.4em\relax Harlow, England: Addison-Wesley, 1999.
%
%\end{thebibliography}

% biography section
% 
% If you have an EPS/PDF photo (graphicx package needed) extra braces are
% needed around the contents of the optional argument to biography to prevent
% the LaTeX parser from getting confused when it sees the complicated
% \includegraphics command within an optional argument. (You could create
% your own custom macro containing the \includegraphics command to make things
% simpler here.)
%\begin{biography}[{\includegraphics[width=1in,height=1.25in,clip,keepaspectratio]{mshell}}]{Michael Shell}
% or if you just want to reserve a space for a photo:

%\begin{IEEEbiography}{Michael Shell}
%Biography text here.
%\end{IEEEbiography}

% if you will not have a photo at all:
\begin{IEEEbiographynophoto}{Pablo Navarrete}
was born in Santiago, Chile. He received the B.Sc. in Physics (2001), B.Sc. in Electrical Engineering (2001) and the Electrical Engineer Degree (2002), from Universidad de Chile at Santiago. He received the Ph.D. degree in Electrical Engineering from Purdue University at West Lafayette, in 2008. He worked as a research intern in CIMNE at Technical University of Catalonia in 2006, and as a visitor student research collaborator at Princeton University at Princeton, NJ, in 2006--2007. He is currently working in the Department of Electrical Engineering at Universidad de Chile. His major interests are the application of signal processing tools in problems of numerical analysis and especially those based on multi--level techniques. His research interests also include wavelet analysis, statistical signal processing and stochastic analysis.
\end{IEEEbiographynophoto}

% insert where needed to balance the two columns on the last page with
% biographies
%\newpage

%\begin{IEEEbiographynophoto}{Jane Doe}
%Biography text here.
%\end{IEEEbiographynophoto}

% You can push biographies down or up by placing
% a \vfill before or after them. The appropriate
% use of \vfill depends on what kind of text is
% on the last page and whether or not the columns
% are being equalized.

%\vfill

% Can be used to pull up biographies so that the bottom of the last one
% is flush with the other column.
%\enlargethispage{-5in}

% that's all folks
\end{document}